\documentclass[12pt]{article}
\usepackage{amsmath, amsfonts, amssymb,amsthm, amscd}

\usepackage{pictexwd,dcpic}
\usepackage{graphicx}
\usepackage{epsf}
\usepackage{epsfig}
\usepackage{rlepsf}
\usepackage{psfrag}
\usepackage{color}
\usepackage{hyperref}
\usepackage{makeidx}

\theoremstyle{definition}

\theoremstyle{remark}

\newtheoremstyle{thm}
  {12pt}
  {12pt}
  {\itshape}
  {\parindent}
  {\scshape}
  {.}
  {5pt}
  {}

\theoremstyle{thm}
\newtheorem*{T6.17}{Theorem \ref{gradientthm1}}
\newtheorem{thm}{Theorem}[section]

\newtheoremstyle{prop}
  {12pt}
  {12pt}
  {\itshape}
  {\parindent}
  {\scshape}
  {.}
  {5pt}
  {}

\theoremstyle{prop}
\newtheorem{prop}[thm]{Proposition}

\newtheoremstyle{lem}
  {12pt}
  {12pt}
  {\itshape}
  {\parindent}
  {\scshape}
  {.}
  {5pt}
  {}

\theoremstyle{lem}
\newtheorem*{L2.6}{Lemma \ref{lem2.4}}
\newtheorem{lem}[thm]{Lemma}

\newtheoremstyle{defn}
  {12pt}
  {12pt}
  {\itshape}
  {\parindent}
  {\scshape}
  {.}
  {5pt}
  {}

\theoremstyle{defn}
\newtheorem{defn}[thm]{Definition}

\newtheoremstyle{examp}
  {12pt}
  {12pt}
  {}
   {\parindent}
  {\scshape}
  {.}
  {5pt}
  {}

\theoremstyle{examp}
\newtheorem{examp}[thm]{Example}

\newtheoremstyle{cor}
  {12pt}
  {12pt}
  {\itshape}
  {\parindent}
  {\scshape}
  {.}
  {5pt}
  {}

\theoremstyle{cor}

\newtheoremstyle{recipe}
  {12pt}
  {12pt}
  {\itshape}
   {\parindent}
  {\scshape}
  {.}
  {5pt}
  {}

\theoremstyle{recipe}

\newtheoremstyle{rem}
  {12pt}
  {12pt}
  {}
   {\parindent}
  {\scshape}
  {.}
  {5pt}
  {}

\theoremstyle{rem}
\newtheorem{rem}[thm]{Remark}

\newcommand{\bp}{\begin{proof}}
\newcommand{\ep}{\end{proof}}

\newcommand{\norm}[1]{\left\Vert#1\right\Vert}

\newcommand{\ssc}{\text{sc}}

\renewcommand{\epsilon}{\varepsilon}

\newcommand{\what}{\widehat}

\newcommand{\cl}{\operatorname{cl}}

\newcommand{\ind}{\operatorname{Ind}}

\newcommand{\supp}{\operatorname{supp}}

\providecommand{\ker}[1]{$\text{ker}\ {#1}$}

\def\norm#1{\mathopen\|#1\mathclose\|}

{\catcode`"=12 \gdef\hex{"}}

\mathchardef\laplace=\hex0001

\mathchardef\nabla=\hex0272

\makeatletter

\def\@@dalembert#1#2{\setbox0\hbox{$#1\mathrm I$}

  \vrule height\ht0 depth\z@ width.04\ht0

  \rlap{\vrule height\ht0 depth-.96\ht0 width.8\ht0}

  \vrule height.1\ht0 depth\z@ width.8\ht0

  \vrule height\ht0 depth\z@ width.1\ht0 }

\def\dalembert{\mathbin{\mathpalette\@@dalembert{}}\,}

\makeatother

\begin{document}

\title{Polyfolds and Fredholm Theory}

\author{ Helmut H. W. Hofer\\  
Institute for Advanced Study
}

\maketitle
\tableofcontents

\section{Introduction}\label{SEC1}

In this paper we discuss the generalized Fredholm theory in polyfolds.  The initial version of the paper, \cite{H2008},  was written in 2008 on the occasion
of a lecture given at the Clay Mathematical Institute (CMI), and described the theory 
as developed in \cite{HWZ1,HWZ2,HWZ3}. When Nick Woodhouse became the CMI president he found a folder of manuscripts from past
CMI  events and came up with the idea of having them published. Six years is a long time.   Since then the theory considerably advanced and a comprehensive discussion of its more  abstract part  is contained in \cite{HWZ7,HWZ8}.  Currently nontrivial applications are being developed, most notably
those to symplectic field theory (SFT), see \cite{FHWZ} and \cite{HWZ6}. In order to account for the developments, the paper was to a large extend rewritten. Rather than
dealing with the general theory, which also allows for boundary with corners, we restrict ourselves to a special case and illustrate it
with a discussion of stable maps, a topic closely related to Gromov-Witten theory.  We also would like to mention the paper \cite{FFW},
where the ideas of polyfold theory are explained as well.  The abstract theory has been applied in \cite{HWZ6} as part of the general construction of SFT.
In  \cite{SZ} it was used to address the Weinstein conjecture in higher dimensions. 
An extensive study of the case with boundary and corners 
is contained in \cite{HWZ7,HWZ8}. These extensions are crucial for an application to SFT. A basic paper in this direction 
is the upcoming \cite{FHWZ} in which the polyfolds relevant for SFT are being constructed. This is the place where the full power 
of the theory becomes apparent.

We shall start by introducing the category of stable maps. The construction of Gromov-Witten invariants can be understood
as the geometric study of perturbations of the full subcategory of $J$-holomorphic stable maps. We take a more general viewpoint and study the whole category. 
Our initial discussion is entirely topological, and  surprisingly the arising objects are so natural
that one is forced to raise the question if these structures go beyond topology. Indeed they do,  and they give rise to a generalized differential geometry as
well as a generalized nonlinear Fredholm theory accompanied by its own blend of nonlinear analysis. The whole package is referred to as polyfold theory.
The polyfold theory  provides a language and a large body of results to address questions arising 
in the study of moduli problems in symplectic geometry. It is clear that its applicability goes far beyond the latter field.
From a nonlinear analysis perspective the symplectic applications are concerned with the study of isomorphism classes of families
of nonlinear first order elliptic equations with varying domains and varying targets. The domains even are allowed to change the topology, and
bubbling-off phenomena will in general occur. This is the source for  compactness and transversality issues which make an algebraic counting of solutions 
very difficult. In the case of the problems arising in symplectic  geometry the polyfold theory overcomes these difficulties.
There is no doubt that the theory should have applications in other parts of nonlinear analysis as well.

It is a basic observation in symplectic geometry/topology that geometric questions can be rephrased
as  questions about  solution spaces of  nonlinear first order elliptic partial differential equation (the solution spaces are by definition the moduli spaces).
Very often just counting solutions suffices to answer a seemingly difficult geometric question.
The construction of the  moduli problems are difficult analytical problems involving analytical limiting behaviors, like bubbling-off, breaking of trajectories and stretching the neck.
The naive solution sets are very often not compact and not being cut out by the differential equation in a generic way. 
In most applications there exist intriguing compactifications of the solution spaces, which, however, are usually not compatible 
with the standard versions of smooth nonlinear analysis based on the notion of Fr\'echet differentiability. Moreover, as a consequence of local symmetries,
 very often
even a generic choice of geometric auxiliary data used  to construct the partial differential equation, will never result in a generic solution set.
The polyfold theory allows to view the partial differential equation within an abstract framework, which provides a Sard-Smale perturbation theory and deals with the transversality issues. The framework is so general that it also encompasses the geometric perturbation theory. Hence one can  proceed geometrically as
long as it is possible, and use the abstract perturbations only if the problems cannot be dealt with geometrically. In particular, whatever has been established
to work classically, will also work in this extended framework.

The analytical limiting  phenomena, which we mentioned above,  even assuming a sufficient amount of genericity, do not look like smooth phenomena if smoothness refers to the usual concept. Indeed, the coordinate changes are from a classical perspective usually nowhere differentiable.
However, it turns out that the notion of smoothness can be relaxed, and a generalization of differential geometry and nonlinear  functional analysis can be developed, for which the limiting phenomena can be viewed as smooth phenomena, even if they are quite often obscured by transversality issues.
In this generalized context the classical nonlinear Fredholm theory can be extended to a much larger class of spaces and operators, which can  deal with the aforementioned problems.

A great example to explain the theory comes from a study of the  category of stable maps associated to a symplectic manifold.
The study of these objects goes back to the seminal work by Gromov on pseudoholomorphic curve theory, see \cite{G}.
Gromov showed that the study of pseudoholomorphic maps is a powerful tool in symplectic geometry, by demonstrating its uses by many examples.
Kontsevich, later pointed out the importance of the notion of stable pseudoholomorphic curves, \cite{Ko}.
Our stable  maps need not to be pseudoholomorphic.  We find several natural constructions 
which at first glance are just of topological nature without exhibiting more regularity.  However, as it turns out, these are the shadows of smooth constructions, 
once the notion of differentiability in finite dimensions, which usually is generalized as Fr\'echet differentiability to infinite-dimensional Banach spaces,  is generalized in a quite  different way. Such a generalization requires an  additional piece of structure, called an sc-structure, which
occurs in interpolation theory, \cite{Tr}, albeit under a different name. In fact, we give a quite different interpretation of such a structure and make clear, that it can be viewed
as a generalization of a smooth structure on a Banach space. We call this generalization sc-smoothness, and the generalization of differentiability of a map we refer to as sc-differentiability.

The very interesting aspect is then the following fact. There are many sc-differentiable maps $r:U\rightarrow U$ satisfying $r\circ r=r$, i.e.
sc-smooth retractions. For Fr\'echet differentiability the image of such a retraction can be shown to be a submanifold of $U$.
Incidentally this is H. Cartan's last mathematical theorem, see \cite{H.Cartan}.
However, the images of sc-smooth retractions can be much more general. 
Most strikingly, they can have locally varying dimensions. Of course, a good notion of differentiability comes with the chain rule so that from $r\circ r=r$ we deduce $Tr=(Tr)\circ (Tr)$. In  other words, the tangent map
of an sc-smooth retraction is again an sc-smooth retraction. If a subset $O$ of an sc-Banach space is the image of an sc-smooth retraction $r$, then $TO=Tr(TU)$  defines
the tangent space, and it turns out, that the definition does not depend on the choice of $r$ as long as $r$ is sc-smooth and has $O$ as its image.
So we obtain quite general subsets of Banach spaces which have tangent spaces. An  sc-smooth map $f:O\rightarrow O'$, where $O\subset E$ and $O'\subset F$ are sc-smooth retracts, is a map such that
$f\circ r:U\rightarrow F$ is sc-smooth, where $r$ is an sc-smooth retraction onto $O$. As it turns out the definition does not depend on the choice of $r$. Further,  one verifies that $Tf:=T(f\circ r)|TO$  defines  a map $TO\rightarrow TO'$ between tangent spaces and  that the definition also does not depend on the choice of $r$.

\begin{figure}[htbp]
\mbox{}\\[2ex]
\centerline{\relabelbox
\epsfxsize 5.0truein \epsfbox{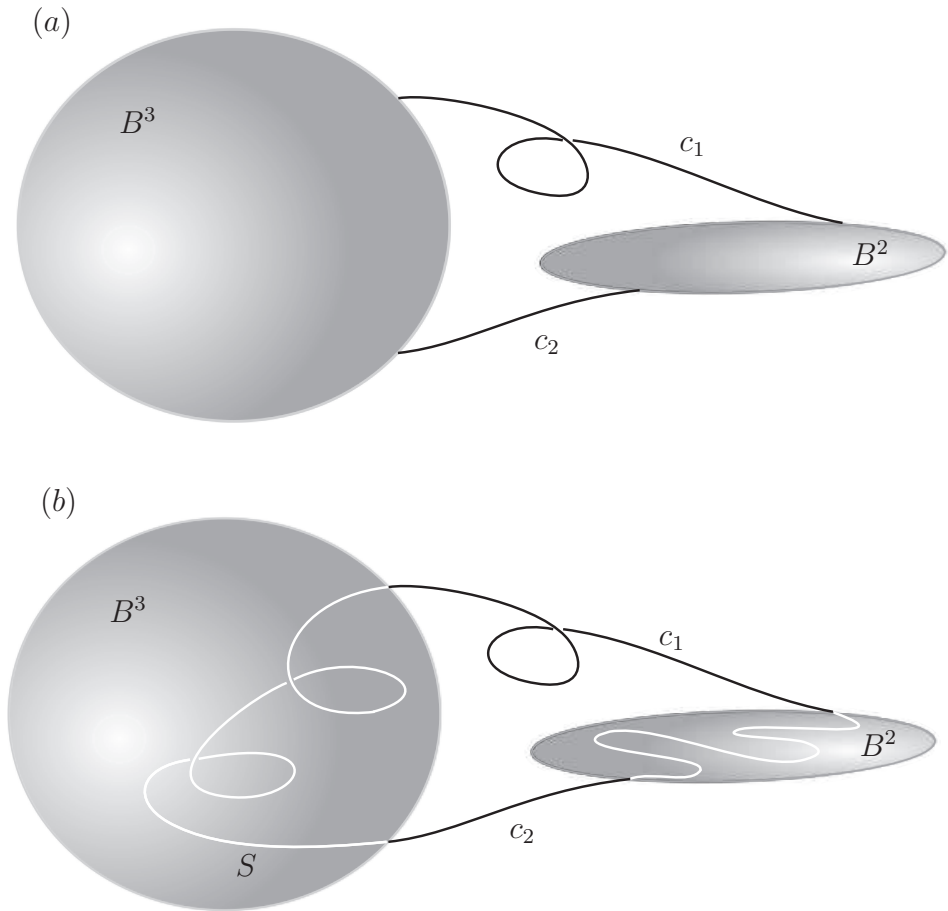}
\relabel {b1}{$B^3$}
\relabel {d1}{$B^2$}
\relabel {b}{$B^3$}
\relabel {d}{$B^2$}
\relabel {s}{$S$}
\relabel {c}{$c_1$}
\relabel {e2}{$c_2$}
\relabel {c1}{$c_1$}
\relabel {e1}{$c_2$}
\relabel {a}{$(a)$}
\relabel {a1}{$(b)$}
\endrelabelbox}
\caption{Figure a) shows a finite-dimensional M-polyfold $X$ which
is homeomorphic to the space consisting of the disjoint union of  an
open three-ball $B^3$ and an open two-ball $B^2$ connected by two
curves $c_1$, $c_2$. Figure b) shows the same M-polyfold containing
a one-dimensional $S^1$-like submanifold $S$. This submanifold could
arise as the zero set  of a transversal section of a strong
M-polyfold bundle $Y$ over $X$, which has varying dimensions.
Namely, over the three-ball  it is two-dimensional, over the two-disk
one-dimensional and otherwise it is trivial. The polyfold theory
would then guarantee a  natural smooth structure on the  solution
set $S$. }\label{porkbarrel}
\end{figure}
In summary, once we have a good notion of differentiability for maps between open sets, we also obtain a notion of differentiability for maps between sc-smooth retracts. However, for the usual notion of differentiability, smooth retracts
are manifolds and one does not obtain anything beyond the usual differential geometry and its standard generalization to Banach manifolds. On the other hand sc-differentiability
opens up new possibilities with serious applications. We generalize differential geometry by generalizing the notion of a manifold to that of an  M-polyfold. These  are metrizable spaces which are locally homeomorphic
to retracts with sc-smooth transition maps.  The theory as described in this paper even gives new objects in finite dimensions, see Figure \ref{porkbarrel}.
Most important for us is the fact,  that the new local models for a vastly generalized  differential geometry allow us to build spaces, which can be used to bring nonlinear partial differential equations which might show bubbling-off phenomena into an abstract geometric framework which allows for a very general nonlinear Fredholm theory with the usual expected properties.\\

\noindent{\bf Acknowledgement:}  The author thanks the Clay Institute for the opportunity to present this material at the 2008 Clay Research Conference.
Some of the initial work was done during a sabbatical at the Stanford mathematics department, and was supported in part by the American Institute of Mathematics. Thanks to P. Albers, Y. Eliashberg, J. Fish,  E. Ionel, K. Wysocki and E. Zehnder for many stimulating discussions.

\section{A Motivating Example and Natural Structures}\label{SEC22}
The following discussion starts with a closer  look  at the category of  stable Riemann surfaces and describes some of the interesting structures
appearing as a result of the classical Deligne-Mumford theory.   These structures and the associated viewpoint can be related to a topic arising in  
Gromov-Witten theory, namely the study of the  category of stable maps which are not necessarily pseudoholomorphic.
This category of stable maps has natural structures. Natural structures in mathematics are usually there for a good reason.
In our case  a deeper analysis reveals the existence of more general smooth models for an extended differential geometry as well as 
a vastly generalized Fredholm theory. The latter can be used to define the Gromov-Witten invariants. The same scheme also works
suitably extended in the  much more complicated framework of stable maps in SFT.

\subsection{The Category of Stable Noded Riemann Surfaces}
A good starting point  for explaining some of the later features  in the polyfold theory is a study of
stable noded Riemann surfaces with marked points. 
\begin{defn}
A stable, possibly noded Riemann surface with marked points, is a tuple $(S,j,M,D)$,
where $(S,j)$ is a closed smooth Riemann surface, $M\subset S$ is a finite collection of points, called marked points,
and a finite collection $D$ of un-ordered pairs $\{x,y\}$ of points $x,y\in S\setminus M$ having the following properties.
\begin{itemize}
\item[(i)]For  $\{x,y\}\in D$ it holds $x\neq y$. If $\{x,y\}\cap \{x',y'\}\neq \emptyset$, then $\{x,y\}=\{x',y'\}$.
\item[(ii)] The topological space $\bar{S}$ obtained by identifying $x\equiv y$ for $\{x,y\}\in D$  is connected.
\item[(iii)] Define the subset $|D|$ by $|D|=\bigcup_{\{x,y\}\in D} \{x,y\}.$
For every connected component $C$ of $S$ having genus $g(C)$ with $n_C$ being the number of points in $C\cap (|D|\cup M)$ it holds $2g(C)+n_C\geq 3$.
\end{itemize}
\end{defn}
Condition (iii) is called the stability condition. In our case the set $M$ is not ordered, but the following discussion 
could be carried out in the ordered case as well. 

We can view $\alpha=(S,j,M,D)$ as an object in a category ${\mathcal R}$
where the morphisms $\Phi:\alpha\rightarrow \alpha'$ are given by
$$
\Phi=(\alpha,\phi,\alpha'):\alpha\rightarrow \alpha'
$$
and $\phi:(S,j)\rightarrow (S',j')$ is a biholomorphic map satisfying $\phi(M)=M'$ and
$\phi(D)=D'$, where $\phi(D):=\{\{\phi(x),\phi(y)\}\ |\ \{x,y\}\in D\}$.
\begin{defn}
We call ${\mathcal R}$ the category of stable, (possibly) noded Riemann surfaces with un-ordered marked points.
\end{defn}
The following result is well known, and is part of the Deligne-Mumford theory. The theory has been described from a more differential geometric perspective 
in \cite{RS}, and is being described in \cite{HWZ-DM} from a polyfold perspective. 
\begin{thm}\label{THH1}
The category ${\mathcal R}$ has the following properties.
\begin{itemize}
\item[(i)]  Every morphism in the category ${\mathcal R}$ is an isomorphism. 
\item[(ii)] The stability condition implies that between any two objects there are at most finitely many morphisms. 
\item[(iii)] The orbit space $|{\mathcal R}|$ of the category ${\mathcal R}$, which is the set of isomorphism classes
of objects in ${\mathcal R}$ carries a natural metrizable topology.
\item[(iv)]  The topological space $|{\mathcal R}|$ carries in a  natural way the structure of a holomorphic orbifold for which each connected component is compact.
\end{itemize}
\end{thm}
A basis for the topology is given in Proposition \ref{propp1}.
We shall need a variation of this result, and shall describe a particular approach to ${\mathcal R}$ which one might view as a toy case for the polyfold approach to stable maps in symplectic manifolds.  We describe some of the results we shall need later on, and refer the reader for more details to \cite{HWZ6,HWZ-DM}, or the references mentioned therein.

Given an object $\alpha=(S,j,M,D)$ in ${\mathcal R}$ it is not just a point in the category, but it has geometry as well. In fact, the additional structures
on the  category come from the fact that every object has its own geometry, which is being exploited by constructing associated objects 
by a plumbing procedure. Recall that the automorphism group $G$ of $\alpha$ is  the finite collection of all morphisms (which are all isomorphisms)
$\Phi:\alpha\rightarrow \alpha$.
\begin{defn}
Let $\alpha=(S,j,M,D)$ be an object in ${\mathcal R}$ with automorphism group $G$.  A small disk structure associated to $\alpha$ and denoted by ${\bf D}$ assigns to every $x\in |D|$ a compact disk-like neighborhood
$D_x\subset S$ around the point $x$, so that the following holds.
\begin{itemize}
\item[(i)] $D_x$ has a smooth boundary and the disks are mutually disjoint.
\item[(ii)] $\bigcup_{x\in |D|} D_x$ is invariant under $G$. 
\item[(iii)] $D_x\cap M=\emptyset$ for all $x\in |D|$.
\end{itemize}
\end{defn}
Given a small disk structure ${\bf D}$ we have for every $\{x,y\}\in D$ an associated noded disk pair $(D_x\cup D_y,\{x,y\})$.
From this data we can construct for every $\{x,y\}$ a natural gluing parameter. This is done as follows.
Denote for $x\in |D|$ by $\what{x}$ an oriented real line  in $(T_xS,j)$. If ${\mathbb S}^1\subset {\mathbb C}$ denotes the unit circle, 
we observe that  given $\theta\in {\mathbb S}^1$, we can define a new oriented real line $\theta\what{x}$, making use of the complex structure $j$.
We introduce an equivalence class of unordered pairs of real oriented lines $\{\what{x},\what{y}\}$ associated to $\{x,y\}\in D$
as follows. We say $\{\what{x},\what{y}\}$ is equivalent to $\{\what{x}',\what{y}'\}$ (here the lines lie over $x$ and $y$) provided there exists 
$\theta\in {\mathbb S}^1$ satisfying
$$
\what{x}=\theta\what{x}'\ \text{and}\ \what{y}=\theta^{-1}\what{y}'.
$$
Equivalence classes will be written as $[\what{x},\what{y}]$.  Clearly the collection of all equivalence classes associated to the pair $\{x,y\}$ can be parametrized by the unit circle via
$$
{\mathbb S}^1\ni \theta\rightarrow [\what{x},\theta\what{y}].
$$
\begin{defn}
An equivalence class $[\what{x},\what{y}]$ is called a decorated nodal pair with underlying  nodal pair $\{x,y\}$.
\end{defn}
Hence above every $\{x,y\}\in D$ there lies an ${\mathbb S}^1$-worth of decorated nodal pairs $[\what{x},\what{y}]$.
We shall refer to $[\what{x},\what{y}]$ as a decoration for $\{x,y\}$.

Consider a formal expression $r[\what{x},\what{y}]$ with $r\in [0,1/2)$.  Two such formal expressions are considered the same, i.e. $r[\what{x},\what{y}]=r'[\what{x}',\what{y}]$
provided one of the following holds:
\begin{itemize}
\item[(i)] $r=r'=0$, or
\item[(ii)]  $r=r'\neq 0$ and $[\what{x},\what{y}]=[\what{x}',\what{y}']$.
\end{itemize}
  \begin{defn}
  A natural gluing parameter associated to the nodal pair $\{x,y\}$ is a formal expression $\mathfrak{a}^{\{x,y\}}=r[\what{x},\what{y}]$ where $r\in [0,1/2)$,
  with the notion of equality as defined above.   In case $r=0$ we shall write $\mathfrak{a}^{\{x,y\}}=0$. 
  \end{defn}
  We define the set  ${\mathbb B}^{\{x,y\}}(1/2)$ by
$$
{\mathbb B}^{\{x,y\}}(1/2)=\{r[\what{x},\what{y}]\ |\ r\in [0,1/2),\ [\what{x},\what{y}]\ \text{a decoration of}\ \{x,y\}\}
$$
and call it the set of natural gluing parameters associated to the nodal pair $\{x,y\}$.  The modulus of $\mathfrak{a}^{\{x,y\}}=r[\what{x},\what{y}]$
is the number $r$ and we shall write $|\mathfrak{a}^{\{x,y\}}|=r$.
\begin{defn}
Given $\alpha$, a  natural  (total) gluing parameter is a map
$$
\mathfrak{a}: D\ni \{x,y\}\rightarrow  \mathfrak{a}^{\{x,y\}}\in {\mathbb B}^{\{x,y\}}(1/2).
$$
The collection of all gluing parameters associated to $\alpha$ will be denoted by ${\mathbb B}^\alpha$.
\end{defn}
The dependency of ${\mathbb B}^\alpha$ on $\alpha$ is as follows. It only depends on the complex multiplication
on the tangent spaces $T_xS$ for $x\in |D|$.
We note the following result having an easy proof which is left to the reader.
\begin{prop}
Assume that $\alpha$ is given. For every $\{x,y\}\in D$ the set ${\mathbb B}^{\{x,y\}}(1/2)$ has in a natural way the structure of a holomorphic manifold. It is characterized by the following 
property. For every choice of oriented real lines $\what{x}$ and $\what{y}$ the map
$$
{\mathbb B}^{\{x,y\}}(1/2)\rightarrow \{z\in {\mathbb C}\ |\ |z|<1/2\}:r[\what{x},\theta\what{y}]\rightarrow r\theta
$$
is biholomorphic. In particular the set ${\mathbb B}^\alpha$, of all natural gluing parameters $\mathfrak{a}:\{x,y\}\rightarrow \mathfrak{a}^{\{x,y\}}$ has a natural holomorphic structure
and is biholomorphic to a product of open disks.
\end{prop}
We shall use the natural gluing parameters to modify $\alpha$ and produce new stable noded Riemann surfaces.
In this context an important concept is that of a gluing profile.
\begin{defn}
A gluing profile is a smooth diffeomorphism $\varphi:(0,1]\rightarrow [0,\infty)$.
\end{defn}
We shall refer to 
$$
\varphi(r)=-\frac{1}{2\pi}\cdot \ln(r)
$$
 as the logarithmic gluing profile, and to 
 $$
 \varphi(r)=e^{\frac{1}{r}}-e
 $$
  as the exponential gluing profile.
Restricting $r[\what{x},\what{y}]$ to the case with $r\in [0,1/2)$ has to do with the choice of these two gluing profiles. For an arbitrary 
gluing profile $\varphi$ we only have to require that $r\in [0,\varepsilon_\varphi)$, where $\varepsilon_\varphi >0$, depending on the profile,  is small enough.

Assume we are given a small disk structure ${\bf D}$ and  a gluing profile $\varphi$.  We can define a gluing or plumbing construction as follows. 
For a gluing parameter $\mathfrak{a}^{\{x,y\}}=r[\what{x},\what{y}]$ associated to $\{x,y\}\in D$, we consider the disks $D_x$,  $D_y$, and the nodal pair $\{x,y\}$. If $\mathfrak{a}^{\{x,y\}}=0$
we keep this data, i.e. with other words the gluing of $(D_x\cup D_y,\{x,y\})$ with $a^{\{x,y\}}=0$ produces
the noded disk $(D_x\cup D_y,\{x,y\})$.
If $\mathfrak{a}^{\{x,y\}}\neq 0$ we compute the value $R=\varphi(r)$ and take a representative $\{\what{x},\what{y}\}$
of the class $[\what{x},\what{y}]$. After this choice,  there exist uniquely determined biholomorphic maps
$$
h_x:(D_x,x)\rightarrow ({\mathbb D},0),\ Th_x(x)\what{x}={\mathbb R},
$$
and
$$
h_y:(D_y,y)\rightarrow ({\mathbb D},0),\ Th_y(x)\what{y}={\mathbb R}.
$$
Here ${\mathbb D}$ is the closed unit disk in ${\mathbb C}$, and ${\mathbb R}\subset T_0{\mathbb D}$ is oriented by $1$.
Recall that the modulus of the annulus $A_{r_1,r_2}:=\{z\in {\mathbb C}\ |\ r_1\leq |z|\leq r_2\}$ is defined by
$$
\text{modulus}(A_{r_1,r_2}) = \frac{1}{2\pi}\cdot \ln(r_2/r_1).
$$
If we are given a Riemann surface $(A,j)$ where $A$ is compact, has smooth boundary components, and is diffeomorphic to an annulus, we know from the uniformization theorem that $(A,j)$ is biholomorphic to $A_{r,1}$ for a uniquely determined $0<r<1$. The modulus of $(A,j)$ by definition is given by
$$
\text{modulus}(A,j)=-\frac{1}{2\pi}\cdot \ln(r).
$$
We note that the right-hand side is precisely $R=\varphi(r)$ for the logarithmic gluing profile. In the next step the choice of the gluing profile matters,
and a remark, why we consider other gluing profiles than the logarithmic one, is  in order. 
\begin{rem}
In later constructions we shall see that the 
sc-smooth structure we are going to construct will depend on the choice of $\varphi$. The sc-smooth structure also influences if certain operators are sc-smooth Fredholm operators. The operators we are interested in will be sc-smooth for the exponential gluing profile, but not for the logarithmic one.
\end{rem}
With some gluing profile $\varphi$ fixed,  we pick concentric compact annuli $A_x\subset D_x$ and 
$A_y\subset D_y$ with one boundary component being the boundary of $D_x$ and $D_y$, respectively,
so that $A_x$ and $A_y$ have modulus $R=\varphi(r)$. We discard the points $D_x\setminus A_x$ and $D_y\setminus A_y$ from $S$ and the nodal pair
$\{x,y\}$. With $A_x\cup A_y\subset S\setminus ((D_x\setminus A_x)\cup (D_y\setminus A_y))$  we identify $A_x$ with $A_y$ by identifying $z\equiv z'$, where  $(z,z')\in A_x\times A_y$ satisfying
$$
h_x(z)\cdot h_y(z')= e^{-2\pi\cdot R}.
$$
The definition of the identification $z\equiv z'$ is independent of the choice of the representatives $\{\what{x},\what{y}\}$ in $[\what{x},\what{y}]$.
Let us denote by ${\mathcal Z}_{\mathfrak{a}^{\{x,y\}}}$ the glued space obtained by identifying the two annuli. If $\mathfrak{a}^{\{x,y\}}=0$ this is the noded $(D_x\cup D_y,\{x,y\})$ and for 
$\mathfrak{a}^{\{x,y\}}\neq 0$
it is the space obtained by identifications just described. Clearly, in the latter case we have natural biholomorphic maps
$$
A_x\rightarrow {\mathcal Z}_{\mathfrak{a}^{\{x,y\}}}\ \text{and}\ A_y\rightarrow  {\mathcal Z}_{\mathfrak{a}^{\{x,y\}}},
$$
by associating to a point its equivalence class. 
Carrying out the above procedure for every noded disk pair $(D_x\cup D_y,\{x,y\})$ defines for a gluing parameter $\mathfrak{a}\in {\mathbb B}^\alpha$ a new noded stable Riemann surface $(S_\mathfrak{a},j_\mathfrak{a},M_\mathfrak{a},D_\mathfrak{a})$ with marked points.
We just described how to obtain $S_\mathfrak{a}$. We denote by $j_\mathfrak{a}$ the natural smooth almost complex structure on the latter.
Observe that $z\equiv z'$ is a holomorphic identification. By $M_\mathfrak{a}$ we just mean $M$ naturally identified with  a subset of $S_\mathfrak{a}$.
The collection of nodal points $D_\mathfrak{a}$ is obtained by first removing from $D$ all $\{x,y\}$ with $\mathfrak{a}^{\{x,y\}}\neq 0$
and then identifying this set with the obvious set of nodal pairs for $S_\mathfrak{a}$. For every $\{x,y\}\in D_\mathfrak{a}$ we still have
$D_x\cup D_y\subset S_\mathfrak{a}$.

At this point we can exhibit a basis of open neighborhoods for a given point $|\alpha|\in |{\mathcal R}|$ for the natural metrizable topology.
For a representative $\alpha$ and an associated gluing parameter $\mathfrak{a}$ we define  $|\mathfrak{a}|$ by
$$
|\mathfrak{a}|=\text{max}_{\{x,y\}\in D} \left | \mathfrak{a}^{\{x,y\}}\right|.
$$
Write $\alpha=(S,j,M,D)$ and fix a small disk structure ${\bf D}$. The set of all smooth almost complex structures $k$ on $S$,  which define the same orientation as $j$, has a natural
metrizable topology measuring the $C^\infty$-distance. Denote a choice of metric by $\rho$.  For a given gluing profile $\varphi$ and $\varepsilon>0$ small enough define
$$
{\mathcal U}(\alpha,{\bf D},\varphi,\varepsilon)=\{|(S_\mathfrak{a},k_\mathfrak{a},M_\mathfrak{a},D_\mathfrak{a})|\ |\ \rho(j,k)<\varepsilon,\ j=k\ \text{on}\ {\bf D},\ |\mathfrak{a}|<\varepsilon\}.
$$
\begin{thm}\label{propp1}
Fix a gluing profile $\varphi$, pick for every class $c\in |{\mathcal R}|$ an object $\alpha_c$, and  choose an associated small disk structure
${\bf D}_c$ for $\alpha_c$. Consider ${\mathcal U}(\alpha_c,{\bf D}_c,\varphi,\varepsilon)$ for $0<\varepsilon<\varepsilon_c$ and denote this collection by $\mathfrak{U}_c$.
Then the union $\mathfrak{U}=\bigcup_{c\in |{\mathcal R}|} \mathfrak{U}_c$ is the basis for a metrizable topology ${\mathcal T}$ on $|{\mathcal R}|$.
The topology ${\mathcal T}$ does not depend on the choices involved in constructing $\mathfrak{U}$.
\end{thm}
The topology ${\mathcal T}$ is the one referred to in Theorem \ref{THH1}.
The automorphism group $G$ acts on the set of gluing parameters ${\mathbb B}^\alpha$ via biholomorphic maps 
$$
G\times {\mathbb B}^\alpha\rightarrow {\mathbb B}^\alpha: (g,\mathfrak{a})\rightarrow g\ast\mathfrak{a},
$$
in a natural way as follows. Namely  $\mathfrak{b}=g\ast \mathfrak{a}$ is defined by 
$$
\mathfrak{b}^{\{g(x),g(y)\}}  = r\cdot [Tg(x)\what{x},Tg(y)\what{y}],
$$
where $r[\what{x},\what{y}]=\mathfrak{a}^{\{x,y\}}$. With the above construction an element $g\in G$ induces a biholomorphic map
$$
g_\mathfrak{a}:(S_\mathfrak{a},j_\mathfrak{a},M_\mathfrak{a},D_\mathfrak{a})\rightarrow (S_\mathfrak{b},j_\mathfrak{b},M_\mathfrak{b},D_\mathfrak{b}),
$$
where $\mathfrak{b}=g\ast \mathfrak{a}$. It also holds
$$
h_{g\ast\mathfrak{a}}\circ g_\mathfrak{a}= (h\circ g)_\mathfrak{a}.
$$
Let us abbreviate $(S_\mathfrak{a},j_\mathfrak{a},M_\mathfrak{a},D_\mathfrak{a})$ by $\alpha_{\mathfrak{a}}$. Then 
$\alpha_\mathfrak{a}=\alpha_{{\mathfrak{a}'}}$ if and only if $\mathfrak{a}=\mathfrak{a}'$. Consider for these objects $\alpha_\mathfrak{a}$
with $|\mathfrak{a}|<1/2$ 
the associated full subcategory ${\mathcal A}={\mathcal A}_{\alpha,{\bf D},\varphi}$  of ${\mathcal R}$. This subcategory is small.
Given the object $\alpha$ with automorphism group $G$,   every choice of a small disk structure (a set of possible choices) produces 
such a subcategory. 
We also consider the translation groupoid $G\ltimes {\mathbb B}^{\alpha}$.
The latter is a small category with object set ${\mathbb B}^\alpha$ and morphism set $G\times {\mathbb B}^{\alpha}$, where $(g,\mathfrak{a})$ is seen as a morphism
$\mathfrak{a}\rightarrow g\ast\mathfrak{a}$. Observe that object set and morphism set have in our case both smooth manifold structures.
There exists the obvious functor
$$
\gamma:G\ltimes {\mathbb B}^{\alpha}\rightarrow {\mathcal A}
$$
which maps the object $\mathfrak{a}$ to $\alpha_\mathfrak{a}$ and the morphism $(g,\mathfrak{a}):\mathfrak{a}\rightarrow g\ast \mathfrak{a}$
to 
$$
(\alpha_\mathfrak{a},g_\mathfrak{a},\alpha_{g\ast\mathfrak{a}}):\alpha_\mathfrak{a}\rightarrow \alpha_{g\ast\mathfrak{a}}.
$$
We note that $\gamma$ is injective on objects. More is true, namely 
from the results of the Deligne-Mumford theory we obtain the following proposition which is true for any gluing profile. It has to be viewed as an intermediate result
neglecting for the moment the possibility of  deforming the (integrable) almost complex structure $j$ associated to $\alpha$.
\begin{prop}
For a suitable $G$-invariant open neighborhood $U$ of $0\in {\mathbb B}^\alpha$ depending on the choice of the gluing profile,
the functor
$$
\gamma: G\ltimes U\rightarrow {\mathcal R}
$$
is injective on objects and fully faithful.
\end{prop}
The above proposition describes  the family we can obtain just from gluing (or plumbing) without considering deformations of the (integrable almost) complex structure.
 We can also obtain from $\alpha$ new objects by changing the  complex structure.
Of course, a priori, we only would like to consider only those deformations which produce non isomorphic objects. There is, however,  an inherent difficulty coming from the automorphisms.
Given the object $\alpha$ in the category ${\mathcal R}$ consider the complex vector space
$\Gamma_0(\alpha)$ consisting of smooth sections of $TS\rightarrow S$, which vanish at the points in $M\cup |D|$.
By $\Omega^{0,1}(\alpha)\rightarrow S$ we denote the complex vector space of smooth sections of $\text{Hom}_{\mathbb R}(TS,TS)\rightarrow S$
which are complex anti-linear. The Cauchy-Riemann operator defines
a linear operator
$$
\bar{\partial}:\Gamma_0(\alpha)\rightarrow \Omega^{0,1}(\alpha).
$$
Define via the formula 
$$
g_a:=1+\sharp D+\sum_C [g(C)-1],
$$
where the sum is taken over all connected components $C$  of $S$, the arithmetic genus $g_a$.
The basic fact is given by the following proposition.
\begin{prop}
The Cauchy-Riemann operator is a complex linear differential operator.
The stability property of $\alpha$ implies that $\bar{\partial}$ is injective 
and as a consequence of Riemann-Roch the cokernel has complex dimension
$3g_a+\sharp M -\sharp D-3$.
\end{prop}
 In the next step we consider the deformation space
of $\alpha$ for fixed combinatorial type. This means we consider the deformation of the complex structure, the marked points, and the 
nodal pairs, but do not remove nodes. The natural deformation space introduced in the following is the infinitesimal version of this,
if we divide out by isomorphisms, but keep track of the automorphisms.

\begin{defn}
The (infinitesimal) natural deformation space (of fixed combinatorial type)  of $\alpha=(S,j,M,D)$ is by definition the complex vector space
$H^1(\alpha)$ defined by
$$
H^1(\alpha)=\Omega^{0,1}/(\bar{\partial}\Gamma_0(\alpha)).
$$
\end{defn}
The construction $\alpha\rightarrow H^1(\alpha)$ is functorial. Namely we associate to an object $\alpha$
a finite-dimensional complex vector space, and to a morphism $\phi:\alpha\rightarrow \alpha'$ a complex vector space isomorphism
$$
H^1(\phi):H^1(\alpha)\rightarrow H^1(\alpha'):[\tau]\rightarrow [\phi\circ \tau \circ\phi^{-1}].
$$
In particular the automorphism group $G$ of the object $\alpha$ has an action on $H^1(\alpha)$.

Assume that $V$ is an open subset of some vector space $E$ and  $V\in v\rightarrow j(v)$ a smooth family of almost complex structures on $S$.
Then  it holds $j(v)^2=-Id$, and differentiating 
at $v$ into the direction $\delta v\in E$  gives the identity 
$$
j(v)(Dj(v)\delta v)=-(Dj(v)\delta v)j(v).
$$
This means that for every $z\in S$ the map  
$$
(Dj(v)\delta v)_z:(T_zS,j(v))\rightarrow (T_zS,j(v))
$$
 is complex anti-linear and 
induces an element $[Dj(v)\delta v]\in H^1(S,j(v),M,D)$. Hence we obtain for every $v\in V$ a linear map
$$
[Dj(v)]:E\rightarrow H^1(S,j(v),M,D).
$$
This map is called the Kodaira-Spencer differential associated to $v\rightarrow j(v)$. 
We define $\alpha_v$ by
$$
\alpha_v=(S,j(v),M,D).
$$
\begin{defn}
 Given a stable $\alpha$ and an associated small disk structure ${\bf D}$ we shall call a smooth family $v\rightarrow j(v)$
 of almost complex structures on $S$  defined on a $G$-invariant open neighborhood $V$ of $0$ in $H^1(\alpha)$ a good deformation
 compatible with ${\bf D}$ provided it has the following  properties.
 \begin{itemize}
 \item[(i)] $j(0)=j$.
 \item[(ii)] $j(v)=j$ on all disks $D_x$ associated to ${\bf D}$.
 \item[(iii)] For every $v \in V$ the Kodaira-Spencer differential 
 $$
 [Dj(v)]:H^1(\alpha)\rightarrow H^1(\alpha_v)
 $$
  is a complex linear isomorphism.
 \item[(iv)] For every $g\in G$ and $v\in V$ the map
 $g:\alpha_v\rightarrow \alpha_{g\ast v}$
 is biholomorphic. 
\end{itemize}
 \end{defn}
 A well-known result with a proof given in \cite{HWZ-DM} is the following.
 \begin{thm}
 Given a stable $\alpha$ in ${\mathcal R}$ with automorphism group $G$, and a small disk structure ${\bf D}$, there always exist good deformations compatible with ${\bf D}$.
 \end{thm}
Suppose we start with a stable $\alpha$ with automorphism group $G$ and have fixed a small disk structure ${\bf D}$.
Then we can, as we have seen, construct the functor
$$
\gamma:G\ltimes {\mathbb B}^\alpha\rightarrow {\mathcal R}.
$$
Now taking a good deformation as described in the previous definition, say $V\ni v\rightarrow j(v)$, we can construct the following functor, where we take either the 
logarithmic or exponential gluing profile.
$$
\Psi: G\ltimes (V\times {\mathbb B}^\alpha)\rightarrow {\mathcal R}. 
$$
To objects $(v,\mathfrak{a})$ we assign
$$
\alpha_{(v,\mathfrak{a})}=(S_\mathfrak{a},j(v)_\mathfrak{a},M_\mathfrak{a},D_\mathfrak{a}).
$$
The morphism $(g,(v,\mathfrak{a})):(v,\mathfrak{a})\rightarrow (g\ast v,g\ast\mathfrak{a})$ is mapped to
$$
(\alpha_{(v,\mathfrak{a})},g_\mathfrak{a}, \alpha_{(g\ast v,g\ast\mathfrak{a})}):\alpha_{(v,\mathfrak{a})}\rightarrow \alpha_{(g\ast v,g\ast \mathfrak{a})}.
$$
The basic result is the following. A proof is given in \cite{HWZ-DM}.
\begin{thm}\label{thmm1}
Given the stable object $\alpha$ in ${\mathcal R}$ with automorphism group $G$, a small disk structure ${\bf D}$, and a good deformation $v\rightarrow j(v)$ compatible with ${\bf D}$, there exists
an open $G$-invariant neighborhood $O$ of $(0,0)\in V\times {\mathbb B}^\alpha$, so that the following holds for a given gluing profile.
\begin{itemize}
\item[(i)] $\Psi:G\times O\rightarrow {\mathcal R}$ is a fully faithfull functor.
\item[(ii)] The map $|\Psi|:{_{G}\backslash} O \rightarrow |{\mathcal R}|$ induced on orbit spaces  defines a homeomorphism onto an open neighborhood
$U$ of $|\alpha|$.
\item[(iii)] For every $(v,\mathfrak{a})$ the partial Kodaira-Spencer differentiable associated to $\alpha_{(v,\mathfrak{a})}$ is an isomorphism.
\item[(iv)] For every point $q\in O$ there exists an open neighborhood $U(q)\subset O$ with the property that every sequence $(q_k)\subset U(q)$
 for which $|\Psi(q_k)|$ converges in $|{\mathcal R}|$ there exists a subsequence of $(q_k)$ converging in $\cl_O(U(q))$.
\end{itemize}
\end{thm} 
For the definition of at the partial Kodaira-Spencer differentiable see \cite{HWZ6} or \cite{HWZ-DM}.
The above construction is possible for every given gluing profile $\varphi$.  
\begin{defn}\label{Basic_def}
Assume a gluing profile $\varphi$ is fixed. A functor $\Psi:G\ltimes O\rightarrow {\mathcal R}$ constructed as described previously, associated to a stable $\alpha$ and a small disk structure ${\bf D}$, and satisfying the properties (i)--(iv) of Theorem \ref{thmm1}
is called a good uniformizer associated to $\alpha$.
\end{defn}
For the following discussion we assume that a gluing profile $\varphi$ has been fixed.
 Let us observe that the construction of $\Psi:G\ltimes O\rightarrow {\mathcal R}$ is associated to a given object $\alpha$ in ${\mathcal R}$
 and 
 $$
 \Psi(0,0)=\alpha\ \text{and}\  \Psi((g,(0,0))=(\alpha,g,\alpha):\alpha\rightarrow \alpha.
 $$
 The constructions requires the choice of a small disk structure ${\bf D}$, the good deformation 
 $v\rightarrow j(v)$, and the open $G$-invariant subset $O$ of $V\times {\mathbb B}^\alpha$. Here ${\mathbb B}^\alpha$ is a well-defined set associated to $\alpha$ and $V$ is an open subset
 of the complex vector space $H^1(\alpha)$ associated to $\alpha$. All in all, our possible choices involved in constructing $\Psi$ for a fixed $\alpha$  constitute a set. Therefore we can make the following definition.
 \begin{defn}
By  $F(\alpha)$  we denote the set of good uniformizers associated to $\alpha$.
\end{defn}
We shall elaborate next about the functorial properties of $F$.
 \begin{defn}
By  ${\mathcal R}^-$ we denote the category which has the same objects as ${\mathcal R}$, but as morphisms
 only the identities.
 \end{defn}
  We can view the assignment
 $$
 \alpha\rightarrow F(\alpha).
 $$
 as a functor defined on ${\mathcal R}^-$ with values in the category of sets  $\text{SET}$. This is the  `minimalistic version' of our construction giving us the functor
  $$
F:{\mathcal R}^-\rightarrow \text{SET}: \alpha\rightarrow F(\alpha)\ \text{(on objects)},
 $$
 which  on morphisms
 $$
{F}(Id_\alpha) = Id_{F(\alpha)}.
  $$
  We said that this is the `minimalistic version' since there is, in fact, much more structure.
Let us explore this now.

Recall that 
for  $\alpha=(S,j,M,D)$  an element $\Psi$ in the set $F(\alpha)$ is constructed by making the following choices.
\begin{itemize}
\item[(i)] A  small disk structure ${\bf D}$  for $(S,j,M,D)$.
\item[(ii)] A deformation of $j$ with specific properties
 $v\rightarrow j(v)$ on $(S,M,D)$, i.e. which is constant on ${\bf D}$, has some symmetry properties, and is defined on a $G$-invariant open neighborhood $V$ of $0$ in $H^1(\alpha)$. 
 \item[(iii)] From this data we obtain the natural gluing parameters ${\mathbb B}^\alpha$ and can construct the family
 $$
\Psi: (v,\mathfrak{a})\rightarrow (S_\mathfrak{a},j(v)_\mathfrak{a},M_\mathfrak{a},D_\mathfrak{a}).
 $$
 \item[(iv)]  A $G$-invariant open neighborhood $O$ of $(0,0)$ in $V\times B^{\alpha}$
 so that the restriction of the above family satisfies a list of properties defining an element in $F(\alpha)$, see Definition \ref{Basic_def}.
 \end{itemize}
Hence the collection $F(\alpha)$ is obtained via a precise instruction which requires choices in its construction.
 When we consider isomorphic objects in ${\mathcal R}$, say $\alpha$ and $\alpha'$, there is a natural way 
 of matching a specific choice for $\alpha$ with a specific choice for $\alpha'$. This matching depends
 on the choice of an isomorphism and has functorial properties, so that one can upgrade 
 $F:{\mathcal R}^-\rightarrow \text{SET}$ to a functor $F:{\mathcal R}\rightarrow \text{SET}$.
 More precisely, 
 assume that 
 $$
 \Phi:=(\alpha,\phi,\alpha'):\alpha\rightarrow \alpha'
 $$
  is a morphism. Then we define the following data for $\alpha'$.
 \begin{itemize}
 \item[(i)] ${\bf D}':=\phi_\ast({\bf D}):=\{\phi(D_z)\ |\ z\in |D|\}$.
 \item[(ii)] The biholomorphic map $\phi$ defines a linear isomorphism
 $\phi_\ast :H^1(\alpha)\rightarrow H^1(\alpha')$ and we take $V'=\phi_\ast (V)$, and  a biholomorphic map $\phi_\ast: {\mathbb B}^\alpha\rightarrow {\mathbb B}^{\alpha'}$.
 \item[(iii)]  For $j'$ we take the map 
 $V'\ni v'\rightarrow j'(v'):=T\phi\circ j(\phi_\ast^{-1} v')\circ T\phi^{-1}.$
 \item[(iv)] $O'=\phi_\ast(O)$
 \end{itemize}
From the data ${\bf D}'$, $v'\rightarrow j'(v')$ and ${\mathbb B}^{\alpha'}$ we construct the element $\Psi'\in F(\alpha')$ corresponding 
to $\Psi$.  This shows that a morphism $\Phi:=(\alpha,\phi,\alpha'):\alpha\rightarrow \alpha'$
defines a bijection 
$$
F(\Phi):F(\alpha)\rightarrow F(\alpha').
$$
Hence we  obtain a functor
$$
F:{\mathcal R}\rightarrow \text{SET}.
$$
Just saying that $F$ is functor from ${\mathcal R}$ into $\text{SET}$  is again quite minimalistic,
since we suppress the various maps $\phi_\ast$ we constructed above. 
 Hence, underlying this bijection 
there is in our case even some fine structure, which we shall describe now.
 
Assume that the automorphism group of  $\alpha=(S,j,M,D)$ is $G$. Then the elements
of $G$ have the form $(\alpha,g,\alpha)$, where $g:\alpha\rightarrow \alpha$ is a biholomorphic map.
For simplicity of notation we identity $g=(\alpha,g,\alpha)$. Assume that 
we are given an isomorphism
 $$
 \Phi=(\alpha,\phi,\alpha'):\alpha\rightarrow \alpha'
 $$
 Having the domains fixed (most importantly the almost complex structures) we identify $\phi=(\alpha,\phi,\alpha')$. Denote by  $G'$ the automorphism group of $\alpha'$, and proceeding as in the case of $G$,  the morphism
$\Phi$ determines
 a group isomorphism 
 $$
 \gamma_\Phi:G\rightarrow G':g\rightarrow \phi\circ g\circ \phi^{-1}.
 $$
  Consider the, by the previous discussion,
 related elements $\Psi\in F(\alpha)$ and $\Psi'\in F(\alpha')$ given by
 $$
 \Psi:G\ltimes O\rightarrow {\mathcal R}\ \ \text{and}\ \ \Psi':G'\ltimes O\rightarrow {\mathcal R}.
 $$
 Related means that $F(\Phi)(\Psi)=\Psi'$.
 The biholomorphic map $\phi$ defines an equivariant diffeomorphism  $O\rightarrow O'$ which satisfies
 $$
\phi_\ast (g\ast q)=\gamma_{\Phi} (g)\ast \phi_\ast(q),
$$
where $q=(v,\mathfrak{a})$. This in particular means that $\phi_\ast$ defines a smooth equivalence 
between the translation groupoids, 
 but also is a bijection on objects and morphisms 
$$
f_\Phi :G\ltimes O\rightarrow G'\ltimes O' :(q,q)\rightarrow (\gamma_{\Phi}(g),\phi_\ast(q)).
$$
From this we see that we have the two functor
$$
\Psi\ \text{and}\ \Psi'\circ f_\Phi 
$$
both defined on $G\ltimes O$. On objects, with
$(v',\mathfrak{a}')=f_\Phi (v,\mathfrak{a})$, we have the biholomorphic map
$\phi_\mathfrak{a}: \alpha_{(v,\mathfrak{a})}\rightarrow \alpha_{(v',\mathfrak{a}')}'$, and
hence the isomorphism
$$
\Phi_{(v,\mathfrak{a})}:=(\alpha_{(v,\mathfrak{a})} ,\phi_\mathfrak{a},\alpha_{f_\Phi(v,\mathfrak{a})}'):\alpha_{(v,\mathfrak{a})}
\rightarrow \alpha_{f_\Phi(v,\mathfrak{a})}'.
$$
Equivalently written this precisely means for an object $(v,\mathfrak{a})$ in $G\ltimes O$
$$
\Phi_{(v,\mathfrak{a})}: \Psi(v,\mathfrak{a})\rightarrow \Psi'\circ f_\Phi(v,\mathfrak{a}).
$$
We note that $\Phi_{(0,0)}=\Phi$ and $f_\Phi(0,0)=(0,0)$. Assume that $(g,q):q\rightarrow g\ast q$ is a morphism in $G\ltimes O$ with $(g',q')=f_\Phi(g,q)$ being the corresponding morphism in $G'\ltimes O'$. 
We obtain the following commutative diagram
$$
\begin{CD}
\Psi(q) @> \Phi_{q}>> \Psi'\circ f_\Phi(q)\\
@V \Psi(g,q) VV   @V \Psi'\circ f_{\Phi}(g,q)VV\\
\Psi(g\ast q) @> \Phi_{g\ast q} >>  \Psi'\circ f_\Phi(g\ast q).
\end{CD}
$$
Hence the map
$$
\Gamma:O\rightarrow {\mathcal R}:q\rightarrow \Gamma(q)=\Phi_q
$$
is a natural transformation (in fact an equivalence)
$$
\Psi\xrightarrow{\Gamma} \Psi'\circ f_\Phi.
$$
This shows that  there is more structure than justing having a functor $F:{\mathcal R}\rightarrow \text{SET}$. Hence we have proved the following.
\begin{thm}
With the functor $F:{\mathcal R}\rightarrow \text{SET}$ as just described, given an isomorphism
$\Phi:\alpha\rightarrow \alpha'$, there exists for every element $\Psi\in F(\alpha)$ with corresponding 
element $\Psi'=F(\Phi)(\Psi)$ an uniquely determined smooth equivalence of categories 
$f_\Phi:G\ltimes O\rightarrow G'\ltimes O'$ and a natural equivalence $\Gamma:\Psi\rightarrow \Psi'\circ f_\Phi$.
\end{thm}
 At this point the good uniformizers in $F(\alpha)$ and $F(\alpha')$ are mostly unrelated with the exception
 when the objects are isomorphic. In the following, in particular for the abstract theory we shall always take the minimalistic viewpoint, even if in all applications there is additional structure along the lines as just described.
 
We  need some additional structure to related $F(\alpha)$ and $F(\alpha')$ even if $\alpha$ and $\alpha'$ are not isomorphic.  This needs some preparation.
\begin{defn}
 Given the functor ${F}:{\mathcal R}^-\rightarrow \text{SET}$ and objects $\alpha$ and $\alpha'$ in ${\mathcal R}^-$, we define for  $\Psi\in {F}(\alpha)$ and $\Psi'\in{F}(\alpha')$ the set ${\bf M}(\Psi,\Psi')$, called the associated transition set between $\Psi$ and $\Psi'$ by
 $$
 {\bf M}(\Psi,\Psi')=\{(q,\Phi,q')\ |\ q\in O,\ q'\in O',\ \Phi:\Psi(q)\rightarrow \Psi'(q')\}.
 $$
 \end{defn}
The construction of the transition set comes with several structure maps. 
 \begin{defn}
The target map $t:{\bf M}(\Psi,\Psi')\rightarrow O'$ and the source map $s:{\bf M}(\Psi,\Psi')\rightarrow O$ are given by
$$
t(q,\Phi,q')=q'\ \text{and}\ \ s(q,\Phi,q')=q.
$$
The inversion map $\iota:{\bf M}(\Psi,\Psi')\rightarrow {\bf M}(\Psi',\Psi)$  is defined by $\iota(q,\Phi,q')=(q',\Phi^{-1},q)$.
 The unit map $O\rightarrow {\bf M}(\Psi,\Psi)$ is given by  $u(q)=(q,1_{\Psi(q)},q)$,  and the multiplication map is defined as
 \begin{eqnarray*}
& {\bf M}(\Psi',\Psi''){_{s}\times_t}{\bf M}(\Psi,\Psi')\rightarrow {\bf M}(\Psi,\Psi''):&\\
&m((q',\Phi',q''),(q,\Phi,q'))=(q,\Phi'\circ\Phi,q'').&
 \end{eqnarray*}
 \end{defn}
 The following results hold for the logarithmic and the exponential gluing profile. There are similar results for other, but not all,  gluing profiles. The logarithmic case is a reformulation of the classical Deligne-Mumford theory, the case of the exponential gluing profile can be reduced to this classical case. In \cite{HWZ-DM} we derive both cases using pde-methods,  but also show that the case with exponential gluing profile can be reduced
 to the classical case.
 \begin{thm}
Given the logarithmic or the exponential gluing profile, the transition sets ${\bf M}(\Psi,\Psi')$ carry a natural metrizable topology and in addition naturally oriented smooth manifold structures, so that
 all the structure maps are smooth, and $s$ and $t$ are orientation preserving local diffeomorphisms. In case of the logarithmic gluing profile
 the smooth manifold structure is underlying to a complex manifold structure and all structure maps are holomorphic.
The  object space of $G\ltimes O$ is oriented by the orientation coming from the fact that $O$ is an open subset of a complex vector space.
The orientations of the  manifolds ${\bf M}(\Psi,\Psi')$ are determined by the fact that the source and target map are orientation preserving.
 \end{thm}
 In particular, at this point we have, after picking the exponential gluing profile $\varphi$,   a natural construction $(F,{\bf M})$ which associates to every object a set ${F}(\alpha)$ of good uniformizers $\Psi$ at $\alpha$, and to a transition set ${\bf M}(\Psi,\Psi')$  an oriented manifold structure. 

\begin{rem}
Let us just note that with $\Psi'$ corresponding to $\Psi$ if $\Phi:\alpha\rightarrow \alpha'$ is given,
the natural transformation $\Gamma:\Psi\rightarrow \Psi'\circ f_\Phi$ defines a diffeomorphism
$$
G\ltimes O\rightarrow {\bf M}(\Psi,\Psi'):(g,q)\rightarrow (q,\Gamma(g\ast q)\circ \Psi(g,q),f_\Phi(g\ast q)).
$$ 
\end{rem}
 
  Here is one thing we can do with this data. Take a family
 ${(\Psi_\lambda)}_{\lambda\in\Lambda}$ ($\Lambda$ a set), of good uniformizers so that
 $$
 |{\mathcal R}| =\bigcup_{\lambda\in\Lambda} |\Psi_\lambda(O_\lambda)|.
 $$
 Then we define an oriented smooth manifold $X$ by 
 $$
 X=\coprod_{\lambda\in\Lambda} O_\lambda,
 $$
 and a smooth oriented manifold ${\bf X}$ by
 $$
 {\bf X}=\coprod_{(\lambda,\lambda')\in\Lambda\times \Lambda} {\bf M}(\Psi_\lambda,\Psi_{\lambda'}).
 $$
We may view $X$ as the collection of objects in a small category and ${\bf X}$ as the set of morphisms. More precisely, if $q_\lambda\in O_\lambda$ and
$q'_{\lambda'}\in O_{\lambda'}$, then the morphisms $q_\lambda\rightarrow q'_{\lambda'}$ are precisely all elements $\Phi\in {\bf M}(\Psi_\lambda,\Psi_{\lambda'})$
with $s(\Phi)=q_\lambda$ and $t(\Phi)=q_{\lambda'}'$.

Using property (iv) it follows immediately
 that ${\mathcal X}=(X,{\bf X})$ is a an \'etale proper Lie groupoid, see \cite{Mj} for a short description
 of the theory. More comprehensive treatments are given in \cite{AR}, and see \cite{hae1,hae2,hae3} for the beginnings of this theory. Using  the $\Psi_\lambda$ we can construct a functor
 $$
 \beta:{\mathcal X}\rightarrow {\mathcal R}
 $$
 which is an equivalence of categories. Namely we map the object $q_\lambda\in O_\lambda\subset X$ to $\Psi_\lambda(q_\lambda)$, so that on objects 
 $$
 \beta:X\rightarrow \text{obj}({\mathcal R}):\beta(q_\lambda)=\Psi_\lambda(q_\lambda),
 $$
 and on morphisms 
 $$
 \beta:{\bf X}\rightarrow \text{mor}({\mathcal R}): \beta(q_\lambda,\phi,q_{\lambda'}')= (\Psi_\lambda(q_\lambda),\phi,\Psi_{\lambda'}(q_{\lambda'}')).
 $$

 Hence we might view ${\mathcal X} $ as a smooth (up to equivalence) model for ${\mathcal R}$.
 Of course, making different choices we obtain $\beta':X'\rightarrow {\mathcal R}$ having the same properties.
 Taking the union of the choices we obtain $\beta'':{\mathcal X}''\rightarrow {\mathcal R}$ and smooth equivalences of \'etale proper Lie groupoids 
 via the inclusions ${\mathcal X}\rightarrow {\mathcal X}''$ and ${\mathcal X}'\rightarrow {\mathcal X}''$. Taking the weak fibered product associated 
 to the diagram
 $$
 {\mathcal X}\rightarrow {\mathcal X}''\leftarrow {\mathcal X}'
 $$
 denoted by ${\mathcal X}\times_{{\mathcal X}''} {\mathcal X}'$ the projections onto the factors 
 are equivalences of \'etale proper Lie groupoids. With other words ${\mathcal X}$ and ${\mathcal X}'$ are Morita equivalent. 
 Moreover $\beta:{\mathcal X}\rightarrow {\mathcal R}$ and $\beta':{\mathcal X}'\rightarrow {\mathcal R}$ are naturally equivalent
 via the data from ${\mathcal X}''$. 
In summary we have constructed up to Morita equivalence smooth models for ${\mathcal R}$. Each pair $({\mathcal X},|\beta|)$ gives an orbifold structure
on $|{\mathcal R}|$ and these orbifold structures are equivalent. The whole collection turns $|{\mathcal R}|$ into a smooth oriented orbifold,
with the structure depending on the gluing profile $\varphi$. If the process is carried out with the logarithmic gluing profile
we again obtain a natural orbifold $|{\mathcal R}|$, which, however this time is holomorphic. The underlying smooth structures for the two constructions  mentioned 
are not the same, but they should be diffeomorphic. The latter is definitely true over the un-noded part.
The reader is referred to \cite{Mj,MM,AR} for more reading on \'etale proper Lie groupoids.

 \subsection{The Category of Stable Maps}\label{SSS2}
 The construction carried out in this subsection is formally very similar to the one we carried out for ${\mathcal R}$. Now 
 the category ${\mathcal R}$ is replaced by a different category ${\mathcal S}$, the category of stable maps. Again $|{\mathcal S}|$ has a natural metrizable topology
 and we shall construct a pair $(F,{\bf M})$ similarly as before.
However, the local models are infinite-dimensional and all ingredients have never any change to be classically smooth. In fact, the occurring maps are generally
 nowhere differentiable in any classical sense. Nevertheless, as it will turn out, there is a smoothness concept, called sc-smoothness, so that the constructions viewed within an associated differential geometry are smooth.

Let $(Q,\omega)$ be a compact  symplectic manifold  without boundary. We  consider maps defined on
Riemann surfaces with image in  $(Q,\omega)$ having various
regularity properties. We shall write
$$
u:{\mathcal O}(S,x)\rightarrow Q
$$
for a mapping germ defined on a Riemann
surface $S$ near $x$.
\begin{defn}
Let $m\geq 2$ be an integer and $\varepsilon>0$.  We say a germ of
continuous map $u:{\mathcal O}(S,x)\rightarrow Q$ is of class
$(m,\varepsilon)$ at the point $x$ if for a smooth chart
$\phi:U(u(0))\rightarrow {\mathbb R}^{2n}$ mapping $u(0)$ to $0$ and
holomorphic polar coordinates $ \sigma:[0,\infty)\times
S^1\rightarrow S\setminus\{x\}$ around $x$, the map
$$
v(s,t)=\phi\circ u\circ \sigma(s,t),
$$
defined for $s$ large, has partial derivatives up to order
$m$, which weighted by $e^{\varepsilon s}$ belong to
$L^2([s_0,\infty)\times S^1,{\mathbb R}^{2n})$ if $s_0$ is sufficiently
large. We say the germ is of class $m$ around a point $z\in S$ provided
$u$ is of class $H^m_{loc}$ near $z$.
\end{defn}
 We observe that the above
definition does not depend on the choices involved, like charts and
holomorphic polar coordinates.
We consider now tuples $\alpha=(S,j,M,D,u)=(\alpha^\ast,u)$, where $\alpha^\ast=(S,j,M,D)$ is a
noded Riemann surface with ordered marked points $M$ and nodal pairs $D$,
and $u:S\rightarrow W$ is a continuous map having some additional regularity properties. We do not assume $\alpha^\ast$ to be stable, i.e. being an object in ${\mathcal R}$ after forgetting the order of $M$. 
\begin{defn}
A {noded Riemann surface with ordered marked points} is a tuple $\alpha^\ast=(S,j,M,D)$, where
$(S,j)$ is a closed Riemann surface, $M\subset S$ a finite collection of { ordered marked points}, and $D$ is a finite collection of un-ordered pairs $\{x,y\}$ of points in $S$, called {nodal pairs}, so that $x\neq y$ and two pairs which intersect are identical. The union of all $\{x,y\}$, denoted by $|D|$ is disjoint from $M$. We call $D$ the set of nodal pairs and $|D|$ the set of nodal points.
\end{defn}
As in the definition of the objects in ${\mathcal R}$ the Riemann surface $S$ might consist of different connected components $C$. We call $C$ a domain component of $S$. The special points on $C$ are as before the points in $C\cap (M\cup|D|)$. We say that $(S,j,M,D)$ is  connected, provided the topological space $\bar{S}$ obtained by identifying $x\equiv y$ in the nodal pairs $\{x,y\}\in D$ is connected.
With our terminology it is  possible that  $\alpha^\ast=(S,j,M,D)$ is  connected  but on the other hand  $S$ may have several connected components, i.e. its domain components.

Next we describe the tuples $\alpha$ in more detail.

\begin{defn} \label{D1}We say that $\alpha=(S,j,M,D,u)$ is a stable map of {class} $({m},{\delta})$
provided the following holds, here $m\geq 2$ and $\delta>0$.
\begin{itemize}
\item[(i)] The underlying topological space obtained by identifying
the two points in any nodal pair is connected.
\item[(ii)] The map $u$ is  of class $(m,\delta)$ around the points in $|D|$ and
of class $m$ around all other points in $S$.
\item[(iii)] $u(x)=u(y)$ for every nodal pair $\{x,y\}\in D$.
\item[(iv)] If a  domain component $C$ of $S$ has genus $g_C$, and
$n_C$ special points so that $2\cdot g_C +n_C\leq 2$, then $\int_C u^\ast\omega >0$. Otherwise
we assume that $\int_C u^\ast\omega\geq 0$. This is called the stability condition.
\end{itemize}
\end{defn}

Next we introduce the category of stable maps of class $(3,\delta_0)$, where $\delta_0$ is a fixed umber in $(0,2\pi)$.
\begin{defn}
Fix a $\delta_0\in (0,2\pi)$. The category of stable maps (of class $(3,\delta_0)$), denoted by
$$
{\mathcal S}^{3,\delta_0}(Q,\omega)
$$
has as objects the tuples  $\alpha$ of class $(3,\delta_0)$ and as morphisms 
$$
\Phi:=(\alpha,\phi,\alpha'):\alpha\rightarrow \alpha',
$$
where $\phi$ is a biholomorphic map $(S,j,M,D)\rightarrow (S',j',M',D')$ satisfying $u'\circ\phi=u$.
\end{defn}
\begin{rem} (i) The category ${\mathcal S}^{3,\delta_0}(Q,\omega)$ has several interesting features.
All the morphisms are  isomorphisms, and between two objects the number of morphisms 
is finite. This is a consequence of the stability condition.

(ii) If we identify isomorphic objects in ${\mathcal S}^{3,\delta_0}(Q,\omega)$ we obtain the orbit class
$|{\mathcal S}^{3,\delta_0}(Q,\omega)|$ which is easily verified to be a set. An element in the orbit set is written as $|\alpha|$.

(iii) The choice of $\delta_0\in (0,2\pi)$ is dictated by the fact that in any analytical treatment of stable maps (Fredholm theory),
one has to derive elliptic estimates near the nodes which involve  self-adjoint operators which have spectrum $2\pi {\mathbb Z}$.
So $(0,2\pi)$ has to be understood as a spectral gap. This requirement is important for the Fredholm theory.
\end{rem}

The first important fact is that the orbit set $|{\mathcal S}^{3,\delta_0}(Q,\omega)|$ carries a natural topology. More precisely:
\begin{thm}\label{th-top}
Given $\delta_0\in (0,2\pi)$ the orbit space $|{\mathcal S}^{3,\delta_0}(Q,\omega)|$
carries a natural Hausdorff, second countable, regular, and hence metrizable, topology.
\end{thm}
We refer the reader to \cite{H2,HWZ6} for the complete construction. However, we shall give some ideas later on.
The key is a recipe to construct a basis for a topology. This recipe involves choices. Nevertheless the resulting 
topology is independent of them and this  is what we mean by being natural.
However, this is only the beginning and there is another natural construction.  What we shall see next will remind the reader 
of the type of construction which already occurred in the treatment of ${\mathcal R}$. In the case of ${\mathcal R}$ we saw uniformizers defined on translation groupoids, where object and morphism sets carried natural smooth manifold structures.  In the relevant constructions for ${\mathcal S}^{3,\delta_0}(Q,\omega)$
our translation groupoids will have much less structure (for the moment), namely they are just
metrizable topological spaces. The reason is that all occurring maps have no chance of being classically smooth and the natural local models
seem very often  far away from open subsets of Banach spaces.

Recall, that given a metrizable space $O$ with the action of a group $G$ by continuous maps,
we can construct the metrizable translation groupoid $G\ltimes O$.  This is a small metrizable category
with objects being the elements in $O$ and the morphisms being the tuples $(g,q)$, where $q\in O$ and $g\in G$.
Here the source of $(g,q)$ is $q$ and the target is $g\ast q$, i.e.
$$
(g,q):q\rightarrow g\ast q.
$$
A metrizable category is a small category, where object and morphism sets are metrizable spaces and all structure maps
are continuous, namely associating to a morphism $\phi$ its source and target, the inversion map $\phi\rightarrow \phi^{-1}$,
the unit map $q\rightarrow 1_q$, as well as the multiplications map $m(\phi,\psi)=\phi\circ \psi$ defined on the appropriate
subspace of the product of two copies of the morphism space.
\begin{defn}
Let $\alpha$ be an object in ${\mathcal S}^{3,\delta_0}(Q,\omega)$ with automorphism group $G$.
A good uniformizer around $\alpha$ is a functor $\Psi:G\ltimes O\rightarrow {\mathcal S}^{3,\delta_0}(Q,\omega)$ having the following properties.
\begin{itemize}
\item[(i)] $\Psi$ is fully faithful and there exists an object $q_0\in O$ with $\Psi(q_0)=\alpha$.
\item[(ii)] Passing to orbit spaces $|\Psi|:{_{G}\backslash }O\rightarrow |{\mathcal S}^{3,\delta_0}(Q,\omega)|$ is a homeomorphism onto an open neighborhood of $|\alpha|$ in $|{\mathcal S}^{3,\delta_0}(Q,\omega)|$.
\item[(iii)] For every object $q\in O$ there exists an open neighborhood $U(q)\subset O$, so that every sequence $(q_k)\subset U(q)$,
for which $|\Psi(q_k)|$ converges in $|{\mathcal S}^{3,\delta_0}(Q,\omega)|$, has a convergent subsequence in $\cl_O(U(q))$.
\end{itemize}
\end{defn}

Given two good uniformizers $\Psi$ around $\alpha$ and $\Psi'$ around $\alpha'$, we can define 
the set ${\bf M}(\Psi,\Psi')$ consisting of all tuples
$(q,\Phi,q')$, where $q\in O$, $q'\in O'$ and $\Phi:\Psi(\alpha)\rightarrow \Psi'(\alpha')$ is a morphism in ${\mathcal S}^{3,\delta_0}(Q,\omega)$.
Associated to ${\bf M}(\Psi,\Psi')$ we have the source map 
$$
s:{\bf M}(\Psi,\Psi')\rightarrow O:(q,\Phi,q')\rightarrow q
$$
and the target map
$$
t:{\bf M}(\Psi,\Psi')\rightarrow O':(q,\Phi,q')\rightarrow q'.
$$
In addition there is the inversion map
$$
\iota:{\bf M}(\Psi,\Psi')\rightarrow {\bf M}(\Psi',\Psi):(q,\Phi,q')\rightarrow (q',\Phi^{-1},q),
$$
the unit map
$$
u:O\rightarrow {\bf M}(\Psi,\Psi):q\rightarrow (q,1_{\Psi(q)},q),
$$
and, given a third uniformizer $\Psi''$, the multiplication map
\begin{eqnarray*}
&m:{\bf M}(\Psi',\Psi''){_{s}\times_t}{\bf M}(\Psi,\Psi')\rightarrow {\bf M}(\Psi,\Psi''):&\\
&m((q',\Phi',q''),(q,\Phi,q'))=(q,\Phi'\circ\Phi,q'').&
\end{eqnarray*}
The maps $s,t,\iota,u$ and $m$ are called the structure maps.

The next natural construction is summarized by the following theorem, where as before $\text{SET}$ is the category of sets.
\begin{thm}\label{thm1}
There exists a natural functor $F:({\mathcal S}^{3,\delta_0}(Q,\omega))^-\rightarrow \text{SET}$ which associates 
to an object $\alpha$ a set $F(\alpha)$ of good uniformizers around $\alpha$. Moreover, $F$ comes with  a natural construction, which associates to every choice $\Psi\in F(\alpha) $
and $\Psi'\in F(\alpha')$ a metrizable topology on ${\bf M}(\Psi,\Psi')$ so that the source and target maps are local homeomorphisms 
and all structure maps are continuous.
\end{thm}
\begin{rem}
It is again possible to lift the functor to ${\mathcal S}^{3,\delta_0}(Q,\omega)$ so that morphisms are mapped to bijections.
One can be  more explicit about  the correspondence of constructions for $\alpha$ and 
the constructions for $\alpha'$, when we are given a morphism $\Phi:\alpha\rightarrow \alpha'$. This is in the spirit of the discussion for ${\mathcal R}$.
\end{rem}
There are many consequences of this result. For example we can pick a family ${(\Psi_\lambda)}_{\lambda\in\Lambda}$, where $\Lambda$ is a set,
of good uniformizers so that 
$$
|{\mathcal S}^{3,\delta_0}(Q,\omega)|=\bigcup_{\lambda\in\Lambda} |\Psi(O_\lambda)|.
$$
Then we can define the objects of a category $X$ as the disjoint union 
$$
X=\coprod_{\lambda\in\Lambda} O_\lambda
$$
and the morphism set ${\bf X}$ by
$$
{\bf X} =\coprod_{(\lambda,\lambda')\in\Lambda\times\Lambda} {\bf M}(\Psi_\lambda,\Psi_{\lambda'}).
$$
We note that $X$ and ${\bf X}$ are both metrizable, and the $\Psi_\lambda$ can be used to define a functor $\beta$
from ${\mathcal X}=(X,{\bf X})$ to ${\mathcal S}^{3,\delta_0}(Q,\omega)$ by
$$
X\rightarrow \text{obj}({\mathcal S}^{3,\delta_0}(Q,\omega)):O_\lambda\ni q\rightarrow \Psi_\lambda(q)
$$
and
\begin{eqnarray*}
&{\bf X}\rightarrow \text{mor}({\mathcal S}^{3,\delta_0}(Q,\omega)):&\\
&{\bf M}(\Psi_\lambda,\Psi_{\lambda'})\ni(q,\Phi,q')\rightarrow \Phi.&
\end{eqnarray*}
\begin{lem}
The functor $\beta$ is an equivalence of categories.
\end{lem}
With other words we can build a topological version of ${\mathcal S}^{3,\delta_0}(Q,\omega)$ up to equivalence.
Of course, the construction of ${\mathcal X}$ involves a choice of a family $(\Psi_\lambda)$. As in the ${\mathcal R}$-case,
if we make a different choice we obtain
$$
\beta':{\mathcal X}\rightarrow {\mathcal S}^{3,\delta_0}(Q,\omega)
$$
and taking the disjoint union of the choices $\beta'':{\mathcal X}''\rightarrow {\mathcal S}^{3,\delta_0}(Q,\omega)$.
The inclusions ${\mathcal X}\rightarrow {\mathcal X}''$ and ${\mathcal X}'\rightarrow {\mathcal X}''$
are local homeomorphisms on the object and morphism spaces, and equivalences of categories.

Starting from this situation and quite similar to the ${\mathcal R}$-case one can define a notion of equivalence  of two pairs $({\mathcal X},\beta)$ and $({\mathcal X}',\beta')$,
so that we might say that ${\mathcal S}^{3,\delta_0}(Q,\omega)$ up to Morita equivalence has a unique topological (metrizable) model.
As it turns out the construction of $F$ in Theorem \ref{thm1} has many more properties. Namely the domains $G\ltimes O$ of the functors $\Psi_\lambda$ are in some 
sense smooth spaces. Not in the usual way, but in a generalized differential geometry. This will also be explained later on. 
At the end of day there is a differential geometric version of Theorem \ref{thm1}, which one might view as the construction of a smooth structure
on the category ${\mathcal S}^{3,\delta_0}(Q,\omega)$, and which allows to build small smooth versions of our category up to Morita equivalence.

\subsection{A Bundle and the CR-Functor}\label{SSS23}
Assume next that a compatible almost complex structure $J$ has been fixed for $(Q,\omega)$.
We consider tuples 
$$
\what{\alpha}=(\alpha,\xi)=(S,j,M,D,u,\xi),
$$
 where $\alpha$  is an object in ${\mathcal S}^{3,\delta_0}(Q,\omega)$, and $\xi(z):T_zS\rightarrow T_{u(z)}Q$ is complex anti-linear
for the given structures $j$ and $J$. Further we assume that 
$$
z\rightarrow \xi(z)
$$
 has Sobolev regularity $H^2$ away from the nodal points. At the nodal points
we assume it to be of class $(2,\delta_0)$. To make this precise,
pick a nodal point $x$ and take positive holomorphic polar coordinates around $x$, say $\sigma(s,t)$ with $x=\lim_{s\rightarrow \infty} \sigma(s,t)$.
Then take a smooth chart $\phi$ around $u(x)$ with $\phi(u(x))=0$.
Finally consider the principal part
$$
(s,t)\rightarrow pr_2\circ T\phi(u(\sigma(s,t)))\circ \xi(\sigma(s,t))\left(\frac{\partial\sigma(s,t)}{\partial s}\right),
$$
which we assume for large $s_0$ to be in $H^{2,\delta_0}([s_0,\infty)\times S^1,{\mathbb R}^{2n})$. The definition of the decay property does not depend on the choice of $\sigma$ and $\phi$.  Denote by $\text{BAN}_G$ the category where the objects are the Banach spaces and the morphisms are topological linear isomorphisms.
The above discussion gives us a functor
$$
\mu:{\mathcal S}^{3,\delta_0}(Q,\omega)\rightarrow \text{BAN}_G.
$$
We associate to an object $\alpha$ the Hilbert space of all $\xi$ of class $(2,\delta_0)$ as described above. To a morphism
$\Phi=(\alpha,\phi,\alpha')$ we associate the topological linear isomorphism
$$
\mu(\alpha)\rightarrow \mu(\alpha'):\xi\rightarrow \xi\circ T\phi^{-1}.
$$
We can now define a new category ${\mathcal E}^{2,\delta_0}(Q,\omega,J)$, whose objects are the tuples
$\what{\alpha}$ of the form
$$
\what{\alpha}=(\alpha,\xi),
$$
where $\alpha$ is an object in ${\mathcal S}^{3,\delta_0}(Q,\omega)$ and $\xi\in \mu(\alpha)$ is a complex anti-linear $TQ$-valued $(0,1)$-form
along $u$ of class $(2,\delta_0)$. A morphism 
$$
\what{\Phi}=(\what{\alpha},\Phi,\what{\alpha}'):\what{\alpha}\rightarrow \what{\alpha}'
$$
is a tuple with $\Phi=(\alpha,\phi,\alpha'):\alpha\rightarrow \alpha'$ being a morphisms in ${\mathcal S}^{3,\delta_0}(Q,\omega)$ so that in addition
$$
\mu(\Phi)(\xi)=\xi'.
$$
There is the projection functor
$$
P:{\mathcal E}^{2,\delta_0}(Q,\omega,J)\rightarrow {\mathcal S}^{3,\delta_0}(Q,\omega),
$$
which on objects maps $(\alpha,\xi)\rightarrow \alpha$ and on morphisms  $((\alpha,\xi),\Phi,(\alpha',\xi'))$ to $\Phi$.
On the object level the preimage under $P$ of an object $\alpha$  is a Hilbert space consisting of the elements $(\alpha,\xi)$. Given $\Phi:\alpha\rightarrow \alpha'$
the morphism $\Phi$ is lifted to a bounded linear isomorphism
$$
\boldsymbol{\Phi}:P^{-1}(\alpha)\rightarrow P^{-1}(\alpha'): (\alpha,\xi) \rightarrow (\alpha',\xi\circ T\phi^{-1}),
$$
where the map $\xi\rightarrow \xi\circ T\phi^{-1}$ is linear.
We may view the two categories  fibering over each other via $P:{\mathcal E}^{2,\delta_0}(Q,\omega,J)\rightarrow {\mathcal S}^{2,\delta_0}(Q,\omega)$ as a `bundle'. Then we have the Cauchy-Riemann section functor
$$
\bar{\partial}_J(\alpha)=\left(\alpha,\frac{1}{2}\cdot\left[ Tu+J(u)\circ Tu\circ j\right]\right),
$$
where $\alpha=(S,j,M,D,u)$. From the results in \cite{HWZ6} we can deduce the following theorem.
\begin{thm}\label{thm2}
The orbit space $|{\mathcal E}^{2,\delta_0}(Q,\omega,J)|$ carries a natural second countable metrizable topology.
The induced maps on orbit spaces $|P|$ and $|\bar{\partial}_J|$ are continuous. The topology on the orbit space
of  the full subcategory associated to the objects $\alpha$ with $\bar{\partial}_J(\alpha)=0$ has compact connected components
and coincides with the (union of all) Gromov compactified moduli space of $J$-pseudoholomorphic curves.
\end{thm}
The previous discussion about the construction $F$ can be extended to cover $P$. Given an  object $\alpha$ the notion of a good uniformizer
generalizes by replacing $O$ by 
a continuous surjective map $p:W\rightarrow O$ between metrizable spaces,  where the fibers are Hilbert spaces (or more generally Banach spaces),
and where in addition we have a continuous $G$-action, linear between the fibers. Then the good bundle uniformizers are given by commutative diagrams
$$
\begin{CD}
G\ltimes W @> \bar{\Psi}>> {\mathcal E}^{2,\delta_0}(Q,\omega,J)\\
@VVV @VVV\\
G\ltimes O @>\Psi>> {\mathcal S}^{3,\delta_0}(Q,\omega)
\end{CD}
$$
with the obvious properties.
The bottom $\Psi$ is as described before and $\bar{\Psi}$ is a linear bounded isomorphism between the fibers. 
Again, as it will turn out,  these constructions will fit into a scheme of generalized differential geometry, and the local representative
of $\bar{\partial}_J$, a $G$-equivariant section of $W\rightarrow O$,  will turn out to be a nonlinear Fredholm section in the extended
framework.  We refer the reader to \cite{FH,HWZ6} for the constructions related to stable maps, and to \cite{HWZ7,HWZ8} for the abstract theory.
In the above extension we have a functorial construction $\bar{F}:({\mathcal S}^{3,\delta_0}(Q,\omega))^-\rightarrow \text{SET}$ covering $F$
in the following way. For every object $\alpha$ the set $\bar{F}(\alpha)$ consists of strong bundle uniformizers $\bar{\Psi}$,
which together with the underlying $\Psi\in F(\alpha)$, fit into the above commutative diagram. Moreover, it comes with an associated construction
of a metrizable topology on the transition set
${\bf M}(\bar{\Psi},\bar{\Psi}')$. There is a projection
$$
{\bf p}:{\bf M}(\bar{\Psi},\bar{\Psi}')\rightarrow{\bf M}({\Psi},{\Psi}'):(k,\what{\Phi},k')\rightarrow (o,\Phi,o')
$$
whose fibers are Hilbert spaces. Moreover, there exists the commutative diagram involving source and target maps
$$
\begin{CD}
K @<s<<  {\bf M}(\bar{\Psi},\bar{\Psi}') @>t>> K\\
@V p VV @V {\bf p} VV  @ V p' VV\\
O@< s<<  {\bf M}(\Psi,\Psi') @> t >> O'.
\end{CD}
$$
We leave the execution of the idea to the reader. The above follows from the results in \cite{HWZ6}.

\subsection{A Basic Construction}
The following discussion is carried out in detail in \cite{HWZ8.7}, see also \cite{HWZ6}.
Consider the Hilbert space $E$ consisting of pair of maps $(u^+,u^-)$
$$
u^\pm:{\mathbb R}^\pm\times S^1\rightarrow {\mathbb R}^N,
$$
where for each pair
 there exists a constant $c\in {\mathbb R}^N$, so that $u^\pm-c$ has partial derivatives
up to order $3$ which weighted by $e^{\delta_0|s|}$ belong to $L^2({\mathbb R}^\pm\times S^1)$. The constant $c$ is called the common asymptotic limit
of $(u^+,u^-)$. 

Given a complex number
$|a|<1/2$ we write $a=|a|\cdot e^{-2\pi i\theta}$. If $a\neq 0$ we define $R=\varphi(|a|)$, where
$\varphi$ is the  exponential gluing profile, which we recall, is defined by
$$
\varphi:(0,1]\rightarrow [0,\infty): \varphi(r)=e^{\frac{1}{r}}-e.
$$
We construct for $0<|a|<1/2$ the finite cylinder $Z_a$ by identifying $(s,t)\in [0,R]\times S^1$ with $(s',t')\in [-R,0]\times S^1$ via
$$
s=s'+R\ \text{and}\  t=t'+\theta. 
$$
Denote by $[s,t]$ the equivalence class of the point $(s,t)\in [0,R]\times S^1$ and by $[s',t']'$ the equivalence class of $(s',t')\in [-R,0]\times S^1$.
Then
$$
[s,t]=[s-R,t-\theta]'
$$
and $Z_a\rightarrow [0,R]\times S^1$ and $Z_a\rightarrow [-R,0]\times S^1$ defined by $[s,t]\rightarrow (s,t)$ and $[s',t']'\rightarrow (s',t')$, respectively,
are holomorphic coordinates, where the targets are equipped with the standard complex structures.

If $a=0$ we define $Z_0={\mathbb R}^+\times S^1\coprod {\mathbb R}^-\times S^1$.
Pick a smooth map $\beta:{\mathbb R}\rightarrow [0,1]$ satisfying $\beta(s)=1$ for $s\leq -1$, $\beta'(s)<0$ for $s\in (-1,1)$ and
$\beta(s)+\beta(-s) =1$ for all $s\in {\mathbb R}$. 
\begin{defn}
Given $(u^+,u^-)\in E$ and $|a|<1/2$ the glued map 
$$
\oplus_a(u^+,u^-):Z_a\rightarrow {\mathbb R}^N,
$$
is defined for $a=0$ by $\oplus_0(u^+,u^-)=(u^+,u^-)$ and for $a\neq 0$ by
$$
\oplus_a(u^+,u^-)([s,t])=\beta\left(s-\frac{R}{2}\right) u^+(s,t) +\left(1-\beta\left(s-\frac{R}{2}\right)\right) u^-(s-R,t-\theta).
$$
\end{defn}
For fixed $a\neq 0$ many different $(u^+,u^-)$ will produce the same glued map. We can remove the ambiguity by introducing the anti-glued map
$\ominus_a(u^+,u^-)$. It is defined on $C_a$ given by $C_0=\emptyset$ and for $0<|a|<1/2$ by
gluing ${\mathbb R}^-\times S^1$ and ${\mathbb R}^+\times S^1$ along $Z_a$,. Note that $Z_a$  can be identified
with an obvious subset of both half cylinders. Namely $C_a$ is given by
$$
C_a=(({\mathbb R}^-\times S^1)\coprod ({\mathbb R}\times S^1))/\sim,
$$
where $(s',t')\equiv (s,t)$ provided $(s',t')\in [-R,0]\times S^1$, $(s,t)\in {\mathbb R}^+\times S^1$
and $s=s'+R$ and $t=t'+\theta$.  We shall write $[s,t]$ or equivalently $[s',t']'$ for the points
in the part being identified. Note that we have for $0<|a|<1/2$  a natural embedding
$$
Z_a\rightarrow C_a:[s,t]\rightarrow [s,t].
$$
For $a\neq 0$ the cylinder $C_a$ has a natural complex manifold structure,
for which the maps on $Z_a\subset C_a$ defined by
$$
[s,t]\rightarrow (s,t)\in {\mathbb R}\times S^1,
$$
and
$$
[s',t']'\rightarrow (s',t')\in {\mathbb R}\times S^1
$$
have natural extensions to biholomorphic maps $C_a\rightarrow {\mathbb R}\times S^1$.
The extensions will be written as above. The transition map
$$
{\mathbb R}\times S^1\rightarrow {\mathbb R}\times S^1:(s,t)\rightarrow (s',t')
$$
is given by $(s,t)\rightarrow (s-R,t-\theta)$.
\begin{defn}
For $(u^+,u^-)\in E$ the anti-glued map $\ominus_a(u^+,u^-):C_a\rightarrow {\mathbb R}^N$ is defined for $a=0$ by 
$\ominus_0(u^+,u^-)=0$, i.e. the only map $\emptyset\rightarrow {\mathbb R}^N$, and otherwise
by
\begin{eqnarray*}
\ominus_a(u^+,u^-)([s,t])
&=&\left(\beta\left(s-\frac{R}{2}\right)-1\right) \left( u^+(s,t)-\text{av}_a(u^+,u^-)\right)\\
&&+\beta\left(s-\frac{R}{2}\right)\left(u^-(s-R,t-\theta)-\text{av}_a(u^+,u^-)\right).
\end{eqnarray*}
Here 
$$
\text{av}_a(u^+,u^-)=\frac{1}{2}\cdot \int_{S^1} (u^+(R/2,t)+u^-(-R/2,t)) dt.
$$
\end{defn}
If we start with $(u^+,u^-)\in E$ and pick any $0<|a|<1/2$ the total gluing map
$$
\boxdot_a(u^+,u^-)=(\oplus_a(u^+,u^-),\ominus_a(u^+,u^-)),
$$
defines a bounded linear isomorphism 
from $E$ to $H^3(Z_a,{\mathbb R}^N)\oplus H^{3,\delta_0}_c(C_a,{\mathbb R}^N)$. Here $ H^{3,\delta_0}_c(C_a,{\mathbb R}^N)$ denotes  the Hilbert space
of maps $u:C_a\rightarrow {\mathbb R}^N$ which are of class $H^3_{loc}$ so that in addition there exists a constant $c\in {\mathbb R}^N$ 
depending on $u$ for which $u\pm c$ has partial derivatives up to order $3$ which weighted by $e^{\delta_0|s|}$ belong to $L^2([s_0,\infty)\times S^1)$ or
$L^2(-\infty,-s_0]\times S^1)$, respectively, for a suitably large $s_0$.
In particular we obtain for $a\neq 0$ a linear topological isomorphism
$$
\ker(\ominus_a)\rightarrow H^3(Z_a,{\mathbb R}^N):(u^-,u^+)\rightarrow \oplus_a(u^+,u^-).
$$
Since $E$ decomposes as the topological direct sum 
$$
E=\ker(\ominus_a)\oplus \ker(\oplus_a)
$$
 there exists an associated
family $a\rightarrow \pi_a$ of bounded projection operators $\pi_a:E\rightarrow E$ which project onto
$\ker(\ominus_a)$ along $\ker(\oplus_a)$.  
\begin{prop}\label{prop-cont}
For every $|a|<1/2$ the linear operator $\pi_a:E\rightarrow E$ is a bounded projection. However, the map
$$
\{a\in {\mathbb C}\ |\ |a|<1/2\}\rightarrow L(E):a\rightarrow \pi_a,
$$
where the space of bounded operators $L(E)$ is equipped with the operator norm is nowhere continuous.
The map
$$
\{a\in {\mathbb C}\ |\ |a|<1/2\}\times E\rightarrow E:(a,(u^+,u^-))\rightarrow \pi_a(u^+,u^-)
$$
is continuous.
\end{prop}
The map $r:\{|a|<1/2\}\times E\rightarrow \{|a|<1/2\}\times E$ defined by
$$
r(a,(u^+,u^-))=(a,\pi_a(u^+,u^-))
$$
is continuous and satisfies $r\circ r=r$. It also preserves the fibers of $\{|a|<1/2\}\times E\rightarrow \{|a|<1/2\}$.  The image of $r$ 
is then a fiber-wise retract. Define $\ker(\ominus_.)$ as the subset of $\{a\in {\mathbb C}\ |\ |a|<1/2\}\times E$ consisting of all tuples
$(a,(u^+,u^-))$ with $\ominus_a(u^+,u^-)=0$. Then we have continuous projection
$$
\ker(\ominus_.)\rightarrow \{a\in {\mathbb C}\ |\ |a|<1/2\},
$$
where the fibers are sc-Hilbert spaces. Let us define the set $\bar{X}^{3,\delta_0}({\mathbb R}^N)$ by
$$
\bar{X}^{3,\delta_0}({\mathbb R}^N):=(\{0\}\times E)\bigcup\left(\bigcup_{0<|a|<1/2} \left(\{0\}\times H^3(Z_a,{\mathbb R}^N)\right)\right).
$$
We obtain a fiber-preserving, fiber-wise linear, bijection
$$
\begin{CD}
\ker(\ominus_.) @>\oplus_. >> \bar{X}^{3,\delta_0}({\mathbb R}^N)\\
@VVV @VVV\\
\{|a|<1/2\} @= \{|a|<1/2\}.
\end{CD}
$$
We equip the right-hand side with the topology ${\mathcal T}$ which makes $\oplus_.$ a homeomorphism. 
We observe that $\ker(\ominus_.)$ depends on the choice of the cut-off function $\beta$, whereas the right-hand side does not depend
on any choice. So one might ask if ${\mathcal T}$ depends on the choice of $\beta$.
\begin{prop}
The definition of the topology ${\mathcal T}$ does not depend on the choice of the cut-off function $\beta$.
Moreover ${\mathcal T}$ is a metrizable topology.
\end{prop}
As a consequence of this proposition  $\bar{X}^{3,\delta_0}({\mathbb R}^N)$ is  in a natural way a metrizable topological space so that the projection
$$
p:\bar{X}^{3,\delta_0}({\mathbb R}^N)\rightarrow \{|a|<1/2\}:[u:Z_a\rightarrow {\mathbb R}^N]\rightarrow a,
$$
is continuous. We can carry out this construction for every ${\mathbb R}^N$. Interestingly a smooth map $f:{\mathbb R}^N\rightarrow {\mathbb R}^M$
induces a continuous map 
$$
\bar{X}^{3,\delta_0}({\mathbb R}^N)\rightarrow \bar{X}^{3,\delta_0}({\mathbb R}^M):u\rightarrow f\circ u.
$$
Summarizing this means that the following assertion holds true.
\begin{prop}
The construction $\bar{X}^{3,\delta_0}$ defines a natural functor from the category of Euclidean spaces with objects ${\mathbb R}^N$
and smooth maps between them, into the category of metrizable topological spaces.
\end{prop}
If $Q$ is a smooth connected manifold without boundary  we can take for a sufficiently large $N$ a smooth proper embedding
$j:Q\rightarrow {\mathbb R}^N$ and define $\bar{X}^{3,\delta_0}(Q)$ to consist of all continuous maps $u:Z_a\rightarrow Q$ such that
$j\circ u\in X^{3,\delta_0}({\mathbb R}^N)$, and equip it with the topology making the map $u\rightarrow j\circ u$ a homeomorphism
onto its image. The resulting topological space $\bar{X}^{3,\delta_0}(Q)$ does not depend on the choice of the embedding.
Indeed, if $k:Q\rightarrow {\mathbb R}^M$ is another proper embedding, then $k\circ j^{-1}$ and $j\circ k^{-1}$ are restrictions of globally 
defined smooth maps, which implies the conclusion. If $Q$ does not have boundary but is not connected we can apply the above to every 
component. Consequently we have the following result.
\begin{prop}
The functor $\bar{X}^{3,\delta_0}$, defined for the spaces ${\mathbb R}^N$ and the smooth maps between them, has a natural extension
to the category of smooth manifolds without boundaries and the smooth maps between them.
\end{prop}
Given a finite-dimensional smooth manifold $\Lambda$ and a smooth map 
$$
\Lambda:L\rightarrow \{a\in {\mathbb C}\ |\ |a|<1/2\}
$$
we can consider the pull-back $\Lambda^\ast X^{3,\delta_0}({\mathbb R}^N)\subset \Lambda\times X^{3,\delta_0}({\mathbb R}^N)$
defined by the diagram
$$
\begin{CD}
@. X^{3,\delta_0}({\mathbb R}^N)\\
@.  @VVV\\
L @>\Lambda>> \{|a|<1/2\}.
\end{CD}
$$
The next result is more complicated to prove. Denote by $j$ the standard almost complex structure on $Z_a$. We assume that for $|a|<1/2$ we are given a family 
$$
a\rightarrow j(a)
$$
of smooth almost complex structures on $Z_a$.  The $j(a)$ live on different domains and we require a certain form of smooth dependence on $a$.  The first requirement is the following where $\varphi$ denotes the exponential gluing profile.
\begin{itemize}
\item[(i)] There exist  $\varepsilon\in (0,1/2)$ and $H>0$ so that 
$$
j(a)=j\ \text{on}\  \{[s,t]\in Z_a\ |\ s\in [H,\varphi(|a|)-H]\}
$$
provided $0<|a|<\varepsilon$, and if $a=0$  it holds that $j(0)|({\mathbb R}^\pm\times S^1)$
equals $j$ on $[H,\infty)\times S^1$ and $(-\infty,-H]\times S^1$. Here we assume that $\varphi(\varepsilon)> 2H$. 
\end{itemize}
With other words, if the neck gets longer the structure will be $j$ in the necks, but can be different near the boundaries.
In the next step we shall have to formulate the smoothness requirements near the boundaries. Consider for $|a|<\varepsilon$
the maps 
$$
[0,\frac{3}{2}H]\times S^1\rightarrow Z_a:(s,t)\rightarrow [s,t]
 \ \text{and}\ \ [-\frac{3}{2}H,0]\times S^1\rightarrow Z_a:(s',t')\rightarrow [s',t']',
 $$
 which we can use to pull back $j(a)$ to obtain the families $j^+(a)$ and $j^-(a)$, respectively defined for $|a|<\varepsilon$.
\begin{itemize}
\item[(ii)] For $\varepsilon$ as in (i) the maps $a\rightarrow j^\pm(a)$, $|a|<\varepsilon$ are smooth.
\end{itemize}
Next we need to say something about the smoothness for $|a|<1/2$. If $\frac{1}{2}\varepsilon<|a|< 1/2$ we can identify $[0,1]\times S^1$
with $Z_a$ via the map $(s,t)\rightarrow [\varphi(|a|)\cdot s,t]$ and via pull-back
obtain the family $j^+(a)$, $\frac{1}{2}\varepsilon<|a|<1/2$ on $[0,1]\times S^1$.
Similarly we can identify $[-1,0]\times S^1$ with $Z_a$ via the map
$(s',t')\rightarrow [-\varphi(|a_0|)s',t']'$ and obtain the pull-back family $a\rightarrow j^-(a)$ on
$[-1,0]\times S^1$.
We require the following.
\begin{itemize}
\item[(iii)] $a\rightarrow j^+(a)$ and $a\rightarrow j^-(a)$ for $\frac{1}{2}\varepsilon<|a|<1/2$ are smooth families.
\end{itemize}
\begin{defn}
A family $a\rightarrow j(a)$, which associates to $|a|<1/2$ a smooth almost complex structure on $Z_a$
and  having the properties (i),(ii), and (iii) as described above will be called a good family (of smooth almost complex structures).
\end{defn}
There are many  equivalent formulations for defining a good family.
The important result is concerned with the continuity of a certain class of maps.
\begin{thm}\label{prop1x}
Assume $j$ and $k$ are two good families as described above and 
$$
\{|a|<1/2\}\rightarrow \{|a|<1/2\}: a\rightarrow b(a)
$$
is a smooth map so that there exists  a core-smooth family (defined below)
$$
a\rightarrow [ \phi_a:(Z_a,j(a))\rightarrow (Z_{b(a)},k(b(a))) ]
 $$
  of biholomorphic maps. Then the map
$$
b^\ast\bar{X}^{3,\delta_0}(Q)\rightarrow \bar{X}^{3,\delta_0}(Q):(a,u)\rightarrow u\circ \phi_a
$$
is continuous.
\end{thm}
The family $a\rightarrow\phi_a$ being core-smooth just requires this map and its inverse to be smooth near the boundaries.
This is well-defined, using  the coordinates $[s,t]$ and $[s',t']'$.  By unique continuation $\phi_a$ is entirely determined  by knowing
$a$, $b(a)$, and its values near the boundaries. So the fact that $\phi_a:Z_a\rightarrow Z_{b(a)}$ is biholomorphic also determines 
the behavior inside, which, of course, is crucial for the validity of the theorem.
\begin{rem}
The proof of this result requires some work. The above map will in fact be sc-smooth once we have introduced this concept.
In \cite{HWZ8.7,HWZ6} proofs of sc-smooth versions of Theorem \ref{prop1x} can be found. 
We note that maps of the type occurring in Theorem \ref{prop1x} usually will not be classically differentiable unless they are constant.
In fact the expression $u\circ\phi_a$, when differentiating with respect to $a$ would loose a derivative of $u$.
\end{rem}

The previous discussion can be extended. In a first step, rather than considering the half-cylinders ${\mathbb R}^\pm\times S^1$,
we consider a noded disk pair $(D_x\cup D_y,\{x,y\})$. From all possible decorations $[\what{x},\what{y}]$ of $\{x,y\}$
we obtain the natural gluing parameters $r[\what{x},\what{y}]$ with $r\in [0,1/2)$ and define the set of gluing parameters associated to $\{x,y\}\in D$ as described before by 
$$
{\mathbb B}^{\{x,y\}}(1/2)=\{r[\what{x},\what{y}]\ |\ r\in [0,1/2),\ [\what{x},\what{y}]\ \text{decorated nodal pair}\}.
$$
 If we pick $\what{x}$ and $\what{y}$, so that $\{\what{x},\what{y}\}$ is a representative of $[\what{x},\what{y}]$,  there are (unique) associated biholomorphic maps $h_x:(D_x,x)\rightarrow ({\bf D},0)$ and similarly $h_y$ so that $Th_x(x)\what{x}={\mathbb R}$ and $Th_y(y)\what{y}={\mathbb R}$. With this data we first obtain the  map
$$
\phi:{\mathbb B}^{\{x,y\}}(1/2)\rightarrow \{z\in {\mathbb C}\ |\ |z|<1/2\}: r[\what{x},\theta\what{y}]\rightarrow r\theta,
$$
which by the definition of the complex structure on ${\mathbb B}^{\{x,y\}}(1/2)$ is biholomorphic.
If $\mathfrak{a}=[\what{x},\theta\what{y}]$ and $h_x$ and $h_y$ are as described above, then $h_x$ and $\theta^{-1}\cdot h_y$
is the uniquely determined choice of maps $(D_x,x)\rightarrow ({\mathbb D},0)$ and $(D_y,y)\rightarrow ({\mathbb D},0)$
with $Th_x(x)\what{x}={\mathbb R}$ and $T(\theta^{-1}\cdot h_y)(y)(\theta\cdot \what{y})={\mathbb R}$.
Hence the identification $z\equiv z'$ with $h_x(z)\cdot h_y(z') =\theta\cdot e^{-2\pi R}$ defines $Z_{r[\what{x},\theta\what{y}]}$, where $R=\varphi(r)$.
With $h_x$ and $h_y$ fixed, let us define biholomorphic maps 
$$
D_x\setminus\{x\}\rightarrow [0,\infty)\times S^1:z\rightarrow (s,t),\ \text{where}\  e^{-2\pi(s+it)}=h_x(z),
$$
and 
$$
D_y\setminus\{y\}\rightarrow (-\infty,0]\times S^1:z'\rightarrow (s',t'),\ \text{where}\ e^{2\pi(s'+it')}=h_y(z').
$$
We would like to find the complex number $a=a(\mathfrak{a})$ so that we obtain a natural induced biholomorphic map from ${\mathcal Z}_\mathfrak{a}$ to
$Z_a$. Clearly the identification $z\equiv z'$ associated to $h_x(z)\cdot (\theta^{-1}\cdot h_y(z'))=e^{-2\pi R}$, which defines ${\mathcal Z}_\mathfrak{a}$
if $\mathfrak{a}=r[\what{x},\theta\cdot \what{y}]$, is equivalent to the identification $(s,t)\equiv (s',t')$ given by
$$
e^{-2\pi(s+it)}\cdot e^{2\pi(s'+it')} =e^{-2\pi i\vartheta}\cdot e^{-2\pi R},
$$
where we write $\theta$ as $e^{-2\pi i\vartheta}$. Hence 
$$
s=s'+R\ \text{and}\ \ t=t'+\vartheta.
$$
This is the identification defining $Z_a$. 
Hence with $a(r[\what{x},\theta\what{y}])= r\theta$ we obtain an induced biholomorphic map
$$
{\mathcal Z}_\mathfrak{a}\rightarrow Z_{a(\mathfrak{a})}: [z=e^{-2\pi(s+it)}]\rightarrow [s,t].
$$
This family of maps allows us to use the definition of $\bar{X}^{3,\delta_0}({\mathbb R}^N)$ to define a set of maps
on the ${\mathcal Z}_\mathfrak{a}$ resulting in a set $\bar{X}^{3,\delta_0}_{\mathcal D}({\mathbb R}^N)$ with 
$$
{\mathcal D}=(D_x\cup D_y,\{x,y\}).
$$
More precisely
$$
\bar{X}^{3,\delta_0}_{\mathcal D}({\mathbb R}^N)=H^{3,\delta_0}_c({\mathcal D},{\mathbb R}^N)\coprod\left[\coprod_{\mathfrak{a}\in {\mathbb B}^{\{x,y\}}(1/2)\ |\ |0<\mathfrak{a}|<1/2\}} H^3(Z_\mathfrak{a},{\mathbb R}^N)\right].
$$
The elements in $ H^{3,\delta_0}_c({\mathcal D},{\mathbb R}^N)$ are pairs $(u^+,u^-)$ of continuous maps defined $D_x$ and $D_y$, respectively, so that $u^+(x)=u^-(y)$ and if $\sigma_x$ are positive polar coordinates on $D_x$ centered at  $x$, and $\sigma_y$ negative ones on $D_y$ centered at $y$, the pair
$(u^+\circ\sigma_x,u^-\circ \sigma_y)$ belongs to $E$.
The previous construction might a priori depend on the choices involved. However making different initial choices of identifying
$\bar{X}^{3,\delta_0}_{\mathcal D}({\mathbb R}^N)$ with $E$ the associated transition families
induce continuous maps 
as a consequence of Theorem \ref{prop1x}. 
 This is a natural construction and the space itself is homeomorphic to the image of a continuous retraction.
Moreover, it follows  immediately from the definition that for a smooth map $f:{\mathbb R}^N\rightarrow {\mathbb R}^M$ the map $w\rightarrow f\circ w$ is continuous. From this we can deduce that the functor $\bar{X}^{3,\delta_0}_{\mathcal D}$ 
has a natural extension to all smooth manifolds $Q$ without boundary. 

At this point we have a well-defined functorial construction which associates 
to a nodal disk pair ${\mathcal D}=(D_x\cup D_y,\{x,y\})$ and smooth manifold without boundary $Q$ a metrizable topological space $\bar{X}^{3,\delta_0}_{\mathcal D}(Q)$ and to a smooth map $f:Q\rightarrow Q'$ a continuous map
$$
\bar{X}^{3,\delta_0}_{\mathcal D}(Q)\rightarrow \bar{X}^{3,\delta_0}_{\mathcal D}(Q'):w\rightarrow f\circ w.
$$
We can go one step further. This step is natural as well.
Consider a noded Riemann surface $(S,j,D)$. We assume that $S$ has no boundary and is compact. We do not assume any stability assumption.
Suppose that we are given a finite group $G$ acting by biholomorphic maps on $(S,j,D)$ and ${\bf D}$ is a small disk structure invariant under $G$. Then 
we can take the natural gluing parameters and consider the family $\mathfrak{a}\rightarrow (S_\mathfrak{a},j_\mathfrak{a},D_\mathfrak{a})$.
We have dealt with problems arising near the nodes and the obtained glued necks. Away from these regions we are just dealing with `static' domains and Sobolev class  $H^3$ functions on them. We obtain the following result from the previous discussion.
\begin{thm}
Given a noded closed Riemann surface $(S,j,D)$ with an action by biholomorphic maps by a finite group and a $G$-invariant small disk structure,
there exists a functorial construction of a metrizable topological space $\bar{X}^{3,\delta_0}_{(S,j,D),{\bf D}}(Q)$ for every smooth manifold $Q$ without boundary.
Moreover $G$ acts continuously on these spaces.
Associated  to a smooth map $f:Q\rightarrow Q'$ there is  a continuous equivariant map
$$
\bar{X}^{3,\delta_0}_{(S,j,D),{\bf D}}(Q)\rightarrow \bar{X}^{3,\delta_0}_{(S,j,D),{\bf D}}(Q'):w\rightarrow f\circ w.
$$
\end{thm}
There are some additional useful constructions which will be used later. Assume that $\Xi\subset S$ is a finite $G$-invariant subset 
in the complement of the disks coming from the small disk structure.  Suppose that we have fixed for every $G$-orbit $[z]$ of a point 
$z\in \Xi$ a codimension two submanifold $H_{[z]}$ in $Q$. Denote the collection of all these constraints by ${\mathcal H}$. Define 
a subset $\bar{X}^{3,\delta_0}_{(S,j,D),{\bf D},{\mathcal H}}(Q)$ of $\bar{X}^{3,\delta_0}_{(S,j,D),{\bf D}}(Q)$ to consist of all $u$ intersecting $H_{[z]}$ at $z$ transversally for all $z\in \Xi$. If we have a map $u$ in $\bar{X}^{3,\delta_0}_{(S,j,D),{\bf D},{\mathcal H}}(Q)$ we find a small neighborhood around 
every $z\in\Xi$ and an open neighborhood $U$ around $u$ so that for every $z\in\Xi$ there is continuous map $U\ni w\rightarrow z_w$, uniquely determined
by the requirement 
 that $w$ intersects $H_{[z]}$ transversally at $z_w$.  Denote by $\Xi'$ a deformation of $\Xi$.
 We can construct  a continuous family of self-maps $\phi_{\Xi',\mathfrak{a}}:(S_\mathfrak{a},D_\mathfrak{a})\rightarrow (S_\mathfrak{a},D_\mathfrak{a})$,
 which in the necks are the identity and which move $\Xi'$ back to $\Xi$ and are supported in a small neighborhood around $\Xi$.. 
 Then the map defined for $w\in U$ given by
 $$
 w\rightarrow w\circ \phi_{\Xi_w,\mathfrak{a}}^{-1},
 $$
 where $\Xi_w$ is the deformation associated to $w$, i.e. the collection of all $z_w$, defines a continuous retraction on some open neighborhood
 $O=O(w)$ so that $r(O) =O\cap \bar{X}^{3,\delta_0}_{(S,j,D),{\bf D},{\mathcal H}}(Q)$. With other words the subset
 $ \bar{X}^{3,\delta_0}_{(S,j,D),{\bf D},{\mathcal H}}(Q)$ of $\bar{X}^{3,\delta_0}_{(S,j,D),{\bf D}}(Q)$ is locally a retract given by an retraction 
 of a particular form. Since the ambient space was modeled locally on  retracts the same is true for  this subset.
 \begin{rem}
All the constructions which we have carried out do not have a chance to be classically smooth. However, as we shall see later, they are in fact smooth construction 
 in the sc-smooth world, i.e. they are differential geometric constructions in an extended differential geometry as described in \cite{HWZ7,HWZ8}.
\end{rem}
We complete the discussion by adding some words about the construction of certain bundles over $\bar{X}^{3,\delta_0}_{(S,j,D),{\bf D},{\mathcal H}}(Q)$.
First of all we would like to change the almost complex structure on the $(S_\mathfrak{a},D_\mathfrak{a})$ by deforming the $j_\mathfrak{a}$.
Let us assume for the moment we are given a smooth family 
$$
\mathfrak{j}:V\ni v\rightarrow j(v)
$$
satisfying  $j(0)=j$ and $j(v)=j$ on the disks of the small disk structure.
Here $V$ is an open subset of some finite-dimensional vector space.  In general we would have an action of $G$ on these spaces,
but at the moment this is not relevant for the point we would like to make.
Then $V\times \bar{X}^{3,\delta_0}_{(S,j,D),{\bf D},{\mathcal H}}(Q)$ is a metrizable topological space and we can view
a point $(v,w)$ in this space as a map defined on $(S_\mathfrak{a},j(v)_\mathfrak{a},D_\mathfrak{a})$, i.e. defined on a deformation of
$(S_\mathfrak{a},j_\mathfrak{a},D_\mathfrak{a})$, the space on which  $w$ was originally defined. We consider next tuples $(v,u,\xi)$,
where $z\rightarrow \xi(z)$ is of class $H^{2,\delta_0}$ and for fixed $z\in S_\mathfrak{a}$ it holds that 
$$
\xi(z):(T_zS_{\mathfrak{a}},j(v)_{\mathfrak{a}})\rightarrow (T_{u(z)}Q,J)
$$
is complex anti-linear. By $\bar{E}^{2,\delta_0}_{(S,j,D),{\bf D},{\mathcal H},\mathfrak{j}}(Q,\omega,J)$ we denote the collection of all these tuples.
We have a natural projection
$$
p:\bar{E}^{2,\delta_0}_{(S,j,D),{\bf D},{\mathcal H},\mathfrak{j}}(Q,\omega,J)\rightarrow V\times \bar{X}^{3,\delta_0}_{(S,j,D),{\bf D},{\mathcal H}}(Q):(v,u,\xi)\rightarrow (v,u).
$$
A similar construction involving retractions allows to equip the above with a natural metrizable topology and the local models 
for these spaces are $K\rightarrow O$, where $O$ is a retraction in an Hilbert  space, and similarly for $K$ in a product of Hilbert spaces,
where in this case the retraction is linear on fibers, with respect to a projection of the product onto the first factor. We refer the reader to \cite{HWZ6,HWZ8.7} for details.

\subsection{A Good Uniformizer}
In the following discussion we shall  need the version of Deligne-Mumford theory of stable Riemann surfaces using the exponential gluing profile $\varphi$.
The discussion is in spirit completely parallel  to the one for  the category ${\mathcal R}$. However, the difference is that all occurring maps in the discussion of ${\mathcal R}$ are either smooth or holomorphic (depending on the gluing profile), but in the stable map discussion they are only continuous. 

We shall now describe the construction of the functors $\Psi$ which were introduced in Subsection \ref{SSS2}.
 Fix an object $\alpha$ in ${\mathcal S}^{3,\delta_0}(Q,\omega)$ with automorphism group $G$ 
and write $\alpha=(S,j,M,D,u)$. As a consequence of the stability condition we obtain the following.
\begin{lem}
The group  $G$ is a finite.
\end{lem}
We need the notion of a domain stabilization.
\begin{defn}
A domain stabilization $\Xi$ for $\alpha$ consists of a finite subset $\Xi\subset S\setminus(M\cup|D|)$ so that the following holds.
\begin{itemize}
\item[(i)] $(S,j,M^\ast,D)$ with $M^\ast$ being the unordered set $M\cup\Xi$ is a stable, perhaps noded,  Riemann surface.
\item[(ii)] $\Xi$ is an invariant set under the action of $G$.
\item[(iii)] If $z,z'\in \Xi$ and $u(z)=u(z')$, then there exists a $g\in G$ with $g(z)=z'$.
\item[(iv)] For every $z\in \Xi$ the map $Tu(z):T_zS\rightarrow T_{u(z)}Q$ is injective and pulls back $\omega$ to a non-degenerate two-form
on $T_zS$ defining the same orientation as $j$.
\end{itemize}
\end{defn}
One easily proves that domain stabilizations always exist as a consequence of the stability assumption.
\begin{lem}
Every object $\alpha$ in ${\mathcal S}^{3,\delta_0}(Q,\omega)$ has a domain stabilization.
\end{lem}
After fixing a domain stabilization we have two stable objects, namely the stable map $\alpha$ in $\text{obj}({\mathcal S}^{3,\delta_0}(Q))$
and the stable noded Riemann surface with unordered marked points $\alpha^\ast=(S,j,M\cup\Xi,D)$, which is an object in ${\mathcal R}$.
The automorphism group $G$ of $\alpha$ is contained in the automorphism group $G^\ast$ of $\alpha^\ast$, and both are finite groups.

We continue with $\alpha^\ast$ and pick a small disk structure ${\bf D}$ which is invariant under $G^\ast$.
Then we fix a good deformation $v\rightarrow j(v)$, with $v\in V$, where $V\subset H^1(\alpha^\ast)$ is an open $G^\ast$-invariant subset.
We can construct the functor  associated to 
$$
(v,\mathfrak{a})\rightarrow (S_\mathfrak{a},j(v)_\mathfrak{a},M^\ast_{\mathfrak{a}}, D_\mathfrak{a}),
$$
where $M^\ast=M\cup\Xi$. Taking a sufficiently small $G^\ast$-invariant  open neighborhood $O$ of $(0,0)$ in $ H^1(\alpha^\ast)\times{\mathbb B}^{\alpha^\ast}$
we obtain a good uniformizer
$$
\Psi^\ast:G^\ast\ltimes O\rightarrow {\mathcal R}.
$$
Denote for $z\in\Xi$ by $[z]$ the $G$-orbit. It follows, that  if $[z]\neq [z']$, then $u(z)\neq u(z')$.  Of course, $u$ takes the same value
on all the elements of $[z]$. For every $[z]$ we fix
a smooth oriented submanifold $H_{[z]}$ in $Q$ of codimension two so that $u(z)\in H_{[z]}$ and
the oriented $\text{range}(Tu(z))\oplus H_{[z]}$ equals the oriented $T_{u(z)}Q$.  Then
$u$ intersects $H_{[z]}$ transversally at $z$. We assume that $H_{[z]}$ does not have a boundary. It is even fine to assume that it is diffeomorphic to an open  ball.
Denote the whole collection ${(H_{[z]})}$ by ${\mathcal H}$. 

Recall  our topological space $\bar{X}^{3,\delta_0}_{(S,j,D),{\bf D}}(Q)$ which contains $u:(S,D)\rightarrow Q$.
From the Sobolev embedding theorem we know that $H^3_{loc}\rightarrow C^1_{loc}$. As a consequence the subset
$\bar{X}^{3,\delta_0}_{(S,j,D),{\bf D},{\mathcal H}}(Q)$, consisting of all $w:(S_\mathfrak{a},D_\mathfrak{a})\rightarrow Q$
which at $z\in\Xi$ intersects $H_{[z]}$ transversally, is well-defined.  We also recall that if $w\in \bar{X}^{3,\delta_0}_{(S,j,D),{\bf D}}(Q)$ 
is close enough to $u$ there is a unique intersection point $z_w$ near $z$ so that $w$ intersects $H_{[z]}$ transversally. 
Moreover $w\rightarrow z_w$ is continuous. One can use this to prove the following before mentioned result.
\begin{prop}
For every map $w$  in  $\bar{X}^{3,\delta_0}_{(S,j,D),{\bf D},{\mathcal H}}(Q)$ there exists an open neighborhood
$U(w)$ in $\bar{X}^{3,\delta_0}_{(S,j,D),{\bf D}}(Q)$ and a continuous map $r:U(w)\rightarrow U(w)$ so that
$r(U(w))=U(w)\cap \bar{X}^{3,\delta_0}_{(S,j,D),{\bf D},{\mathcal H}}(Q)$.
\end{prop}
With other words the subspace associated to the constraints is a local continuous retract.
This part is very important and ties in with the DM-theory.
With our original object being the stable map  $\alpha=(S,j,M,D,u)$ it holds that $(S,D,u)\in \bar{X}^{3,\delta_0}_{(S,j,D),{\bf D},{\mathcal H}}(Q)$.
We also have the object $\alpha^\ast$ in ${\mathcal R}$ and with the small disk structure we obtained the good uniformizing family
$$
(\mathfrak{a},v)\rightarrow (S_\mathfrak{a},j(v)_\mathfrak{a},(M\cup\Xi)_\mathfrak{a},D_\mathfrak{a})
$$
which was used to construct the good uniformizer
$$
\Psi: G^\ast\ltimes O^\ast\rightarrow {\mathcal R}.
$$
Next we form  the topological product space  $V\times \bar{X}^{3,\delta_0}_{(S,j,D),{\bf D},{\mathcal H}}(Q)$
and bring the two separate discussions together.
We define for its elements the map
$$
(v,(S_\mathfrak{a},D_\mathfrak{a},w))\rightarrow (S_\mathfrak{a},j(v)_\mathfrak{a},M_\mathfrak{a},D_\mathfrak{a},w).
$$
The group $G$ acts on  the topological space $V\times \bar{X}^{3,\delta_0}_{(S,j,D),{\bf D},{\mathcal H}}(Q)$ and the above map defines a functor 
$$
{\Psi}: G\ltimes (V\times \bar{X}^{3,\delta_0}_{(S,j,D),{\bf D},{\mathcal H}}(Q))\rightarrow {\mathcal S}^{3,\delta_0}(Q).
$$
The fundamental observation is the following, which follows from the results in \cite{HWZ6}.
\begin{thm}\label{IMP}
For a suitable  $G$-invariant open neighborhood ${O}$ of
$$ 
(0,(S,D,u))\in V\times \bar{X}^{3,\delta_0}_{(S,j,D),{\bf D},{\mathcal H}}(Q)
$$
the following properties hold.
\begin{itemize}
\item[(i)] ${\Psi}:G\ltimes {O}\rightarrow  {\mathcal S}^{3,\delta_0}(Q)$ is full and faithful.
\item[(ii)] The map induced on orbit spaces is a homemomorphism onto an open neighborhood of $|\alpha|$.
\item[(iii)] If $(v,(S_\mathfrak{a},D_\mathfrak{a},w))\in {O}$, then $(\mathfrak{a},v)\in O^\ast$, where we recall the good unifomizer
for ${\mathcal R}$ denoted by $\Psi:G^\ast\ltimes O^\ast\rightarrow {\mathcal R}$.
\item[(iv)] For every $q:=(v,(S_\mathfrak{a},D_\mathfrak{a},w))\in {O}$ there exists an open neighborhood
$U({q})\subset {O}$ so that every sequence $({q}_k)$ in $U({q})$, for which $(|{\Psi}({q}_k)|)$ converges in $|{\mathcal S}^{3,\delta_0}(Q)|$, has a subsequence 
${({q}_{k_\ell})}$ which converges in $\cl_{O}(U({q}))$.
\end{itemize}
\end{thm}
We note that the construction of ${\Psi}$ requires several choices. However starting with an object $\alpha$
we obtain  a set worth ${F}(\alpha)$ of functors having the properties stated in the theorem. 
Hence again we obtain a functorial construction
$$
{F}:({\mathcal S}^{3,\delta_0}(Q,\omega))^-\rightarrow \text{SET}.
$$
\begin{rem}
Moreover, if $\Phi:\alpha\rightarrow \alpha'$
is an isomorphism the possible choices made for $\alpha$ can be bijectively identified with choices for $\alpha'$.
This means $F$ is, in fact, a functor ${\mathcal S}^{3,\delta_0}(Q,\omega)\rightarrow \text{SET}$.
Further there is a precise geometric relationship between $\Psi$ and $\Psi'=F(\Phi)(\Psi)$, similar to the ${\mathcal R}$-case.
\end{rem}
Assume that we are given two such functors ${\Psi}$ and ${\Psi}'$ associated to
stable maps $\alpha$ and $\alpha'$ and suppose 
there exist ${q}_0$ and ${q}_0'$ for which there exists an isomorphism
$$
\Phi_0:{\Psi}({q}_0)\rightarrow {\Psi}'({q}_0').
$$
Then $\Phi_0$ has an underlying biholomorphic map $\phi_0$ and, as we shall see, taking a variation ${q}'$ of ${q}'_0$
there is a uniquely determined core-continuous (see \cite{HWZ6}) germ ${q}'\rightarrow \phi_{q'}$,  and a germ of homeomorphism
${q}\rightarrow {q}({q}')$ so that $\Phi_{{q}'}:=({\Psi}({q}({q'})),\phi_{q'},{\Psi}({q}'))$ satisfies
$$
\Phi_{{q}'}:{\Psi}({q}({q'}))\rightarrow {\Psi}({q'}).
$$
One can use this fact to construct a topology on the transition sets ${\bf M}({\Psi},{\Psi'})$.
For this topology the source and target maps will be local homeomorphisms and all structure maps will be continuous.
One can find the details in \cite{HWZ6}, where this program is carried out in the sc-smooth case, which we shall discuss
later on. The current topological discussion is easier.

In summary, there is  a natural construction $({F},{\bf M})$ for ${\mathcal S}^{3,\delta_0}(Q,\omega)$. Here $F$ associates
to an object $\alpha$ a set ${F}(\alpha)$ of good uniformizers defined on metrizable translation groupoids,
and ${\bf M}$ associates to the transition sets ${\bf M}(\Psi,\Psi')$ metrizable topologies, having the property that the structural maps become continuous
and the source and target maps even local homeomorphisms.

The procedure we just outlined can be extended to deal with 
 $$
 P:{\mathcal E}^{2,\delta_0}(Q,\omega,J)\rightarrow {\mathcal S}^{3,\delta_0}(Q,\omega).
 $$
 In this case the construction $({F},{\bf M})$ can be extend to a construction $({\bar{F}},{\bf \bar{M}})$ covering the construction 
 $({F},{\bf M})$, where ${\bar{F}}:({\mathcal S}^{3,\delta_0}(Q,\omega))^-\rightarrow \text{SET}$, but associates to an object $\alpha$ a good bundle uniformizer
 $$
 \begin{CD}
 G\ltimes K @>\bar{\Psi}>> {\mathcal E}^{2,\delta_0}(Q,\omega,J)\\
 @VpVV    @V P VV\\
 G\ltimes O   @>\Psi>>  {\mathcal E}^{3,\delta_0}(Q,\omega).
 \end{CD}
 $$
 It requires variations of the previous constructions and the necessary details can be found in \cite{HWZ6,H2}.

One might raise the question, given the naturality of the construction, if there is perhaps more structure to be found.
In particular, since perturbations of the  pseudoholomorphic part of ${\mathcal S}^{3,\delta_0}(Q,\omega)$ for a given compatible almost complex structure
are used to define Gromov-Witten invariants. The answer is a resounding `yes'. In fact all the constructions which occured are smooth constructions in an extension of
differential geometry which relies on a more flexible notion of differentiablity in Banach spaces.  In this differential geometry 
there is a much larger library of smooth finite-dimensional or infinite-dimensional local models for smooth spaces.
These local models can have locally varying dimensions but still have tangent spaces. Moreover, there is a nonlinear Fredholm theory
with the usual expected properties.
We shall describe this in the next section.

\section{Smoothness}\label{SER2}
It is an amazing fact that the construction $({F},{\bf M})$ is not just topological, but in fact a smooth construction
within a suitable framework of smoothness, which is quite different from the classical one. It will be described next.
We  keep in mind the occurrence 
of continuous retractions when constructing the domains of the good uniformizers. 

\subsection{Sc-Structures and Sc-Smooth Maps}
Assume we are  given two Banach spaces $E$ and $F$ for which we have as vector spaces a continuous  inclusion $E\subset F$.
In interpolation theory, \cite{Tr}, general methods are developed to construct Banach spaces which interpolate
between $E$ and $F$. We take the concept of a scale (with suitable properties) from interpolation theory, but give it a new interpretation as a generalization of a smooth structure.
This study was initiated in \cite{HWZ1,HWZ2,HWZ3}, and leads to a quite unexpectedly rich theory.  In  \cite{HWZ7,HWZ8} we streamlined the presentation 
and added many further developments, which have not been published before.
\begin{defn}\label{d1}
Let $E$ be a Banach space. An sc-smooth structure (or {sc-structure} for short) for $E$ consists
of a nested sequence of Banach spaces $E_0\supset E_1\supset E_2\supset \cdots$ with $E_0=E$ so that
\begin{itemize}
\item[(1)]  The inclusion  $E_{i+1}\rightarrow E_i$ is  a compact operator.
\item[(2)] $E_\infty=\bigcap_i E_i$ is dense in every $E_m$.
\end{itemize}
\end{defn}

\begin{examp}\label{example32} A typical example is $E=L^2({\mathbb R})$
with the sc-structure given by $E_m:=H^{m,\delta_m}({\mathbb R})$, where $H^{m,\delta_m}({\mathbb R})$ is the Sobolev space
of functions in $L^2$ so that the derivatives up to order $m$ weighted by $e^{\delta_m |s|}$ belong to $L^2$.
Here $\delta_m$ is a strictly increasing sequence starting with $\delta_0=0$.
\end{examp}

If $E$ and $F$ are sc-Banach spaces, then $E\oplus F$ has a natural sc-structure given by
$$
(E\oplus F)_m=E_m\oplus F_m.
$$
Every finite-dimensional vector space has a unique sc-structure, namely the constant one,
where $E_i=E$. If $E$ is infinite-dimensional the constant sequence violates (1) of Definition \ref{d1}.
We continue with some considerations about linear sc-theory.
\begin{defn}
Let $E$ be an sc-Banach space and $F\subset E$ a linear subspace. We call $F$ an { sc-subspace} provided the filtration $F_i=F\cap E_i$
turns $F$ into an sc-Banach space. If $F\subset E$ is an sc-Banach space, then we say that it has an sc-complement, provided there exists an
sc-subspace $G$ such $F_i\oplus G_i=E_i$ as topological linear sum for all $i$.
\end{defn}

Let us note that a finite-dimensional subspace $F$ of $E$ has an sc-complement if and only if $F\subset E_\infty$, see \cite{HWZ1}.

The linear operators of interest are those linear operators $T:E\rightarrow F$,
which map $E_m$ into $F_m$ for all $m$, such that $T:E_m\rightarrow F_m$ is a bounded linear operator.
We call $T$ an { sc-operator}. An { sc-isomorphism} $T:E\rightarrow F$ is a bijective sc-operator so that its inverse is also an sc-operator.
Of particular interest are the linear {sc-Fredholm operators}.

\begin{defn}
A sc-operator $T:E\rightarrow F$ is said to be {sc-Fredholm} provided there exist sc-splittings $E=K\oplus X$ and $F=Y\oplus C$ so that $C$ and $K$ are smooth
and  finite-dimensional, $Y=T(X)$, and $T:X\rightarrow Y$ defines a linear sc-isomorphism.
\end{defn}
We note that the above implies that $E_m=X_m\oplus K$ and $F_m=T(X_m)\oplus C$ for all $m$.  The Fredholm index is by definition
$$
\ind(T)=\dim(K)-\dim(C).
$$
Let us also observe that for every $m$ we have a linear Fredholm operator (in the classical sense)  $T:E_m\rightarrow F_m$, which in particular have the same index and identical kernels.

Next we begin with the preparations to introduce the notion of an sc-smooth map.
If $U\subset  E$ is an open subset we can define an sc-structure for $U$ by the nested sequence
${(U_i)}_{i=0}^\infty$ given by $U_i=E_i\cap U$. We note that $U_\infty=\bigcap U_i$ is dense in every $U_m$.
Considering $U$ with its sc-structure we see that $U_{i_0}\subset E_{i_0}$ also admits an sc-structure defined by
$$
{(U_{i_0})}_m:=U_{i_0+m}.
$$
We write $U^{i_0}$ for $U_{i_0}$ equipped with this sc-structure.
Given two such sc-spaces $U$ and $V$ we write $U\oplus V$ for $U\times V$ equipped with the
obvious sc-structure. Now we can give the rigorous definition of the tangent $TU$ of an open subset $U$ in an sc-Banach space $E$.

\begin{defn}
The { tangent} $TU$ of an  open subset $U\subset E$ of the sc-Banach space $E$ is defined by $TU=U^1\oplus E$.
\end{defn}
We note that
$$
(TU)_i=U_{1+i}\oplus E_i.
$$
Continuing Example \ref{example32} we have
$$
TL^2({\mathbb R}) = H^{1,\delta_1}({\mathbb R})\oplus L^2({\mathbb R})\ \ \hbox{and}\ \
(TL^2({\mathbb R}))_i = H^{i+1,\delta_{i+1}}({\mathbb R})\oplus H^{i,\delta_i}({\mathbb R}).
$$
\begin{defn}
Given two open subsets $U$ and $V$ in sc-Banach spaces, a map $f:U\rightarrow V$ is said to be of class ${ sc}^0$ provided
for every $m$ the map $f$ maps $U_m$ into $V_m$ and the map $f:U_m\rightarrow V_m$ is continuous. 
\end{defn}
The following example takes a little bit of work.
\begin{examp}
Take $L^2({\mathbb R})$ with the previously defined sc-structure and define
$$
\Phi:{\mathbb R}\oplus L^2({\mathbb R}) \rightarrow L^2({\mathbb R}):(t,u)\rightarrow \Phi(t,u),
$$
where $\Phi(t,u)(s)=u(s+t)$. Then $\Phi$ is $\ssc^0$. 
\end{examp}

Next we define the notion of an $\ssc^1$-map.

\begin{defn}
Let $U\subset   E$ and $V\subset F$ be open subsets  in sc-Banach spaces.
An  $\ssc^0$-map $f:U\rightarrow V$ is said to be ${ sc}^1$ provided for every $x\in U_1$ there exists a continuous linear operator $Df(x):E_0\rightarrow F_0$ so that the following holds.
\begin{itemize}
\item[(1)] For $h\in E_1$ with $x+h\in U$ we have
$$
\lim_{\parallel h\parallel_1\rightarrow 0} \frac{1}{\parallel h \parallel_1}\cdot \parallel f(x+h)-f(x)-Df(x)h\parallel_0 =0.
$$
\item[(2)] The map $Tf$ defined by $Tf(x,h)=(f(x),Df(x)h)$ for $(x,h)\in TU$ defines an  $\ssc^0$-map $Tf:TU\rightarrow TV$.
    \end{itemize}
    \end{defn}
Inductively we can define the notion of a  $sc^k$-map. A map is $sc^\infty$ provided it is $sc^k$ for all $k$. The following result shows that the chain-rule holds.
This is quite unexpected since this fact looks not compatible with (1) of the previous definition. However, one is saved by the compactness
of the inclusions stipulated by our definition of sc-structure.
\begin{thm}[Chain-Rule]
Assume that $U,V$ and $W$ are  open subsets in sc-Banach spaces, and $f:U\rightarrow V$ and $g:V\rightarrow W$
are $\ssc^1$-maps. Then $g\circ f$ is $\ssc^1$ and $T(g\circ f)=(Tg)\circ (Tf)$. The same holds for sc$^k$ and sc-smooth maps.
\end{thm}
\begin{examp}
One can show, see \cite{HWZ8.7}, that
the map $\Phi$ from the previous example is sc-smooth. Classically it is nowhere differentiable.
\end{examp}

\subsection{Sc-Smooth
 Spaces and M-Polyfolds}

Now we are in the position to introduce new local models for smooth spaces.
The interesting thing about sc-smoothness is the fact that there are many smooth retractions
with complicated images, so that one obtains a large `library' of smooth local models for spaces.  This library is large enough  to describe problems
occurring when studying the nonlinear Cauchy-Riemann operators in symplectic geometry,  which shows analytical limiting behavior allowing for bubbling-off and similar analytical phenomena. 

\begin{defn}
Let $U\subset E$ be an open subset of the sc-Banach space $E$.
A map $r:U\rightarrow U$ is called a {$\ssc^\infty$-retraction} provided it sc-smooth and $r\circ r=r$.
\end{defn}
The chain rule implies that for a $\ssc^\infty$-retraction $r$ its tangent map $Tr$ is again an $\ssc^\infty$-retraction.
We call the image $O=r(U)$ of an $\ssc^\infty$-retraction $r:U\rightarrow U$ an {$\ssc^\infty$}-{retract}.
The crucial definition is the following.
\begin{defn}
A {local sc-model} is a pair $(O,E)$, where $E$ is an sc-Banach space and $O\subset E$ an  $\ssc^\infty$-retract
given as the image of an  sc-smooth retraction $r:U\rightarrow U$ defined on an open subset  $U$ of $E$.
\end{defn}

The following lemma is easily established.
\begin{lem}
Assume that $(O,E)$ is a local sc-model and $r$ and $s$ are sc-smooth retractions defined on open subsets $U$ and $V$
of $E$, respectively,  having $O$ as the image. Then $Tr(TU)=Ts(TV)$.
\end{lem}

In view of this lemma we can define the tangent of a local sc-model which again is a local sc-model as follows.
\begin{defn}
The { tangent} of the local sc-model $(O,E)$ is defined by
$$
T(O,E):= (TO,TE),
$$
where $TO=Tr(TU)$ for any sc-smooth retraction $r:U\rightarrow U$ having $O$ as the image,
where $U$ is  open in $E$.
\end{defn}
\begin{rem}Let us observe that if $(O,E)$ is a local sc-model and $O'$ an open subset of $O$, then $(O',E)$ is again a local sc-model.
Indeed, if $r:U\rightarrow U$ is an  sc-smooth retraction with $O=r(U)$, then define $U'=r^{-1}(O')$,
 and $r'=r|U':U'\rightarrow U'$ is an sc-smooth retraction with image $O'$.
\end{rem}

A map $f:O\rightarrow O'$ between two local sc-models is { sc-smooth (or $\ssc^k$)} provided $f\circ r:U\rightarrow E'$ is sc-smooth (or $\ssc^k$). One easily verifies that the definition
does not depend on the choice of $r$. We can define the tangent map $Tf:TO\rightarrow TO'$ by
$$
Tf:= T(f\circ r)|Tr(TU).
$$
As it turns out this is well-defined and does not depend on the choice of $r$ as long as it is compatible with $(O,E)$.

\begin{thm}[Chain Rule]
Assume that $(O,E)$, $(O',E')$ and $(O'',E'')$ are local sc-models and $f:O\rightarrow O'$ and $g:O'\rightarrow O''$
are $\ssc^1$. Then $g\circ f:O\rightarrow O''$ is $\ssc^1$ and
$$
T(g\circ f) =(Tg)\circ (Tf).
$$
Moreover if $f,g$ are sc$^k$ so is $g\circ f$. The same for $sc^\infty$.
\end{thm}
The following Remark \ref{2.18} explains how the current account is related to \cite{HWZ1,HWZ2,HWZ3}.
\begin{rem}\label{2.18} In the series of papers \cite{HWZ1,HWZ2,HWZ3} we developed a generalized Fredholm theory in a slightly more restricted situation, which however is more than enough for the applications.  Namely rather than considering $\ssc$-smooth retractions and $\ssc$-smooth retracts, splicings and open subsets of splicing cores were considered, which one can view as a special case.
A splicing consists of an open subset $V$  in some sc-Banach space $W$ and a family of bounded linear projections
$\pi_v:E\rightarrow E$, $v\in V$, where $E$ is another sc-Banach space, so that the map
$$
V\oplus E\rightarrow E:(v,e)\rightarrow\pi_v(e)
$$
is sc-smooth. Then the associated splicing core is $K$, defined by
$$
K=\{(v,e)\in V\oplus E \ |\ \pi_v(e)=e\}.
$$
 Clearly $V\oplus E$ is an open subset  in $W\oplus E$ and $r(v,e):=(v,\pi_v(e))$ is an sc-smooth retraction. The associated retract is, of course, the splicing core $K$. If $O$ is an open subset of $K$ we know that it is again an sc-smooth retract.  Let us note that in all our applications the retractions are obtained from splicings. The above modifications have been implemented in \cite{HWZ7,HWZ8}.
\end{rem}
We demonstrate next how the definition of a manifold can be generalized.
Let $Z$ be a metrizable topological space. A  chart for $Z$ is a tuple $(\varphi,U,(O,E))$, where $\varphi:U\rightarrow O$ is a homeomorphism
and $(O,E)$ is a local sc-model. We say that two such charts are sc-smoothly compatible provided
$$
\psi\circ\varphi^{-1}:\varphi(U\cap V)\rightarrow \psi(U\cap V)
$$
is sc-smooth and similarly for $\varphi\circ \psi^{-1}$. Here $(\psi,V,(P,F))$ is the second chart. Note that the sets $\varphi(U\cap V)$ and $\psi(U\cap V)$
are sc-smooth retracts for sc-smooth retractions defined on open sets in $E$ and $F$, respectively. An sc-smooth atlas for  $Z$
consists of a family of sc-smoothly compatible charts so that their domains cover $Z$. Two sc-smooth atlases are compatible provided their union is an sc-smooth atlas. This defines
an equivalence relation.

\begin{defn}
Let $Z$ be a metrizable space. An  {sc-smooth structure} on $Z$ is given by an sc-smooth atlas. Two sc-smooth structures
are {equivalent} if the union of the two associated atlases defines again an sc-smooth structure.
A M-polyfold  is a metrizable space $Z$ together with an equivalence class of sc-smooth structures.
\end{defn}
We note that these spaces have a natural filtration $Z_0\supset Z_1\supset Z_2\supset \cdots $. The points in $Z_i$ one should view as the points
of some regularity $i$.
The sc-smooth spaces are a very general type of space on which one can define sc-smooth functions.

 It is possible to generalize many of the constructions
from differential geometry to these spaces. If we have an sc-smooth partition of unity we can define Riemannian metrics and consequently a curvature tensor.
Note, however, that curvature would only be defined  at points of regularity at least $2$. The existence of an sc-smooth partition of unity depends on the sc-structure.

The {tangent space} at a point of level at least one is defined in the same way as one defines them for Banach manifolds, see \cite{LA}. Namely one considers
tuples $(z,\varphi,U,(O,E),h)$, where $z\in Z_1$ and $(\varphi,U,(O,E))$ a chart, so that $z\in U$, and $h\in T_{\varphi(z)}O$. Two such tuples, say the second is $(z',\varphi',U',(O',E'),h')$, are said to be equivalent
provided $z=z'$ and $T(\varphi'\circ \varphi^{-1})(\varphi(z))h=h'$.
An equivalence class $[(z,\varphi,U,(O,E),h)]$ then by definition is a tangent vector at $z$. The tangent space at $z\in Z_1$ is denoted by $T_zZ$ and we define
$TZ$ as
$$
TZ=\bigcup_{z\in Z_1} \{z\}\times T_zZ.
$$
One can show that $TZ$ has a natural M-polyfold structure so that the natural map $TZ\rightarrow Z^1$ is sc-smooth, see \cite{HWZ7}.

\begin{examp} Consider the metrizable space $Z$ given as the subspace of ${\mathbb R}^2$ defined by
$$
Z=\{(s,t)\in {\mathbb R}^2\ |\ t=0\ \hbox{if}\ s\leq 0\}.
$$
Then $Z$ admits the structure of an M-polyfold . In order to see this, one constructs a topological embedding into ${\mathbb R}\oplus L^2({\mathbb R})$, where
$L^2({\mathbb R})$ has the previously introduced sc-structure, in such a way that the image is an sc-smooth retract. Here the idea of an sc-smooth splicing comes in handy! Take a smooth, compactly supported map $\beta:{\mathbb R}\rightarrow [0,\infty)$ with $\int \beta(t)^2 ds=1$.
Denote by $f_s$, for $s\in (0,\infty)$ the unit length element in $L^2$ defined by
$$
f_s(t)=\beta(t+e^\frac{1}{s}).
$$
For $s\in (-\infty,0]$ we define $f_s=0$, and denote by $\pi_s$ the $L^2$-orthogonal projection onto the subspace spanned by $f_s$.
Then a somewhat lengthy computation shows that
$$
r:{\mathbb R}\oplus L^2\rightarrow {\mathbb R}\oplus L^2:r(s,u)=(s,\pi_s(u))
$$
is an $\ssc$-smooth retraction, with obvious image $O$ being
$$
\{(s,t\cdot f_s)\ |\ (s,t)\in {\mathbb R}^2\}.
$$
Hence $(O,{\mathbb R}\oplus L^2)$ is a local sc-model.
We note that it has varying dimension.
The map
$$
Z\rightarrow {\mathbb R}\oplus L^2:(s,t)\rightarrow (s, t\cdot f_s)
$$
is a homeomorphic embedding onto $O$. The map is clearly continuous and
injective and has image $O$.  Define ${\mathbb R}\oplus L^2\rightarrow {\mathbb R}^2$ by
$$
(s,x)\rightarrow \left(s,\int_{\mathbb R} x(t)f_s(t)dt\right).
$$
This map is continuous and its restriction to $O$ is the inverse of the previously defined map.
Hence we obtain the structure of an  M-polyfold on $Z$. This gives us the first example of a finite-dimensional space, with varying dimension, which has
a generalized manifold structure. We also note that the induced filtration is constant, so that a tangent space is defined at all points. This is due to the fact that the local model $O$ lies entirely in the smooth part of ${\mathbb R}\oplus L^2$. It is instructive to study sc-smooth curves
$\phi:(-\varepsilon,\varepsilon)\rightarrow O$ satisfying $\phi(0)=(0,0)$. Modifications of the above construction allow us to put smooth structures on the spaces
shown in Figure \ref{porkbarrel}. The M-polyfold does not allow an sc-smooth embedding into any ${\mathbb R}^N$, since then it would have be a smooth manifold by \cite{H.Cartan} . However, as seen in the construction, it can be sc-smoothly embedded into an infinite-dimensional space.
\end{examp}

\subsection{Strong Bundles}
The notion of a strong bundle is designed to give additional structures in the Fredholm theory, which guarantee a compact perturbation and transversality theory.
The crucial point is the fact that there will be a well-defined vector space of perturbations, which have  certain compactness properties.
On the other
hand these perturbations are plentiful enough to allow for different versions of Sard-Smale type theorems, \cite{smale}, in the Fredholm theory.

Let us start with a non-symmetric product $U\triangleleft F$, where $U$ is an open subset in some  sc-Banach space $E$, and $F$ is also an sc-Banach space.
By definition, as a set $U\triangleleft F$ is the product  $U\times F$, but in addition it has a double filtration
$$
(U\triangleleft F)_{m,k}=U_m\oplus F_k
$$
defined for all pairs $(m,k)$ satisfying $0\leq k\leq m+1$. We view $U\triangleleft F\rightarrow U$ as a bundle with base space $U$ and fiber $F$,
where the double filtration has the interpretation that above a point $x\in U$ of regularity $m$ it makes sense to talk about fiber regularity
of a point $(x,h)$ up to order $k$ provided $k\leq m+1$.
At this point it is not clear why one introduces this non-symmetric product coming with a non-symmetric double filtration. We refer the reader to the later Example  \ref{examp} explaining why it is introduced.

Given $U\triangleleft F$, we might consider the associated sc-spaces $U\oplus F$ and $U\oplus F^1$. 
Of interest  for  us are the maps
$$
\Phi:U\triangleleft F\rightarrow V\triangleleft G
$$
of the form
$$
\Phi(u,h)=(\varphi(u),\phi(u,h))
$$
which are linear in $h$. 
\begin{defn}
We say that the map $\Phi$ as described above is of class $\ssc^0_\triangleleft$, provided it induces $\ssc^0$-maps
$U\oplus F^i\rightarrow V\oplus G^i$ for $i=0,1$.
\end{defn}
We define the tangent $T(U\triangleleft F)$ by
$$
T(U\triangleleft F) = (TU)\triangleleft (TF).
$$
Note that the order of the factors is different from the order in $T(U\oplus F)$. One has to keep this in mind. Indeed,
$$
T(U\triangleleft F)=U_1\oplus E\oplus F_1\oplus F\ \ \hbox{and}\ \ T(U\oplus F)=U_1\oplus F_1\oplus E\oplus F.
$$
\begin{defn}
A map $\Phi:U\triangleleft F\rightarrow V\triangleleft G$ is of class $\ssc^1_\triangleleft$ provided
the maps $\Phi:U\oplus F^i\rightarrow V\oplus G^i$ for $i=0,1$ are $\ssc^1$. Taking the tangents of the latter,
gives after rearrangement, the $\ssc^0_\triangleleft$-map
$$
T\Phi:(TU)\triangleleft (TF)\rightarrow (TV)\triangleleft (TG).
$$
Iteratively we can define what it means that a map is $\ssc^k_\triangleleft$ for $k=1,2,\ldots$ and we can also define $\ssc_\triangleleft$-smooth maps.
\end{defn}

Given $U\triangleleft F\rightarrow U$, a sc-smooth section $f$ is map of the form
$x\rightarrow (x,\bar{f}(x))$ such that the induced map $U\rightarrow U\oplus F$ is sc-smooth. In particular, $f$ is `horizontal' with respect to the filtration, i.e. a point on level $m$ is mapped to a point of bi-level $(m,m)$. This can be considered as a convention, and it is precisely this convention which is responsible for the filtration constraint $k\leq m+1$. There is another class of sections called $\ssc^+$-sections. These are sc-smooth sections
of $U\triangleleft F\rightarrow U$ which induce  sc-smooth maps
$U\rightarrow U\oplus F^1$.  In particular, if $s$ is an $\ssc^+$-section of  $U\triangleleft F\rightarrow U$ and  $s(x)=(x,\bar{s}(x))$  for $x\in U_m$  then  $\bar{s}(x)\in F_{m+1}$.
This type of sections will be important for the perturbation theory.
Indeed, it is a kind of compact perturbation theory since the inclusion $F_{m+1}\rightarrow F_m$ is compact.
We give an example before we generalize an earlier discussion about retracts and retractions to bundles of the type  $U\triangleleft F\rightarrow U$.

\begin{examp}\label{examp}
Let us denote by $E$ the Sobolev space $H^1(S^1,{\mathbb R}^n)$ of loops.
We define an sc-structure by $E_m= H^{1+m}(S^1,{\mathbb R}^n)$. Further we define
$F=L^2(S^1,{\mathbb R}^n)=H^0(S^1,{\mathbb R}^n)$ which we filter via
$F_m=H^m(S^1,{\mathbb R}^n)$. Finally we introduce $E\triangleleft F\rightarrow E$. Then we can view the map $f:x\rightarrow\dot{x}$ as an sc-smooth section of
$E\triangleleft F\rightarrow E$.
In particular,  $f$ maps $E_m$ into $E_m\oplus F_m$.
We observe that the filtration of $F$ is picked in such a way that the first order differential operator $x\rightarrow \dot{x}$ is an sc-smooth section,
in particular,  it is horizontal, i.e. the choices are made in such a way that they comply with our convention that sc-smooth sections are index preserving.
The map $x\rightarrow x$ can be viewed as an $\ssc^+$-section. Then $x\rightarrow \dot{x}+x$ is an sc-smooth section obtained from the sc-smooth section $x\rightarrow \dot{x}$ via the perturbation by an $\ssc^+$-section.
Consider now a smooth vector bundle map
$$
\Phi:{\mathbb R}^n\oplus {\mathbb R}^n\rightarrow
{\mathbb R}^n\oplus {\mathbb R}^n
$$
of the form
$$
\Phi(x,h)=(\varphi(x),\phi(x)h),
$$
where $\varphi:{\mathbb R}^n\rightarrow {\mathbb R}^n$ is a diffeomorphism and for every $x\in {\mathbb R}^n$ the map $\phi(x):{\mathbb R}^n\rightarrow {\mathbb R}^n$ is a linear isomorphism. Then we define for $(x,h)\in E\oplus F$ the element $\Phi_\ast(x,h)(t)=(\varphi(x(t)),\phi(x(t))h(t))$.
Note that if $x\in E_m$ and $h\in F_{k}$ for $k\leq m+1$, then $\Phi_\ast(x,h)=:(y,\ell)$ satisfies $y\in E_m$ and $\ell\in F_k$.
However, if $x\in E_m$ and $y\in F_k$ for some $k>m+1$ we cannot conclude that $\ell\in F_k$. We can only say that $\ell\in F_{m+1}$. Now one easily verifies that
$$
\Phi_\ast:E\triangleleft F\rightarrow E\triangleleft F
$$
is $\ssc_\triangleleft$-smooth. This justifies our constraint $k\leq m+1$ for the double filtration. 
\end{examp}

\begin{defn}
An  { $\ssc^\infty_\triangleleft$-retraction} is an  $\ssc_\triangleleft$-smooth map
$$
R:U\triangleleft F\rightarrow U\triangleleft F
$$
with the property $R\circ R=R$.
\end{defn}
Of course, $R$ has the form $R(u,h)=(r(u),\phi(u,h))$ with $r$ being an $\ssc$-smooth retraction and $\phi(u,h)$ linear in the fiber.
Given $R$, we can define its image $K=R(U\triangleleft F)$ and $O=r(U)$. Then we have a natural projection map
$$
p:K\rightarrow O.
$$
We may view this as the local model for a strong bundle. Observe that $K$ has a double filtration and $p$ maps points of regularity
$(m,k)$ to  points of regularity $m$.
\begin{defn}
The tuple $(K,E\triangleleft F)$, where $K$ is a subset of $E\triangleleft F$, so that there exists an
$\ssc^\infty_\triangleleft$-retraction $R$ defined on $U\triangleleft F$, where $U\subset E$ is open  and $K=R(U\triangleleft F)$,
is called a {local strong bundle model}.
\end{defn}
Starting with $(K,E\triangleleft F)$ we have the projection $K\rightarrow E$ and denote its image by $O$
and the induced map by $p:K\rightarrow O$. One can define $T(K,E\triangleleft F)$ by
$$
T(K,E\triangleleft F)=(TK, TE\triangleleft TF),
$$
where $TK$ is the image of $TR$. As before we can show that the definition does not depend on the choice of $R$.

Now we are in the position to define the notion of a strong bundle. Let $p:W\rightarrow X$ be a surjective continuous map between two metrizable spaces,
so that for every $x\in X$ the space $W_x:=p^{-1}(x)$ comes with the structure of a Banach space.  A {strong bundle chart} is a tuple
$(\Phi,p^{-1}(U),E\triangleleft F))$, where $\Phi:p^{-1}(U)\rightarrow K$ is a homeomorphism, covering a homeomorphism $\varphi:U\rightarrow O$,
which between each fiber is a bounded linear operator
$$
\begin{CD}
p^{-1}(U) @>\Phi >> K\\
@V p VV @ VVV\\
U @>\varphi >> O.
\end{CD}
$$
We call two such charts { $\ssc_\triangleleft$-smoothly equivalent} if the associated transition maps
are $\ssc_\triangleleft$-smooth.  We can define the notion of a { strong bundle atlas} and can define the notion of equivalence of two such atlases.
\begin{defn}
Let $p:W\rightarrow X$ be as described before. A { strong bundle structure} for $p$ is given by a strong bundle atlas. Two strong bundle structures are {equivalent}
if the associated atlases are equivalent. Finally $p$ equipped with an equivalence class of strong bundle atlases is called a {strong bundle}.
\end{defn}

Let us observe that a strong bundle $p:W\rightarrow X$ admits a double filtration $W_{m,k}$ with $0\leq k\leq m+1$. By forgetting part of this double filtration we observe that $W(0)$,
which is $W$ filtered by $W(0)_m:=W_{m,m}$, has in a natural way the structure of an  M-polyfold. The same is true for $W(1)$ which is the space $W_{0,1}$ equipped with the filtration $W(1)_m:=W_{m,m+1}$. Obviously the maps $p:W(i)\rightarrow X$ for $i=0,1$ are sc-smooth.

The previously introduced notions of $\ssc$-smooth sections and $\ssc^+$-sections for $U\triangleleft F\rightarrow U$ generalize as follows.
\begin{defn}
Let $p:W\rightarrow X$ be a strong bundle over the M-polyfold $X$ (without boundary).
\begin{itemize}
\item[(i)] A {sc-smooth section} of the  strong bundle $p$ is an  $\ssc^0$-map $s:X\rightarrow W$ with $p\circ s=Id_X$ such that $s:X\rightarrow W(0)$ is sc-smooth.
The vector space of all such sections is written as $\Gamma(p)$.
\item[(ii)] A sc$^+$-section section of the strong bundle $p$  is a sc$^0$-map $s:X\rightarrow W(1)$  with $p\circ s=Id_X$, which in addition is sc-smooth.
The vector space of $\ssc^+$-sections is denoted by $\Gamma^+(p)$. 
\end{itemize}
\end{defn}
In some sense $\ssc^+$-sections are compact perturbations, since the inclusion map
$W(1)\rightarrow W(0)$ is fiber-wise compact.   They are very important for  the perturbation theory.

\subsection{A Special Class of Sc-Smooth Germs}
The next goal   is to define a suitable notion of Fredholm section of a strong bundle. 
The basic fact about the usual Fr\'echet differentiability is the following. If  $f:U\rightarrow F$ is a smooth map (in the usual sense) between an open neighborhood
$U$ of $0\in E$ with target the  Banach space $F$,  and satisfying $f(0)=0$, then we can describe the solution set of $f=0$ near $0$ by an implicit function theorem
provided $df(0)$ is surjective and the kernel of $df(0)$ splits, i.e has a topological linear complement. So smoothness and some properties of the linearized operator at a solution give us always
qualitative knowledge about the solution set near $0$.  On the other hand $f:U\rightarrow F$ being only sc-smooth, $df(0)$ being surjective
and its kernel having an sc-complement is not enough to conclude much about the solution space near $0$. However, as we shall see there is a large class
of sc-smooth maps for which a form of the implicit function theorem holds. In applications the class is large enough to explain gluing constructions (\`a la Taubes and Floer) as   smooth implicit function theorems in the sc-world.

One of the issues which has to be addressed at some point is the fact, that the
spaces we are concerned with, have locally varying dimensions. Though it might sound like a major issue it will turn out that there is a simple way to deal with  these type of problems. It is a crucial observation, that in applications base and fiber dimension change coherently. The sc-formalism incorporates this with a minimum amount of technicalities. One should remark that our presentation is slightly more general than the one given in \cite{HWZ2}.

Let us begin with some notation. As usual $E$ is an sc-Banach space.  We shall write $\mathcal{O}(E,0)$ for an unspecified nested sequence
$U_0\supset U_1\supset U_2\supset \cdots$, where every $U_i$ is an  open neighborhood of $0\in  E_i$. Note that this differs from previous notation where $U_i=E_i\cap U$. When we are dealing with germs we always have the new definition in mind. A {sc-smooth germ}
$$
f:{\mathcal O}(E,0)\rightarrow F
$$
is a map defined on $U_0$  so that for points $x\in U_1$ the tangent map
$Tf:U_1\oplus E_0\rightarrow TF$ is defined, which again is a germ
$$
Tf:{\mathcal O}(TE,0)\rightarrow TF.
$$
We introduce a {basic class} $\mathfrak{C}_{basic}$ of germs of maps as follows.
\begin{defn}
An element in $\mathfrak{C}_{basic}$ is an sc-smooth germ
$$
f:\mathcal{O}({\mathbb R}^{n}\oplus W,0)\rightarrow ({\mathbb R}^N\oplus W,0)
$$
for suitable $n$ and $N$, and an sc-Banach space $W,$  so that the following holds. If $P:{\mathbb R}^N\oplus W\rightarrow W$ is the projection, then $P\circ f$ has the form
$$
P\circ f(r,w)=w-B(r,w)
$$
for $(r,w)\in U_0\subset  {\mathbb R}^{n}\oplus W$. Moreover, for every $\varepsilon>0$ and $m\in {\mathbb N}$
we have
$$
\norm{B(r,w)-B(r,w')}_m\leq \varepsilon\cdot \norm{w-w'}_m
$$
for all $(r,w), (r,w')\in U_m$ close enough to $(0,0)$ on level $m$.
\end{defn}
In \cite{HWZ2} the class of basic germs was slightly more general
in the sense that it was not required that $f(0)=0$ in its definition. However, all important results were then proved under the additional assumption that $f(0)=0$.
In the applications to SFT and the other mentioned theories one can
bring the occurring nonlinear elliptic differential operators
even at bubbling-off points (modulo a filling, which is a crucial concept in the polyfold theory and will be explained shortly) via sc-smooth coordinate changes into the above form, see \cite{HWZ6}, \cite{H2} for Gromov-Witten theory, and  \cite{FHWZ} for the operators in SFT. It is important to note that if $f$ is sc-smooth so that $f(0)=0$ and $Df(0)$ is sc-Fredholm, it is generally not true that after a change of coordinates $f$ can be pushed forward to an element which belongs to $\mathfrak{C}_{basic}$.

As shown in \cite{HWZ2,HWZ7}, basic germs admit something like an infinitesimal
smooth implicit function theorem near $0$ (this is something intrinsic to sc-structures) which for certain maps can be `bound together' to a local implicit function theorem. 

To explain this, assume that $U\subset E$ is an open neighborhood of $0$
and $f:U\rightarrow F$ is an sc-smooth having the following properties, where $U_i=E_i\cap U$.
\begin{itemize}
\item[(i)]  $f(0)=0$ and $Df(0)$ is a surjective sc-Fredholm operator.
\item[(ii)] $f$ is regularizing. This means if $x\in U_m$ and $f(x)\in F_{m+1}$, then $x\in U_{m+1}$.
\item[(iii)]  Viewing $f$ as a section of $U\triangleleft F\rightarrow U$, near every smooth point $x$, and
for a suitable $\ssc^+$-section with $s(x)=f(x)$ the germ 
$$
f-s:{\mathcal O}(E,0)\rightarrow F
$$
is conjugated to a basic germ. 
\end{itemize}
Under these conditions there is a local implicit function theorem near $0$, which guarantees a local solution set
of dimension being the Fredholm index of $Df(0)$ at $0$, and in addition guarantees a natural manifold structure on this solution set.

The infinitesimal implicit function theorem refers to the following phenomena for basic germs. If $f\in \mathfrak{C}_{basic}$, then $Pf(a,w)=w-B(a,w)$,
where $B$ is a family of contractions on every level $m$ near $(0,0)$.
Hence, using Banach's fixed point theorem we find a germ $\delta_m$
solving $\delta_m(a)=B(a,\delta_m(a))$ on level $m$ for $a$ near $0$.
By uniqueness a solution on level $m$ also solves the problem on lower levels. This implies that we have a solution germ $a\rightarrow (a,\delta(a))$ of $Pf(a,w)=0$. The infinitesimal sc-smooth implicit function theorem gives the nontrivial fact that the germ
$$
\delta:{\mathcal O}({\mathbb R}^{n},0)\rightarrow (W,0)
$$
is an sc-smooth germ, see \cite{HWZ2,HWZ7}.

In summary, as we shall discuss in more detail later, if we have a regularizing sc-smooth section which
around every smooth point is conjugated mod a suitable $\ssc^+$-section to a basic germ, then the `infinitesimal' implicit function theorems around points $y$ near $x$, combine together to give a `local' implicit function theorem near a point $x$ where the linearization is surjective. We refer the reader to \cite{HWZ7} for a comprehensive discussion.

\subsection{Sc-Fredholm Sections}
Assume next that $p:K\rightarrow O$ is a strong local bundle, i.e.
$(K,E\triangleleft F)$  is a local strong bundle model.
Suppose $f:{\mathcal O}(O,x)\rightarrow K$ is a germ which we shall write as $[f,x]$.
\begin{defn}
A {filling} for the germ $[f,x]$ consists of the following data.
\begin{itemize}
\item[(1)] An  sc-smooth germ $\bar{f}:{\mathcal O}(E,x)\rightarrow F$.
\item[(2)] A choice of strong bundle retraction $R:U\triangleleft F\rightarrow U\triangleleft F$ such that $K$ is the image of $R$.
\end{itemize}
Viewing $f$ as a map $O\rightarrow F$
such that $\phi(y)f(y)=f(y)$, where $R(y,h)=(r(y),\phi(y)h)$, we assume that the data satisfies the following properties:
\begin{itemize}
\item[(1)] $\bar{f}(y)=f(y)$ for all $y\in O$ near $x$.
\item[(2)] $\bar{f}(y)=\phi(r(y))\bar{f}(y)$ for $y$ near $x$ in $U$ implies that $y\in O$.
\item[(3)] The linearisation of the map
$$
y\rightarrow (Id-\phi(r(y)))\bar{f}(y)
$$
at $x$ restricted to the $\ker(Dr(x))$ defines a linear topological isomorphism $\ker(Dr(x))\rightarrow \ker(\phi(x))$.
\end{itemize}
The germ $[f,x]$ is said to be { fillable} provided there exists
a germ of strong bundle map $\Phi$, covering a germ of (local)  sc-diffeomorphism $\varphi$, so that the push-forward germ $[\Phi_\ast(f),\varphi(x)]$ has a filling. A { filled version} of $[f,x]$ is an sc-smooth germ $[\bar{g},\bar{x}]$ obtained as a filling of a suitable push-forward.
\end{defn}
In the definition the meaning  that $\Phi$ is a germ of strong bundle map covering $\varphi$ is the following. For the given $(K,U\triangleleft F)$, $p:K\rightarrow O$, with $x\in O\subset U$,
there exists $(K',U'\triangleleft F')$, $p':K'\rightarrow O'$, with $x'\in O'\subset E'$ and open neighborhoods $x\in V\subset O$, $x'\in V'\subset O'$,
so that the following is a commutative diagram associated to a strong bundle isomorphism
$$
\begin{CD}
p^{-1}(V) @>\Phi >> {(p')}^{-1}(V')\\
@Vp VV   @V p' VV\\
V @>\varphi >> V'.
\end{CD}
$$
Moreover the size of $V$ and then of $V'$ is unspecified small, but fixed. Also in this case $V_i=O_i\cap V$. Only for the solution germs the neighborhoods
in higher regularity shrink.

If $[f,x]$ has a filling $[\bar{f},x]$, the local study of $f(y)=0$ with $y\in O$ near $x$ is equivalent to the local study of $\bar{f}(y)=0$ where $y\in U$ close to $x$. Let us note that if $f(x)=0$ the linearisation $f'(x):T_xO\rightarrow K_x$ has the same kernel as $\bar{f}'(x):T_xU\rightarrow F_x$ and the cokernels are naturally isomorphic.

\begin{defn}
If $f$ is an sc-smooth section of a strong M-polyfold bundle $p:W\rightarrow X$ with $\partial X=\emptyset$, and $x$ is a smooth point, we say that the germ $[f,x]$ admits a filled version, provided in a suitable  local coordinate representation  $[f,x]$ admits a filled version as defined in the previous definition. We always may assume that the filled version has the form $g:{\mathcal O}(E,0)\rightarrow F$.
\end{defn}
We recall the notion of a regularizing section, which we already mentioned before.
\begin{defn}
Let $p:W\rightarrow X$ be a strong bundle over the M-polyfold $X$ (without boundary) and $f$ an sc-smooth section. We say that $f$ is {regularizing} provided for a point $x\in X$ the assertion $f(x)\in W_{m,m+1}$ implies that $x\in X_{m+1}$.
\end{defn}

Note that for a regularizing section $f$ a solution $x$ of $f(x)=0$ belongs necessarily to $X_\infty$. If $f$ is regularizing and $s\in\Gamma^+(p)$, then $f+s$ is regularizing.
Now we come to the crucial definition.
\begin{defn}
We call the sc-smooth section $f$ of the strong bundle $p:W\rightarrow X$, over an M-polyfold with $\partial X=\emptyset$, an  {sc-Fredholm section}, provided $f$ is regularizing, and around every smooth point $x$ the germ $[f,x]$ has a filled version $[g,0]$ so that for a suitable germ of $\ssc^+$-section $s$ with $s(0)=g(0)$ the germ $[g-s,0]$ is conjugated to an element in $\mathfrak{C}_{basic}$. We denote the collection of all sc-Fredholm sections of $p$ by ${\mathcal F}(p)$.
\end{defn}

\begin{rem}
An sc-Fredholm section according to the above definition is slightly more general than the sc-Fredholm sections defined in \cite{HWZ2}. An additional advantage of the current definition is the stability result
saying,  that
given an sc-Fredholm section for $p:W\rightarrow X$ and an $\ssc^+$-section $s$, then $f+s$ is an sc-Fredholm section for $p:W\rightarrow X$. With the version given in \cite{HWZ2} one can only conclude that $f+s$ is an sc-Fredholm section of $p^1:W^1\rightarrow X^1$. In applications the difference is only `academic'. However, as far as a presentation is concerned this new version is more pleasant. see \cite{HWZ7}.
\end{rem}

The following stability result is crucial for the perturbation theory and rather tautological. In the setup of \cite{HWZ2} it was a nontrivial theorem. However,
some of the burden is now moved to the implicit function theorem, see \cite{HWZ7} for the proofs in this new setup.
\begin{thm}[Stability]
Let $p:W\rightarrow X$ a strong bundle over the M-polyfold $X$.
Then given $f\in {\mathcal F}(p)$ and $s\in \Gamma^+(p)$ we have that
$f+s\in {\mathcal F}(p)$.
\end{thm}

Fredholm sections allow for an implicit function theorem.

\begin{thm}
Assume that $p:W\rightarrow X$ is a strong bundle over the M-polyfold $X$ without boundary. Let
 $f$ be a sc-Fredholm section and $x$ a smooth point such that $f(x)=0$ and $f'(x):T_xX\rightarrow W_x$ is surjective. Then the solution set near $x$ carries in a natural
way the structure of a smooth manifold with dimension being the Fredholm index of $f'(x)$. In addition there exists an open neighborhood $V$ of $x$, so that for every $y\in V$ with $f(y)=0$
the linearisation $f'(y)$ is surjective. Moreover,  its kernel can be identified with the tangent spaces of the solution set at $y$.
\end{thm}

Finally we introduce the notion of an auxiliary norm and give a useful compactness result.
\begin{defn}\label{AUXN}
Assume that $p:W\rightarrow X$ is a strong bundle over the M-polyfold $X$ without boundary. An auxiliary norm for $p$
is a map $N:W\rightarrow {\mathbb R}^+\cup\{+\infty\}$ having the following properties.
\begin{itemize}
\item[(i)] The restriction  $N|W_{0,1}$ is real valued and continuous, and on $W\setminus W_{0,1}$ the map $N$  takes the value $\infty$.
\item[(ii)] $N$ restricted to any fiber of $p:W_{0,1}\rightarrow X$ is a complete norm.
\item[(iii)] If $(h_k)\subset W_{0,1}$ is a sequence such that $p(h_k)\rightarrow x_0$ in $X$, and $N(h_k)\rightarrow 0$, then $h_k\rightarrow 0_{x_0}$ 
in $W_{0,1}$.
\end{itemize}
\end{defn}
The existence of an auxiliary norm can be establishes using continuous partitions of unity.
\begin{prop}
Given a strong bundle $p:W\rightarrow X$ over the M-polyfold $X$ without boundary there exists an auxiliary norm $N$.
For two given auxiliary norms $N_1$ and $N_2$ there exist continuous maps $f_1,f_2:X\rightarrow (0,\infty)$ such that
$$
f_1\cdot N_1 \leq N_2\leq f_2\cdot N_1.
$$
\end{prop}
Now we are in the position to state a useful compactness result for sc-Fredholm sections which is important for the perturbation theory, see 
\cite{HWZ7}.
\begin{thm}
Let $p:W\rightarrow X$ be a strong bundle over the M-polyfold $X$ without boundary. Suppose 
$f$ is an sc-Fredholm section for which $f^{-1}(0)$ is compact. 
\begin{itemize}
\item[(i)] Given an auxiliary norm $N$ there exists an open neighborhood
$U$ of $f^{-1}(0)$ so that for every sc$^+$-section $s\in\Gamma^+(p)$ with support in $U$, satisfying $N(s(x))\leq  1$ for all $x\in X$,
 the set $(f+s)^{-1}(0)$ is compact.
 \item[(ii)] If $X$ admits sc-smooth partitions of unity we find for every $\varepsilon\in (0,1]$ an sc$^+$-section $s$ with support in $U$
 so that $N(s(x))<\varepsilon$ for all $x$, and the set $M=(f+s)^{-1}(0)$ has the structure of a compact smooth manifold without boundary,
 so that the linearization $f'(m):T_mM\rightarrow W_m$ for all $m\in M$ is surjective and the tangent space $T_mM$ can be canonically
 identified with $\ker(f'(m))$.
 \end{itemize}
 \end{thm}

\section{Polyfold Structures and Consequences}
At this point we have discussed a smooth theory  for the category of stable Riemann surfaces based on the construction of good uniformizers,
and a topological theory associated to the category of stable maps. In order  transform the latter discussion into one taking place in a smooth world,
we have generalized the finite- or infinite-dimensional classical differential geometry to a more general sc-smooth differential geometry and described some of the aspects of an associated nonlinear functional analysis. The differential geometric/nonlinear functionalanalytic theory is discussed in great detail in \cite{HWZ7,HWZ8} and gives many more ideas about the framework outlined in Section \ref{SER2}. In the following we show how to use this theory to equip certain categories with sc-smooth structures
and illustrate the ideas with the stable map example.
\subsection{Polyfold Structures on Certain Categories}
 It is useful to start with a general definition  taken from \cite{HWZ8}.\begin{defn}
A good category with metrizable orbit space is given by  a pair $({\mathcal C},{\mathcal T})$,
where ${\mathcal C}$ is a category with the following properties.
\begin{itemize}
\item[(i)] Between any two objects  there are only finitely many morphisms and every morphism is an isomorphism.
\item[(ii)] The orbit space $|{\mathcal C}|$ is a set.
\end{itemize}
Moreover ${\mathcal T}$ is a metrizable topology on $|{\mathcal C}|$. We call $({\mathcal C},{\mathcal T})$ a GCT
(G=good, C=category, T=topology).
\end{defn}
By our discussion in the previous sections the categories ${\mathcal R}$ and ${\mathcal S}^{3,\delta_0}(Q,\omega)$ are GCT's.
Precisely for such categories a construction of type $(F,{\bf M})$ is important and very often naturally exists  if ${\mathcal C}$ is a category coming from geometric considerations. Next we give a precise definition of the type of previously discussed constructions in the polyfold framework.
\begin{defn}
Let ${\mathcal C}$ be  a GCT. A good uniformizer for ${\mathcal C}$ around an object $c$ with automorphism group $G$, written as
$$
\Psi: G\ltimes O\rightarrow {\mathcal C},
$$
consists of a M-polyfold $O$ with an sc-smooth action of $G$ on $O$, where $G\ltimes O$ is the associated translation groupoid,
and $\Psi$ is a functor with the following properties.
\begin{itemize}
\item[(i)] $\Psi$ is full and faithful.
\item[(ii)] There exists a $q_0\in O$ with $\Psi(q_0)=c$.
\item[(iii)] Passing to orbit spaces, $|\Psi|:{_{G}\backslash}O\rightarrow |{\mathcal C}|$ is a homeomorphism onto an open neighborhood
of $|\alpha|$.
\item[(iv)] Given $q\in O$ there exists an open neighborhood $U(q)\subset O$ with the property, that every sequence $(q_k)\subset U(q)$
for which $(|\Psi(q_k)|)$ converges in $|{\mathcal C}|$, has a convergent subsequence with limit in $\cl_O(U(q))$.
\end{itemize}
\end{defn}
As in the continuous case we can define for two good uniformizers the transition set ${\bf M}(\Psi,\Psi')$ and have the usual structural maps.
Here comes the crucial definition, which is the minimalistic version.
Denote for the category ${\mathcal C}$ by ${\mathcal C}^-$ the category having the same
objects, but as morphisms only the identities.
\begin{defn}
Let ${\mathcal C}$ be a GCT. A polyfold structure for ${\mathcal C}$ is a pair $(F,{\bf M})$, where $F:{\mathcal C}^-\rightarrow \text{SET}$ is a functor associating to an object $c$ a set $F(c)$ of good uniformizers, and where ${\bf M}$ associates to two good uniformizers, say $\Psi\in F(c)$ and $\Psi'\in F(c')$, a M-polyfold structure to the transition set ${\bf M}(\Psi,\Psi')$, so that all structural maps are sc-smooth, and the source and target map are local sc-diffeomorphisms.
\end{defn}
Instead of denoting the polyfold structure by  $(F,{\bf M})$ we just write $F$, i.e. 
$F\equiv(F,{\bf M}).$
Let us emphasize that $F(c)$ is a set of good uniformizers, i.e. a specifically picked collection of good uniformizers constructed by a given procedure.
Since classically smooth manifolds are in particular M-polyfolds we see that the constructions associated to ${\mathcal R}$ equip it with a polyfold structure.
\begin{rem}
There are several useful points one should make.\\
\noindent (a) Note that from $F$ the category ${\mathcal C}$ can be recovered.
First of all, since $F$ is defined on ${\mathcal C}^-$, we know all the objects of ${\mathcal C}$.
For two  objects $c$ and $c'$ pick $\Psi\in F(c)$ and $\Psi'\in F(c')$, and let $q\in O$, $q'\in O'$
with $F(q)=c$ and $F(q')=c'$.  Then we can identify the elements $(q,\phi,q')$ in  ${\bf M}(\Psi(q),\Psi'(q'))$
with the morphisms in the  set $\text{mor}_{\mathcal C}(c,c')$ via $(q,\phi,q')\rightarrow \phi$.\\
\noindent (b) We also note that having a polyfold structure on ${\mathcal C}$ there exists a functor 
$$
\text{reg}:{\mathcal C}\rightarrow {\mathbb N}_0\cup\{\infty\}=:{\mathbb N}^\infty_0,
$$
which associates to an object its regularity. Here ${\mathbb N}_0^\infty$ has as morphisms only has the identities.
Given an object $c$ pick $\Psi\in F(c)$ and $q_0\in O$ so that $\Psi(q_0)=c$. Then define
$$
\text{reg}(c)=\sup \{k\in {\mathbb N}_0\ |\ q_0\in O_k\}.
$$
This is well-defined. independent of the choices involved, and a morphism invariant.
We can define the full subcategories ${\mathcal C}_r$ of ${\mathcal C}$ associated to the objects 
$c$ with $\text{reg}(c)\geq r$.\\
\noindent (c)  Given a polyfold structure $F$ for ${\mathcal C}$ we consider ${\mathcal C}_1$
and define $F^1:{\mathcal C}_1\rightarrow \text{SET}$ by
$$
F^1(c)=\{\Psi^1:G\ltimes O^1\rightarrow {\mathcal C}_1\ |\ \Psi\in F(c)\}.
$$
One can equip $|{\mathcal C}_1|$ with a uniquely determined metrizable topology which makes
every element in $F^1(c)$ a good uniformizer.
Then one can show that $F^1$ defines a polyfold structure on ${\mathcal C}_1$. We  shall
denote by ${\mathcal C}^1$ the category ${\mathcal C}_1$ equipped with the polyfold structure.
\end{rem}

We also need a bundle version. As in Subsection \ref{SSS23} denote by $\text{BAN}_G$ the category whose objects are Banach spaces and the morphisms are 
invertible topological linear isomorphisms.
 Suppose we are given a category ${\mathcal C}$ and a functor
$$
\mu:{\mathcal C}\rightarrow \text{BAN}_G.
$$
 Then we can build a new category ${\mathcal E}={\mathcal E}_\mu$ whose objects are pairs $(c,h)$, where $c$ is an object in ${\mathcal C}$ and $h$ is a vector in $\mu(c)$. The morphisms $(c,h)\rightarrow (c',h')$ are the  lifts of the morphisms $\phi:c\rightarrow c'$, namely
$$
(\phi,\mu(\phi),h):(c,h)\rightarrow (c',h')
$$
provided $s(\phi)=c$, $t(\phi)=c'$ and $\mu(\phi)(h)=h'$. We shall abbreviate $\what{\phi}=(\phi,\mu(\phi))$, which is a linear topological isomorphism
$$
\what{\phi}:P^{-1}(s(\phi))\rightarrow P^{-1}(t(\phi)).
$$
The class of all morphisms for the category ${\mathcal E}_\mu$ we shall denote by $\boldsymbol{\mathcal E}_\mu$.
Since  $|{\mathcal C}|$ is a set
the same is true for $|{\mathcal E}|$. We denote by $P=P_\mu:{\mathcal  E}_\mu\rightarrow {\mathcal C }$ the functor which on objects  $(c,h)\rightarrow c$ and
on morphisms $(\phi,\mu(\phi),h)\rightarrow \phi$. We also have the source map $s:\boldsymbol{\mathcal E}_\mu\rightarrow {\mathcal E}_\mu$ 
defined by $s(\phi,\mu(\phi),h)=(s(\phi),h)$ and the target map $t(\phi,\mu(\phi),h)=(t(\phi),\mu(\phi)(h))$.
\begin{defn}
A bundle GCT is given by a tuple $({\mathcal C},\mu,{\mathcal T}_\mu)$, where ${\mathcal C}$ is a GCT,  $\mu:{\mathcal C}\rightarrow \text{BAN}_G$
is a functor, and ${\mathcal T}_\mu$ a metrizable topology on $|{\mathcal E}_\mu|$, so that  $|P_\mu|:|{\mathcal E}_\mu|\rightarrow |{\mathcal C}|$ is continuous
and open. 
\end{defn}
We shall introduce the  notion of strong bundle uniformizers.
The strong bundle uniformizer are build on strong  bundles $p:K\rightarrow O$  over M-polyfolds, equipped
with an action of a  finite group $G$ acting by sc-smooth strong bundle isomorphisms inducing an action of $G$ on $O$ by sc-diffeomorphisms,
so that $p$ is equivariant.
\begin{defn}
A good strong bundle uniformizer for the bundle GCT $({\mathcal C},\mu,{\mathcal T}_\mu)$ around an object $c$ in ${\mathcal C}$
 is a functor
$\bar{\Psi} : G\ltimes K\rightarrow {\mathcal E}$ covering a functor $\Psi:G\ltimes O\rightarrow {\mathcal C}$ so that the following holds.
\begin{itemize}
\item[(i)] The  following diagram is commutative diagram
$$
\begin{CD}
G\ltimes K @>\bar{\Psi} >> {\mathcal E}_\mu\\
@V p VV  @V P_\mu VV\\
G\ltimes O @> \Psi >> {\mathcal C}.
\end{CD}
$$
and $\Psi(q_0)=c$ for some $q_0\in O$.
\item[(ii)] $\bar{\Psi}$ is full and faithful.
\item[(iii)] $|\bar{\Psi}|:|K|\rightarrow |{\mathcal E}_\mu|$ is a homeomorphism onto an open subset of $|{\mathcal E}_\mu|$ of the form
$|P|^{-1}(U)$, where $U=|\Psi(O)|$.
\item[(iv)] The map $\Psi$ is fiber-wise a topological linear Banach space isomorphism.
\item[(v)] Given $q\in O$ there exists an open neighborhood $U(q)\subset O$ so that  every sequence $(h_k)\subset K$,
with $q_k:=p(h_k)\in U(q)$, for which $(|\bar{\Psi}(h_k)|)$ converges in $|{\mathcal E}_\mu|$, has a convergent subsequence in $p^{-1}(\cl_O(U(q)))$.
\end{itemize}
\end{defn}
We observe that automatically $\Psi:G\ltimes O\rightarrow {\mathcal C}$ has to be a good uniformizer.
Given two strong bundle uniformizers $\bar{\Psi}$ and $\bar{\Psi}'$, we can similarly as before define the transition set
$$
{\bf M}(\bar{\Psi},\bar{\Psi}')=\{(h,\bar{\Phi},h')\ |\ h\in K,\ h'\in K',\ \bar{\Phi}\in \text{mor}(\bar{\Psi}(h),\bar{\Psi}'(h'))\}.
$$
Observe that we have a natural map
$$
{\bf M}(\bar{\Psi},\bar{\Psi}')\rightarrow {\bf M}(\Psi,\Psi'):(h,\bar{\Phi},h')\rightarrow (p(h),\Phi,p'(h')).
$$
where a fiber has a a natural Banach space structure.
 We note that $s$ and $t$ in the following diagram are fiber-wise linear (the top row)
 \begin{eqnarray}\label{diagg}
\begin{CD}
K @< s <<  {\bf M}(\bar{\Psi},\bar{\Psi}') @>t >>  K'\\
@V p VV   @VVV @V p' VV\\
O @< s<< {\bf M}(\Psi,\Psi') @> t >> O'
\end{CD}
\end{eqnarray}

The main definition is now that of strong polyfold bundle structure for $({\mathcal C},\mu,{\mathcal T}_\mu)$.
\begin{defn}
A strong polyfold bundle structure for $({\mathcal C},\mu,{\mathcal  T}_\mu)$ is given by $(\bar{F},{\bf M})$
where $\bar{F}:{\mathcal C}^-\rightarrow \text{SET}$ is a functor associating to every object $c$ in ${\mathcal C}$
a set of good strong bundle uniformizers and to every transition 
$$
{\bf M}(\bar{\Psi},\bar{\Psi}')\rightarrow {\bf M}(\Psi,\Psi')
$$
a strong bundle structure so that source and target map define local strong bundle isomorphisms and all structural maps are strong bundle maps, see Diagram
\ref{diagg}.
\end{defn}

\subsection{Tangent Category and Differential Forms}
Assume ${\mathcal C}$ is  a GCT  equipped with a polyfold structure $(F,{\bf M})$. Then we have a filtration 
and inclusion functors, since as a consequence of this construction we can talk about the regularity of an object
$$
{\mathcal C}_\infty ....  \rightarrow {\mathcal C}_{i+1}\rightarrow {\mathcal C }_i...\rightarrow {\mathcal C}_0={\mathcal C}.
$$
We can define the tangent category $T{\mathcal C}$ together with a projection functor $T{\mathcal C}\rightarrow {\mathcal C}^1$.
Here ${\mathcal C}^1$ is the category ${\mathcal C}_1$ with filtration ${\mathcal C}^1_i:={\mathcal C}_{i+1}$ (and its polyfold structure).
Consider tuples $(c,\Psi,(q,h))$, where $c$ is an object in ${\mathcal C}_1$, $\Psi\in F(c)$, say $\Psi:G\ltimes O\rightarrow {\mathcal C}$, $\Psi(q)=c$, and $h\in T_qO$. We shall introduce the notion of equivalence of two such tuples. For this consider a second one, say 
$(c',\Psi',(q',h'))$. Take  suitable open neighborhoods $U(q,1_c,q')\subset {\bf M}(\Psi,\Psi')$, $U(q)\subset O$ and $U(q')\subset O'$, so that the source and target maps
$$
U(q)\xleftarrow{s} U(q,1_c,q')\xrightarrow{t} U(q')
$$
are sc-diffeomorphisms and define 
$$
L:U(q)\rightarrow U(q') : L(p)=(t\circ (s|U(q,1_c,q'))^{-1})(p).
$$
We call  the two tuples equivalent, written as
$$
(c,\Psi,(q,h))\simeq (c',\Psi',(q',h'))
$$
provided the following holds
$$
c=c'\ \text{and}\ \ TL(q)(h)=(q',h').
$$
 We denote an equivalence class by $[(c,\Psi,(q,h))] $ and view them as objects in a category denoted by $T{\mathcal C}$. 
 On the object level we have the projection functor $\tau:T{\mathcal C}\rightarrow {\mathcal C}^1$ given by
 $$
 \tau([c,\Psi,(q,h))]=c.
 $$
 We observe that $\tau^{-1}(c)$ is a Banach space. Given a morphism $\phi:c\rightarrow c'$ we can define a topological linear isomorphism
 $$
 T\phi:\tau^{-1}(c)\rightarrow \tau^{-1}(c')
 $$
 as follows
 \begin{eqnarray}\label{hofer-x}
T\phi ([(c,\Psi,(q,h))]) =[(c',\Psi',TL(q,h)].
\end{eqnarray}
 The morphisms in $T{\mathcal C}$ are given by the tuples
 $$
\Phi:= ([(c,\Psi,(q,h))],T\phi,[(c',\Psi',(q',h'))])
 $$
 where $\phi:c\rightarrow c'$  and 
 $T\phi([(c,\Psi,(q,h))]) =[(c',\Psi',(q',h'))]$. Here 
 $$
 s(\Phi)=[(c,\Psi,(q,h))]\ \text{and}\ t(\Phi)=[(c',\Psi',(q',h'))].
 $$
  The projection functor $\tau$ on the level of morphisms is defined by
 $$
 \tau( ([(c,\Psi,(q,h))],T\phi,[(c',\Psi',(q',h'))]) = \phi.
 $$
 At this point we have shown the following.
 \begin{lem}
Let ${\mathcal C}$ be a   GCT equipped  with  a polyfold structure $(F,{\bf M})$.  Then there is a well-defined tangent category $T{\mathcal C}$ with objects being the equivalence classes
$[(c,\Psi,(q,h))]$, with $c$ being an object in ${\mathcal C}^1$, $(q,h)\in T_qO$, where $\Psi(q)=c$ and 
$\Psi\in F(c)$. The morphisms are given by the tuples
$$
\Phi=([(c,\Psi,(q,h))],T\phi,[(c',\Psi',(q',h'))]),
$$
where $T\phi([(c,\Psi,(q,h))]=[(c',\Psi',(q',h'))]$.
In addition the projection functor $\tau:T{\mathcal C}\rightarrow {\mathcal C}^1$
is defined by $\tau(\Phi)=\phi$ on morphisms and $\tau([(c,\Psi,(q,h))]=c$ on objects.
\end{lem}
We shall show that we can equip $T{\mathcal C}$  with a polyfold structure as well.  
This needs some preparation.
Fix an object $c_0$ in ${\mathcal C}^1$
 and pick $\Psi_0\in F(c_0)$ which is given by
 $$
 \Psi_0: G\ltimes O\rightarrow {\mathcal C}.
 $$
 There is an element $q_0\in O_1$ with $\Psi_0(q_0)=c_0$. For $q\in O_1$ let $c=\Psi_0(q)$ and note that $q$ is an object in ${\mathcal C}_1$.
 Pick $\Psi\in F(c)$, say $\Psi:G'\ltimes O'\rightarrow {\mathcal C}$ and $q'_0\in O'$ with $\Psi(q_0')=c=\Psi_0(q)$. Then $q_0'\in O_1'$.
 Now we us the transition M-polyfold ${\bf M}(\Psi_0,\Psi)$
 and take open neighborhoods $U(q)$, $U(q_0')$ and $U(q,{1_c},q_0')$ so that we have a diagram of sc-diffeomorphisms 
 $$
 U(q)\xleftarrow{s} U(q,1_c,q_0')\xrightarrow{t}U(q_0').
 $$
 We define 
 $$
 L:U(q)\rightarrow U(q_0'):  p\rightarrow L(p):=L=t\circ {(s|U(q,1_c,q'))}^{-1}(p).
 $$ 
 The we define for $(q,h)\in T_qO$, which is an object in $G\ltimes TO$
 $$
 T\Psi_0(q,h) = [(c,\Psi, TL(q,h))]
 $$
 which belongs to $\tau^{-1}(c)$. 
 For a morphism $(g,(q,h)):(q,h)\rightarrow (g\ast q, g\ast h)$ in $G\ltimes TO$ we define the morphism
 $$
 T\Psi_0(g,(q,h)):T\Psi_0(q,h)\rightarrow T\Psi_0(g\ast q,g\ast h)
 $$
 by
\begin{eqnarray*}
 &&T\Psi_0(g,(q,h))\\
 &=& ([(c,\Psi,TL(q,h))],T(\Psi_0(g,q)),T(\Psi_0(g,q))([(c,\Psi,TL(q,h))]))\\
 &=&(T\Psi_0(q,h),T(\Psi_0(g,q)),T(\Psi_0(q,q))(T\Psi_0(q,h))).
 \end{eqnarray*}
 Here $T(\Psi_0(g,q))$ is the tangent associated to the morphism $\Psi_0(g,q):\Psi_0(q)\rightarrow \Psi_0(g\ast q)$ defined as in  (\ref{hofer-x}).
 
We define for an object  $c_0$ in ${\mathcal C}^1$.
 $$
\bar{F}(c_0)=\{ T\Psi: G\ltimes TO\rightarrow T{\mathcal C}\ |\ \Psi\in F(c_0)\},
$$
and view this definition as a functor $\bar{F}:{({\mathcal C}^1)}^-\rightarrow \text{SET}$. 
Then we define 
$$
TF:(T{\mathcal C})^-\rightarrow \text{SET} : TF = \bar{F}\circ \tau,
 $$
i.e. $TF([(c_0,\Psi_0,(q_0,h_0))]):=\bar{F}(c_0)$. Given $T\Psi\in TF(([(c_0,\Psi_0,(q_0,h_0))])$
and $T\Psi'\in ([(c'_0,\Psi'_0,(q_0',h_0'))]$ we can build ${\bf M}(T\Psi,T\Psi')$. 
  The basic result proved in \cite{HWZ8} is the following  theorem.
 \begin{thm}
 Given a polyfold structure $F$ for the ${\mathcal C}$ there exists an associated natural polyfold structure $TF$ for $T{\mathcal C}\rightarrow {\mathcal C}^1$ covering the lifted one for ${\mathcal C}^1$
 given by $F^1$, see the diagram below.
 \end{thm}
  More precisely given $\Psi\in F(c)$, say 
$\Psi:G\ltimes O\rightarrow {\mathcal C}$ the associated $\Psi^1$ fits into the following commutative diagram
$$
\begin{CD}
G\ltimes TO @>T\Psi >> T{\mathcal C}\\
@V \tau_O VV   @V \tau VV \\
O^1 @> \Psi^1 >> {\mathcal C}^1.
\end{CD}
$$
Here we already have an example of an sc-smooth functor. Namely $\tau:T{\mathcal C}\rightarrow {\mathcal C}^1$ has for every $T\Psi$ a  sc-smooth representative, namely $\tau_O$.
One of the main points of having polyfold structures is to say that certain functors defining algebraic 
structures are sc-smooth, or that they are sc-Fredholm functors, which allows a perturbation theory.
Also if ${\mathcal C}$ is equipped with a polyfold structure we can define sc-differential forms
as certain kind of functors.

We can define the $k$-fold product category $T{\mathcal C}\times..\times T{\mathcal C}$ which projects 
to a $k$-fold product of ${\mathcal C}^1$ with itself and pull-back by the multi-diagonal which we denote by
$$
\oplus_{i=1}^k T{\mathcal C}\rightarrow {\mathcal C}^1.
$$
The preimages of objects are $k$-fold products of  Banach spaces.
Viewing ${\mathbb R}$ as a category with only the identities as morphisms we are interested in certain functors
$$
\omega:\oplus_{i=1}^k T{\mathcal C}\rightarrow {\mathbb R},
$$
which are multi-linear and skew-symmetric on the fibers. We assume $T{\mathcal C}$ equipped with its canonical
polyfold structure $TF$. Given $T\Psi\in \bar{F}(c)$ we can pull back $\omega$ via $\oplus_{i=1}^k T\Psi$
and obtain $\omega_\Psi:TO\oplus..\oplus TO\rightarrow {\mathbb R}$. In \cite{HWZ7} we have introduced the notion of an sc-smooth differential form. 
\begin{defn}
The functor $\omega$ is said to be sc-smooth provided there exist a family ${(\Psi_\lambda)}_{\lambda\in\Lambda}$ ($\Lambda$ a set) of good unifomizers associated to $F$ so that the collection of sets $|\Psi_\lambda(O_\lambda)|$ covers $|{\mathcal C}|$, and in addition all the $\omega_{\Psi_\lambda}$ are sc-smooth differential forms.
\end{defn}
Several remarks are in order.
\begin{rem}\noindent(1) One can show that the  definition does not depend on the choice of the family $(\Psi_\lambda)$, see \cite{HWZ7}.\\

\noindent(2) We view $\omega$ as associated to ${\mathcal C}$, despite the fact that it is defined on $T{\mathcal C}$ which lies over  ${\mathcal C}^1$, and therefore
 seemingly only involves ${\mathcal C}^1$,
However, by using the good uniformizers for ${\mathcal C}$ in the definition of $T{\mathcal C}$ incorporates the polyfold structure on ${\mathcal C}$ in a subtle way.\\

\noindent(3) One is tempted to call $\omega$ as defined above a sc-smooth differential form on ${\mathcal C}$. If we call $\omega$ according to
(2) as associated to ${\mathcal C}$ it will turn out that  $d\omega$, the exterior differentiation can only be defined as a form associated to  ${\mathcal C}^1$.
Since we have the system of inclusion functors $...\rightarrow {\mathcal C}^{i+1}\rightarrow {\mathcal C}^i\rightarrow...{\mathcal C}$
one can take a direct limit for forms. The set (see below) of all $[\omega]$ defined by the  direct limit for all $k\geq 0$
turns out to be invariant under $d$ defined by $d[\omega]=[d\omega]$ so that it is better to call $[\omega]$ an sc-smooth differential form.
Using this we obtain a de Rham complex associated to ${\mathcal C}$ as we shall see below. For the moment we shall call
$\omega$ an sc-differential form, and leave the name sc-smooth differential form for the result of a further construction.
\end{rem}
The collection of all sc-smooth functors $\omega$ is a set, since it is completely determined by the set $(\omega_{\Psi_\lambda})$. Using the inclusion functors ${\mathcal C}^{i+1}\rightarrow {\mathcal C}^i$, we can pull-back a functor $\omega$ defined on $\oplus_{i=1}^k T({\mathcal C}^i)$ to
 $\oplus_{i=1}^k T({\mathcal C}^{i+1})$. This is nothing else but restricting $\omega$ to tangent
 vectors of specified higher regularity. Denote by $\Omega^k{(\mathcal C}^i)$ the set of differential
 $k$-forms on ${\mathcal C}^i$, i.e. defined on $\oplus_{j=1}^k T{\mathcal C}^i$. This set has the obvious structure as a real vector space. Then we have the direct system
 $$
\rightarrow  \Omega^k({\mathcal C}^i)\rightarrow \Omega^k({\mathcal C}^{i+1})\rightarrow..
 $$
 and denote the direct limit by $\Omega^k_\infty({\mathcal C})$, and its elements by $[\omega]$.
 If $(X,\beta)$ is an ep-groupoid constructed from the polyfold structure $F$ on ${\mathcal C}$
 we can use the equivalence $\beta:X\rightarrow {\mathcal C}$ to pull-back $[\omega]$
 since $T\beta:TX\rightarrow T{\mathcal C}$ is well-defined. In fact this pull-back completely determines
 $[\omega]$ and the pull back is compatible with the exterior derivative defined on $X$, see \cite{HWZ5},
 \cite{HWZ7}, and \cite{HWZ8}.
 As it turns out the exterior differential is well-defined, so that we obtain the de Rham complex
 $$
 (\Omega_\infty^\ast({\mathcal C},F),d).
 $$
 So in particular there exists a de Rham cohomology.
We refer the reader for details of this theory to \cite{HWZ7,HWZ8} and \cite{HWZ5}.

\subsection{Finite-dimensional, Branched,  Weighted Subpolyfolds}
Suppose that ${\mathcal C}$ is a GCT equipped with a polyfold structure $F$.
We are interested in certain full subcategories which arise 
when studying  Fredholm functors later on. In a first step consider the non-negative rational numbers ${\mathbb Q}^+$ as
objects in a category with the morphisms being the identities. Of interest for us are certain functors 
$$
\Theta:{\mathcal C}\rightarrow {\mathbb Q}^+.
$$
In order to define this class of functors we need the definition of a submanifold $M$ of a M-polyfold
$X$.
\begin{defn}\label{DEX1}
Let $X$ be a M-polyfold and $M$ a subset. We say $M$ is a submanifold provided for every $m\in M$ there exists 
an open neighborhood $U=U(m)\subset X$ and an sc-smooth map $r:U\rightarrow U$ having the following properties.
\begin{itemize}
\item[(i)] $r(U_i)\subset U_{i+1}$ for all $i$ and $r:U\rightarrow U^1$ is sc-smooth.
\item[(ii)] $r\circ r=r$
\item[(iii)] $r(U)=U\cap M$.
\end{itemize}
\end{defn}
First of all we note that $r$ is an sc-smooth retraction, so that $M$ is a sub-polyfold. But the stronger requirement that $r:U\rightarrow U^1$ is sc-smooth,
in fact, implies that the M-polyfold structure on $M$ induced from $ X$ is the sc-smoothly equivalent to the structure of a finite-dimensional smooth manifold, see \cite{HWZ7}. This, of course, justifies
that we call $M$ a submanifold in the first place. More precisely we have the following result.
\begin{prop}
Let $M$ be a submanifold of the M-polyfold $X$ in the sense of Definition \ref{DEX1}.
Then $M$ is a sub-polyfold and its sc-smooth structure induced from $X$ is sc-smoothly equivalent to a classical manifold structure.
\end{prop}
\begin{rem}
A hint that this is true is given by the following consideration. If
$r(x)=x$ it follows that  $x\in X_\infty$, so that for every $m\in M$ the tangent space $T_mM$ is defined
and an sc-Banach space.
Since $Tr(m):T_mX\rightarrow T_mX$ is an sc$^+$-operator its image is compact. Since $Tr(m)$ restricted 
to its image is the identity the image must be finite-dimensional. 
\end{rem}
If $M$, viewed as manifold,  is equipped with an orientation we shall call it an oriented submanifold.

\begin{defn}
Suppose $\Theta:{\mathcal C}\rightarrow {\mathbb Q}^+$ is a functor and $c$ an object. Pick $\Psi\in F(c)$ so that we can take the 
functor $\Theta\circ \Psi: G\ltimes O\rightarrow {\mathbb Q}^+$.  Let $\Psi(q_0)=c$ and assume there exists an open neighborhood $U(q_0)$ in $O$,
finitely many submanifolds ${(M_i)}_{i\in I}$ and positive rational numbers ${(\sigma_i)}_{i\in I}$ such that
$$
\Theta\circ \Psi(q)=\sum_{\{i\in I\ |\ q\in M_i\}} \sigma_i
$$
for all $q\in U(q_0)$. We say $\Theta$ has a smooth finite-dimensional representation with respect to $\Psi$ at $c$. We also say that  the representation is $n$-dimensional provided every $M_i$ is $n$-dimensional.
\end{defn}
If $\Psi'\in F(c)$ with $\Psi'(q_0')=c$  we consider ${\bf M}(\Psi,\Psi')$ and find open neighborhoods so that the source and target map 
are sc-diffeomorphisms 
$$
V(q_0)\xleftarrow{s} V(q_0,1_c,q_0')\xrightarrow{t} V(q_0').
$$
We may assume that $V(q_0)\subset U(q_0)$ and can map the $M_i\cap V(q_0)$ to $M_i'\subset V(q_0')$. 
Then 
$$
\Theta\circ \Psi'(q')=\sum_{\{i\in I\ |\ q\in M_i'\}} \sigma_i
$$
for $q'\in V(q_0')$. Hence if we have a smooth finite-dimensional representation at $c$ for some $\Psi\in F(c)$ it holds for
all uniformizers in $F(c)$.  The same argument goes through if $\phi:c\rightarrow c'$ is an isomorphism and shows if we have a smooth finite-dimensional
representation at $c$ we also have this property at an isomorphic $c'$. Observe that if  the representation at $c$ with respect 
to $\Psi$ is $n$-dimensional this will be true for every isomorphic $c$ as well.
 We can therefore say that $\Lambda$ is smooth at $|c|$ and has a $n$-dimensional representation at $|c|$.
Since the collection of all $|c|$ is a set,  the following makes sense.
\begin{defn}
Let ${\mathcal C}$ be a GCT equipped with a polyfold structure,  and $\Theta:{\mathcal C}\rightarrow {\mathbb Q}^+$ a functor.  We call $\Theta$ a smooth, weighted, branched subpolyfold
of dimension $n$, provided $\Theta$ has a smooth $n$-dimensional representation at every $|c|\in |{\mathcal C}|$.
\end{defn}
There is also a notion of orientation for $\Theta$, see \cite{HWZ5,HWZ8}.
\begin{defn}
Let ${\mathcal C}$ be a polyfold category and $\Theta:{\mathcal C}\rightarrow {\mathbb Q}^+$  a smooth, weighted, branched subpolyfold
of dimension $n$. We say $\Theta$ is closed provided the orbit space associated to all objects $c$ with $\Theta(c)>0$
is a compact subset of $|{\mathcal C}|$.
\end{defn}
An important result is that we can integrate sc-differential forms over closed,  smooth, oriented, weighted branched subpolyfold
of dimension $n$, see \cite{HWZ8,HWZ5}. 
\begin{thm}
Let ${\mathcal C}$ be a GCT equipped with a polyfold structure $F$, and $\Theta:{\mathcal C}\rightarrow {\mathbb Q}^+$ a  closed,  smooth, oriented, weighted, branched subpolyfold
of dimension $n$. Suppose further $[\omega]$ is a n-dimensional sc-smooth differential form on ${\mathcal C}$. Then there is a well-defined integral
$$
\oint_\Theta[\omega],
$$
called the branched integral. 
\end{thm}
The integral is characterized uniquely by certain  properties, and in the somewhat more general context with boundary with corners, even a version of Stokes Theorem holds. The basic observation is that $\Theta$ defines a compact subset of $|{\mathcal C}|$
equipped with,  what is called a Lebesgue $\sigma$-algebra of measurable sets, and
a measurable ${\mathbb Q}^+$-valued  weight function $w$ . The differential form $[\omega]$ induces a signed
measure $\mu_{[\omega]}$ and the integral is given by $\int wd\mu_{[\omega]}$. This construction is formidable and given in \cite{HWZ5} for ep-groupoids, but generalizes immediately to our context, since the ep-groupoid version is compatible with Morita equivalence, i.e. generalized isomorphisms associated to diagrams of sc-smooth equivalences between ep-groupoids, see \cite{HWZ8}.

\subsection{Fredholm Theory and Transversality}
Suppose we are given a bundle GCT $({\mathcal C},\mu,{\mathcal T})$ equipped with a strong polyfold bundle structure $\bar{F}:{\mathcal C}^-\rightarrow\text{SET}$, 
and a  section functor
$$
f:{\mathcal C}\rightarrow {\mathcal E}_\mu.
$$
If $c$ is an object in ${\mathcal C}$ and $\bar{\Psi}$ a strong bundle uniformizer around $c$ we obtain the commutative diagram.
$$
\begin{CD}
G\ltimes K @>\bar{\Psi} >> {\mathcal E}_\mu\\
@V p VV @V P_\mu VV\\
G\ltimes O @>\Psi >> {\mathcal C}
\end{CD}
$$
Since $f$ maps an object to an element in the Banach space associated to $c$, i.e. $P^{-1}(c)$ it follows
that an object $q\in O$ is mapped to $f\circ \Psi(q)\in P^{-1}(\Psi(q))$. However $\bar{\Psi}:p^{-1}(q)\rightarrow P^{-1}(\Psi(q))$ is a linear isomorphism
and it follows that $f$ has a local representative  $f_\Psi$ which is a section of $p$.

\begin{defn}
Let $f$ be a section functor of $P:{\mathcal E}_\mu \rightarrow {\mathcal C}$, where $({\mathcal C},\mu,{\mathcal T})$
is a bundle GCT equipped with a strong polyfold bundle structure $\bar{F}$. We say $f$ is an sc-Fredholm functor
provided there exists a family $(\bar{\Psi}_\lambda)_{\lambda\in\Lambda}$, $\Lambda$ a set, 
of good strong bundle uniformizers, so that  $(\Psi_\lambda)_{\lambda\in\Lambda}$ covers $|{\mathcal C}|$, and 
 every $f_{\bar{\Psi}_\lambda}$ is an sc-Fredholm section of
$K_\lambda\rightarrow O_\lambda$.
\end{defn}
This definition does not depend on the choice of  the strong bundle uniformizers taken from $\bar{F}$.

Of course, like in the classical situation, one is interested in a perturbation theory,
which allows an sc-Fredholm section to be brought into a general position. Since we are dealing with functors, 
there is the added difficulty that symmetry, i.e. compatiblility with morphisms, and transversality are 
competing issues. 

In order to achieve transversality we locally break the symmetry under a small perturbation by an sc$^+$-section,
but keep track of the symmetry by  introducing a finite family (correlated with the initial pertubation) of local sc-Fredholm problems invariant
under the symmetry. In order to have still the right counts of the solutions the problems in the local family need to be weighted.
Of course, these local modifications have to be done coherently, so that overlapping families can be patched together in a suitable way.

In this context it is very  important to understand  how big perturbations can be, in order to guarantee that the perturbed Fredholm section is again proper.
  In order to formulate some results we need
some auxiliary structures. View ${\mathbb R}^+\cup \{+\infty\}=[0,+\infty]$ as a category only having the identities as morphisms.
\begin{defn}
Let $({\mathcal C},\mu,{\mathcal T})$ be a bundle GCT equipped with a strong polyfold bundle structure $\bar{F}:{\mathcal C}^-\rightarrow \text{SET}$.
An {auxiliary norm} $N$ is a functor $N:{\mathcal E}_\mu \rightarrow {\mathbb R}^+\cup\{+\infty\}$ with the following 
properties.  There exists a family of good strong bundle uniformizers ${(\bar{\Psi}_\lambda)}_{\lambda\in \Lambda}$, $\Lambda$ a set, so that
the underlying ${(|\Psi_\lambda(O_\lambda)|)}_{\lambda\in\Lambda}$ cover $|{\mathcal C}|$ and for every $\lambda\in\Lambda$
$$
N\circ \bar{\Psi}_\lambda:K_\lambda\rightarrow {\mathbb R}^+\cup\{+\infty\}
$$
 is an auxiliary norm according to Definition \ref{AUXN}.
\end{defn}
The definition does not depend on the family of good strong bundle uniformizers.
\begin{defn}
Let $({\mathcal C},\mu,{\mathcal T})$ be a bundle GCT equipped with a strong polyfold bundle structure $\bar{F}:{\mathcal C}^-\rightarrow \text{SET}$. A functor $\Lambda:{\mathcal E}_\mu\rightarrow {\mathbb Q}^+$ is called an sc$^+$-smooth multisection functor provided there exists a set ${(\bar{\Psi}_\lambda)}_{\lambda\in A}$ of good strong bundle uniformizers so that the underlying $(\Psi_\lambda)$ have the property
$$
|{\mathcal C}|=\bigcup_{\lambda\in A} |\Psi_\lambda(O_\lambda)|.
$$
with the following property:
\begin{itemize}
\item[($\bullet$)] For every $\Lambda\circ \bar{\Psi}_\lambda: G_\lambda\ltimes K_\lambda\rightarrow {\mathbb Q}^+$ and 
given $q\in O_\lambda$ there exist an open neighborhood $U(q)\subset O_\lambda$, and sc$^+$-sections ${(s_i)}_{i\in I}$ defined for 
$K_\lambda|U(q)$, and rational weights $\sigma_i>0$, $i\in I$, with $\sum_{i\in I}\sigma_i=1$, so that for $e\in K_\lambda|U(q)$
$$
\Lambda\circ \bar{\Psi}_\lambda(e)=\sum_{\{i\in I\ |\ e=s_i(p_\lambda(e))\}}\sigma_i.
$$
\end{itemize}
\end{defn}
\begin{rem}
The definition is independent of the choice of $(\bar{\Psi}_\lambda)$ as long as the associated open sets $(\Psi(O_\lambda))$ cover the orbit space
of ${\mathcal S}^{3,\delta_0}(Q,\omega)$. The functor $\Lambda$ induces a map
$|\Lambda|:|{\mathcal E}_\mu|\rightarrow {\mathbb Q}^+$.
\end{rem}

\begin{defn}
The support of an sc$^+$-smooth multisection functor $\Lambda$ is the full subcategory 
$\text{supp}(\Lambda)$ associated to all objects  $e$ in ${\mathcal E}_\mu$ with  $\Lambda(e)>0$.
The domain support $\text{dom-supp}(\Lambda)$ is the full subcategory associated to the closure of the open subset 
$U$ of $|{\mathcal C}|$ consisting of all $|c|$ so that there exists $|e|\in |P|^{-1}(|c|)$, $|e|\neq 0$, with $|\Lambda|(|e|)>0$.
\end{defn}
We note that each fiber of $|P|$ has a distinguished $0$ element.
Given $\Lambda:{\mathcal E}_\mu\rightarrow {\mathbb Q}^+$ and an auxiliary norm $N:{\mathcal E}_\mu\rightarrow [0,\infty]$
we can measure the size of $\Lambda$. Namely for an object $c$ there are finitely many points $e_i\in P_\mu^{-1}(c)$,
for which $\Lambda(e_i)>0$. Moreover the bi-regularity is at least $(0,1)$, so that $N(e_i)$ is finite. Therefore we can define
$\max\{N(e_i)\ |\ i\in I\}$ and obtain a functor ${\mathcal C}_\mu\rightarrow [0,\infty)$ by
$$
c\rightarrow \max\{N(e_i)\ |\ i\in I\}.
$$
We can pass to orbit space and obtain a continuous map  $N_\Lambda:|{\mathcal C}|\rightarrow [0,\infty)$.

Now we are in the position to describe a series of results in sc-Fredholm theory for our categorical setup.
Without going into much detail there is the notion of orientation for an sc-smooth Fredholm section functor, see \cite{HWZ6,HWZ7,HWZ8}. 
The first result is concerned with a compactness assertion.
\begin{thm}
Let $({\mathcal C},\mu,{\mathcal T})$ be a bundle GCT equipped with a strong polyfold bundle structure $\bar{F}:{\mathcal C}^-\rightarrow \text{SET}$.
Suppose $N:{\mathcal E}_\mu\rightarrow {\mathbb R}^+\cup \{+\infty\}$ is an auxiliary norm,  and $f$ an sc-Fredholm section functor of $P_\mu$, having the property that the orbit space $|f^{-1}(0)|\subset |{\mathcal C}|$ is compact. Then there exists an open neighborhood
$U$ of $|f^{-1}(0)|$ in $|{\mathcal C}|$ so that for every sc$^+$-multisection functor $\Lambda :{\mathcal E}_\mu\rightarrow {\mathbb Q}^+$ with domain support in ${\mathcal C}_U$ and satisfying $N_\Lambda(|c|)\leq 1$ for all $|c|\in |{\mathcal C}|$, the orbit space associated to $\text{supp}(\Lambda\circ f)$ is compact. 
\end{thm}
Note that $\Lambda\circ f :{\mathcal C}\rightarrow {\mathbb Q}^+$ is a functor and $\supp(\Lambda\circ f)$ consists of all objects $c$ with $\Lambda\circ f(c)>0$.
\begin{defn}
If $(U,N)$ is a pair consisting of an auxiliary norm $N$, and an open neighborhood $U$ of the orbit space associated to $f^{-1}(0)$, so that the conclusion of the previous theorem
concerning compactness holds, we shall say that $(U,N)$ controls compactness.
\end{defn}
 In order to construct sc$^+$-multisection functors one needs that the underlying polyfold structure on ${\mathcal C}$
 admits sc-smooth partitions of unity. For example if everything is build on Hilbert scales, or at least the zero-level, then these are available, see \cite{HWZ7,HWZ8}.
 For the following denote by  $H^\ast_{dR}({\mathcal C})$ the de Rham cohomology.
\begin{thm}
Let $({\mathcal C},\mu,{\mathcal T})$ be a bundle GCT equipped with a strong polyfold bundle structure $\bar{F}:{\mathcal C}^-\rightarrow \text{SET}$.
Suppose $N:{\mathcal E}_\mu\rightarrow {\mathbb R}^+\cup \{+\infty\}$ is an auxiliary norm,  and $f$ an sc-Fredholm section functor of $P_\mu$, having the property that the orbit space $|f^{-1}(0)|\subset |{\mathcal C}|$ is compact. We assume that the induced polyfold structure on ${\mathcal C}$ admits sc-smooth partitions of unity. Let $U$ be an open neighborhood around the orbit space associated to $f^{-1}(0)$, so that $(U,N)$ controls compactness.
Then the following holds.
\begin{itemize}
\item[(i)] Given any $\varepsilon\in (0,1)$ there exists a sc$^+$-multisection functor $\Lambda$ with domain support in ${\mathcal C}_U$ and $N_\Lambda(|c|)<\varepsilon$ for all $|c|$, so that 
$$
\Theta:=\Lambda\circ f:{\mathcal C}\rightarrow {\mathbb Q}^+
$$
 is a smooth, closed, weighted, branched subpolyfold of dimension $n$.
\item[(ii)] If $f$ is oriented, i.e. $(f,\mathfrak{o})$,  then $\Theta$ is naturally oriented.
\item[(iv)] Given $[\omega]\in H^n_{dR}({\mathcal C})$ the branched integral $\oint_\Theta[ \omega]$ does not depend  on the choice of $\Lambda$ as long as it is generic and satisfies the above conditions.
\item[(v)] The value of the integral is independent of the choice $(U,N)$ as long as it is admissible.
\end{itemize}
\end{thm}
As a consequence of this theorem the oriented sc-Fredholm section functor $(f,\mathfrak{o})$ defines a linear functional
$$
I_{(f,\mathfrak{o})}:H^\ast_{dR}({\mathcal C})\rightarrow {\mathbb R},
$$
via
$$
I_{(f,\mathfrak{o})}(\boldsymbol{ [ }[\omega]\boldsymbol{ ]})=\oint_{\Lambda\circ f}[\omega],
$$
for $\Lambda$ having support sufficiently small and being generic. This functional will stay the same under even large oriented deformations
$$
t\rightarrow (f_t,\mathfrak{o}_t)
$$
as long as the orbit space of the solution set satisfies some compactness properties. During the deformation we  can even change the $\mu_t$  as long as it is done sc-smoothly for an overall strong polyfold bundle structure. 
There are also appropriate versions where the underlying category changes (sc-smoothly) as well.

\subsection{The Stable Map Example}
In some sense we just need to take a fresh look at what we did in Section \ref{SEC22} and verify, implementing the discussion from Section \ref{SER2},
 that the constructions are sc-smooth. We shall concentrate on the category ${\mathcal S}^{3,\delta_0}(Q,\omega)$ rather than the full problem 
$P:{\mathcal E}^{3,\delta_0}(Q,\omega,J)\rightarrow {\mathcal S}^{3,\delta_0}(Q,\omega)$. In the case of ${\mathcal S}^{3,\delta_0}(Q,\omega)$ our discussion so far gives a precise construction with established topological properties.
From the discussion in \cite{HWZ6}, based on results in \cite{HWZ8.7}, see the forthcoming \cite{H2} for a comprehensive treatment, it follows that indeed all the constructions are sc-smooth. 
\begin{rem}
It is very good exercise using  \cite{HWZ8.7} and some of the results in \cite{HWZ6} to fill in the technical details. 
\end{rem}

We give a few useful  comments. In the constructions mentioned in Theorem \ref{IMP} we already used the exponential gluing profile.
Given the strictly increasing sequence $\delta=(\delta_0,\delta_1,..)$ with all $0<\delta_i<2\pi$, it can be shown that 
$\bar{X}^{3,\delta_0}_{(S,j,D),{\bf D},{\mathcal H}}(Q)$ has a M-polyfold structure. It only depends on $\delta$ and we shall abbreviate 
this topological space equipped with this M-polyfold structure by $\bar{X}^{3,\delta}_{(S,j,D),{\bf D},{\mathcal H}}(Q)$. 
As a consequence  the product 
$$
V\times \bar{X}^{3,\delta}_{(S,j,D),{\bf D},{\mathcal H}}(Q)
$$
 also has a M-polyfold structure.  For this structure the automorphism group
$G$ acts by sc-diffeomorphisms. Taking a  suitable open $G$-invariant neighborhood $O$ of $(0,(S,D,u))$ we obtain
$$
\Psi: G\ltimes O\rightarrow {\mathcal S}^{3,\delta_0}(Q,\omega)
$$
which on objects is given by
$$
(v,(S_\mathfrak{a},D_\mathfrak{a},w))\rightarrow (S_\mathfrak{a},j(v)_\mathfrak{a},M_\mathfrak{a},D_\mathfrak{a},w).
$$
Also recall that we had the underlying good uniformizers for ${\mathcal R}$
$$
G\ltimes O^\ast\rightarrow {\mathcal R},
$$
which on  objects map
$$
(v,\mathfrak{a})\rightarrow (S_{\mathfrak{a}},j(v)_\mathfrak{a},M^\ast_\mathfrak{a},D_\mathfrak{a}).
$$
Our aim is to define a polyfold structure on ${\mathcal S}^{3,\delta_0}(Q,\omega)$ utilizing the previous construction of $F$.
Of course, we need to be able to equip ${\bf M}(\Psi,\Psi')$ with a M-polyfold structure.
In  order to achieve this we have to make an additional requirement which is achieved by possibly restricting 
$\Psi$ to a smaller $O$. This leads to the definition of a good (polyfold) uniformizer.
\begin{defn}\label{UNI}
Assume a sequence $\delta$ and the exponential gluing profile are given as previously described.
Let $\alpha=(S,j,M,D,u)$ be an object in ${\mathcal S}^{3,\delta_0}(Q,\omega)$ with automorphism group $G$. A good polyfold uniformizer $\alpha$ for
${\mathcal S}^{3,\delta_0}(Q,\omega)$ is a functor $\Psi:G\ltimes O\rightarrow {\mathcal S}^{3,\delta_0}(Q,\omega)$
where $O$ is a $G$-invariant open neighborhood of $(0,(S,D,u))$ in $V\times  \bar{X}^{3,\delta}_{(S,j,D),{\bf D},{\mathcal H}}(Q)$ so that the following holds.
\begin{itemize}
\item[(i)] $\Psi$ is fully faithful and $\Psi(v,(S,D,u))=(S,j,M,D,u)$.
\item[(ii)] Passing to orbit spaces, $|\Psi|:{_G\backslash}O\rightarrow |{\mathcal S}^{3,\delta_0}(Q,\omega)|$ is a homeomorphism onto an open neighborhoodof $|\alpha|$.
\item[(iii)] The collection of all $(v,\mathfrak{a})$ occuring in elements $(v,(S_\mathfrak{a},D_\mathfrak{a},w))\in O$ are contained
in $O^\ast$, and $O^\ast\ni (v,\mathfrak{a})\rightarrow (S_\mathfrak{a},j(v)_\mathfrak{a},M_\mathfrak{a}^\ast,D_\mathfrak{a})$ defines a good uniformizer for ${\mathcal R}$.
\item[(iv)] For every object $q\in O$ there exists an open neighborhood $U(q)\subset O$, so that every sequence $(q_k)\subset U(q)$,
for which $|\Psi(q_k)|$ converges in $|{\mathcal S}^{3,\delta_0}(Q,\omega)|$, has a convergent subsequence in $\cl_O(U(q))$.
\end{itemize}
\end{defn}
\begin{rem}
\noindent(a) The condition (iii) is very important and it summarizes four conditions from the definition of a good uniformizer for ${\mathcal R}$. Most important is that the partial Kodaira-Spencer differentials are isomorphisms. This is extensively used when putting a M-polyfold structure on ${\bf M}(\Psi,\Psi')$.
It is crucial for making $s$ and $t$ local sc-diffeomorphisms, see \cite{HWZ6}.\\

\noindent(b) Let $F(\alpha)$ consists of the good uniformizers previously constructed and  their domain $G\ltimes O$ equipped
with the  M-polyfold structures associated to a choice of $\delta$. In addition we assume that the set of all $(v,\mathfrak{a})$ coming from the
$(v,(S_\mathfrak{a},D_\mathfrak{a},w))$ lies in $O^\ast$, so that we have the good uniformizer for ${\mathcal R}$ defined on $G\ltimes O^\ast$.
With $F$ modified as just described it is possible to lift ${\bf M}$ to a functor which associates to $(\Psi,\Psi')$ not only a metrizable space, but in fact a 
M-polyfold structure. The details given in the ep-groupoid seen are presented in \cite{HWZ6}.
\end{rem}

The relevant  construction 
$$
F:{({\mathcal S}^{3,\delta_0}(Q,\omega))}^-\rightarrow \text{SET}
$$
 associates to an object $\alpha$ the set $F(\alpha)$ consisting of good
polyfold uniformizers as described in Definition \ref{UNI}. This leads to the following important result.
\begin{thm}
Given the exponential gluing profile $\varphi$ and an increasing sequence of weights $0<\delta_0<\delta_1<...<2\pi$, the construction $({F},{\bf M})$ given in 
Theorem \ref{IMP}, with the modification just mentioned, defines a polyfold structure, when the local models are equipped with the sc-structures associated to these weights.
\end{thm}

We can also construct a bundle GCT $({\mathcal S}^{3,\delta_0}(Q,\omega),\mu,{\mathcal T})$ as follows. Given an object $\alpha=(S,j,M,D,u)$
we associate to it the Hilbert space $\mu(\alpha)$ consisting of all $(TQ,J)$-valued $(0,1)$-forms $\xi$ along $u$ of class $(2,\delta_0)$.
In particular $\xi(z):(T_zS,j)\rightarrow (T_{u(z)}Q,J)$ is complex anti-linear. A morphism $\Phi:=(\alpha,\phi,\alpha'):\alpha\rightarrow \alpha'$ defines a  linear
topological isomorphism
$$
\mu(\Phi):\mu(\alpha)\rightarrow \mu(\alpha'):  \xi\rightarrow \xi\circ T\phi^{-1}.
$$
This allows us to define the category ${\mathcal E}_\mu$ with objects being the $(\alpha,\xi)$, where $\alpha$ is an object in ${\mathcal S}^{3,\delta_0}(Q,\omega)$,
and $\xi$ is a vector in $\mu(\alpha)$. We shall denote ${\mathcal E}_\mu$ by ${\mathcal E}^{2,\delta_0}(Q,\omega,J)$.
One can define a  metrizable topology ${\mathcal T}$ for $|{\mathcal E}^{2,\delta_0}(Q,\omega,J)|$ and carry out a construction  $(\bar{F},{\bf M})$, 
$$
\bar{F}:({\mathcal S}^{3,\delta_0}(Q,\omega))^-\rightarrow \text{SET}
$$
equipping $P$ with a strong polyfold bundle structure.
Here  $\bar{F}(\alpha)$ consists of good strong bundle uniformizers, where an element $\bar{\Psi}$ in the latter fits into a commutative diagram of functors with certain properties
 $$
 \begin{CD}
 G\ltimes K @>\bar{\Psi}>> {\mathcal E}^{2,\delta_0}(Q,\omega,J)\\
 @V p VV @V P VV\\
 G\ltimes O @>\Psi >> {\mathcal S}^{3,\delta_0}(Q,\omega).
 \end{CD}
 $$
 On objects $\Psi$ has the form 
 $$
 \Psi(v,(S_\mathfrak{a},D_\mathfrak{a},w))=(S_\mathfrak{a},j(v)_\mathfrak{a},M_\mathfrak{a},w)
 $$
 and $\bar{\Psi}$ is given by
 $$
 \bar{\Psi}((v,(S_\mathfrak{a},D_\mathfrak{a},w),\xi)=(S_\mathfrak{a},j(v)_\mathfrak{a},M_\mathfrak{a},w,\xi).
 $$
 The  transition structure
$$
{\bf M}(\bar{\Psi},\bar{\Psi}')\rightarrow {\bf M}(\Psi,\Psi')
$$
has an sc-smooth strong bundle structure. These structures, as already explained, can be viewed as some kind of sc-smooth bundle structure for 
$P:{\mathcal E}^{2,\delta_0}(Q,\omega,J)\rightarrow {\mathcal S}^{3,\delta_0}(Q,\omega)$.  The latter, equipped with this strong bundle structure, which depends on the weight sequence $\delta$, is written as
$$
P:{\mathcal E}^{2,\delta}(Q,\omega,J)\rightarrow {\mathcal S}^{3,\delta}(Q,\omega).
$$

The local representative of the Cauchy-Riemann section takes the form
\begin{eqnarray}\label{display1}
&O\rightarrow K:&\\
&(v,(S_\mathfrak{a},D_\mathfrak{a},w))\rightarrow \left(v,\left(S_\mathfrak{a},D_\mathfrak{a},w\right),\frac{1}{2}\left[ Tw+J(w)\circ Tw\circ j(v)_\mathfrak{a}\right]\right).&\nonumber
\end{eqnarray}
It has been proved in \cite{HWZ6} that the section just defined is sc-Fredholm.
\begin{thm}
The local representative of the section functor $\bar{\partial}_J$ of $P$ with respect to $\bar{\Psi}\in \bar{\bf F}(\alpha)$ for any object 
$\alpha$, as given in (\ref{display1}) is sc-Fredholm. Hence $\bar{\partial}_J$ is an sc-Fredholm section of
$P:{\mathcal E}^{2,\delta}(Q,\omega,J)\rightarrow {\mathcal S}^{3,\delta}(Q,\omega)$.
\end{thm}
In the case of our Gromov-Witten example
consider the full subcategory  ${\mathcal S}_{g,m,A}^{3,\delta}(Q,\omega)$ of ${\mathcal S}^{3,\delta}(Q,\omega)$, where $g\geq 0$ and $m\geq 0$ are integers, and $A\in H_2(Q,{\mathbb Z})$,
consisting of all stable maps of arithmetic genus $g$ with $m$ marked points in the homology class $A$. 
\begin{thm} 
The following holds. 
\begin{itemize}
\item[(i)]  The orbit space $|{\mathcal S}_{g,m,A}^{3,\delta}(Q,\omega)|$ is an open  and closed subset of $|{\mathcal S}^{3,\delta}(Q,\omega)|$
and consequently ${\mathcal S}_{g,m,A}^{3,\delta}(Q,\omega)$ has an induced polyfold structure.  Moreover, 
${\mathcal E}_{g,m,A}^{2,\delta}(Q,\omega,J): ={\mathcal E}^{2,\delta}(Q,\omega,J)|{\mathcal S}_{g,m,A}^{3,\delta}(Q,\omega)$ has an induced strong polyfold bundle structure.
\item[(ii)] The orbit space of $(f|{\mathcal S}_{g,m,A}^{3,\delta}(Q,\omega))^{-1}(0)$ is compact.
\item[(iii)]  For every $1\leq i\leq m$  the evaluation map
\begin{eqnarray}\label{xst1}
ev_i:{\mathcal S}_{g,m,A}^{3,\delta}(Q,\omega)\rightarrow Q
\end{eqnarray}
at the $i$-th marked point is an  sc-smooth functor in the sense that $\text{ev}_i\circ \Psi$ is sc-smooth if $\Psi$ is a good uniformizer.
\item[(iv)] If $2g+m\geq 3$ then the forgetful functor
\begin{eqnarray}\label{xst2}
\sigma:{\mathcal S}_{g,m,A}^{3,\delta}(Q,\omega)\rightarrow {\mathcal R}_{g,m}^{ord}
\end{eqnarray}
into the stable Riemann surface category with ordered mark points is  sc-smooth as well, where we can take good uniformizers for ${\mathcal S}^{3,\delta}_{g,m,A}(Q,\omega)$ and
${\mathcal R}_{g,m}^{ord}$ to see that the local representative is sc-smooth.
\end{itemize}
\end{thm}
The consequence of (iii)  is that the pull-back of a smooth differential form on $Q$ gives an sc-smooth  differential form on ${\mathcal S}_{g,m,A}^{3,\delta}(Q,\omega)$.  For  the  forgetful map  the pull-back of a smooth differential form defines an sc-smooth differential form on ${\mathcal S}_{g,m,A}^{3,\delta}(Q,\omega)$ as a consequence of (iv).

Moreover $\bar{\partial}_J$ restricted to ${\mathcal S}_{g,m,A}^{3,\delta}(Q,\omega)$ is an sc-smooth Fredholm functor with a compact solution set,
and it has a natural orientation, giving us $(\bar{\partial}_J,\mathfrak{o})$ and its restrictions $(\bar{\partial}_{J,(g,k,A)},\mathfrak{o})$.
As a consequence we have the linear maps $I_{(g,m,A)}:=I_{(\bar{\partial}_{J,(g,k,A},\mathfrak{o})}$
$$
I_{(g,m,A)}:H^\ast_{dR}({\mathcal S}_{g,k,A}^{3,\delta}(Q,\omega))\rightarrow {\mathbb R}.
$$
These maps give precisely the data from which one can define the Gromov-Witten invariants, see \cite{HWZ6}.

\begin{rem}  There is, of course, a large literature on Gromov-Witten invariants  and further developments, see \cite{FO,FOOO,LiT,LuT,R1,R2,Tian,McD,MW12,MW14,MWss}. All the methods differ. In some sense, all the approaches had to come up with a fix for the fact that classical Fredholm theory doesn't work.
\end{rem}

The theory which we described here provides a very powerful language to deal with moduli problems in symplectic geometry.
In \cite{FHWZ} we shall use this approach to define a Fredholm setup for SFT. The language for doing so is that developed in \cite{HWZ7,HWZ8}.
The necessary nonlinear analysis comes from \cite{HWZ8.7}.

In \cite{H2} the ideas described here are carried out on the level of a graduate text describing the construction of Gromov-Witten theory and the verification of its axioms, see \cite{Man},  in detail.


\begin{thebibliography}{99}
\bibitem{AR} A. Adem, J. Leida\ and\ Y. Ruan, Orbifolds and Stringy
Topology, Cambridge Tracts in Mathematics, 171. Cambridge University Press, Cambridge, 2007. xii+149 pp.
\bibitem{BSV} V. Borisovich, V. Zvyagin\ and\ V. Sapronov, Nonlinear Fredholm maps and Leray-Schauder degree, Russian
Math. Survey's 32:4 (1977), p 1-54.
\bibitem{BEHWZ} F.~Bourgeois, Y.~Eliashberg, H.~Hofer,
K.~Wysocki and E.~Zehnder,  Compactness Results in Symplectic Field
Theory, {\em Geometry and Topology}, Vol. 7, 2003, pp.799-888.
\bibitem{H.Cartan} H. ~Cartan, Sur les
r\'etractions d'une vari\'et\'e, C. R. Acad.Sc. Paris, t. 303, Serie I, no 14, 1986, p. 715.
\bibitem{CRS} K. Cieliebak, I. Mundet i Riera\ and\ D. A. Salamon,
 Equivariant moduli problems, branched manifolds,
 and the Euler class,
Topology {\bf 42} (2003), no.~3, 641--700.
\bibitem{DK} S. Donaldson\ and\ P. Kronheimer, The geometry of
four-manifolds, Oxford Mathematical Monographs. Oxford Science
Publications. The Clarendon Press, Oxford University Press, New
York, 1990.
\bibitem{EGH} Y. Eliashberg, A. Givental\ and\ H. Hofer, Introduction to Symplectic Field Theory,
 Geom. Funct. Anal. {\bf 2000}, Special Volume, Part II, 560--673.
\bibitem{El} H. Eliasson, Geometry of manifolds of maps, J. Differential Geometry {\bf 1}(1967), 169--194.
\bibitem{FFW} O. Fabert, J. W. Fish, R. Golovko, and K. Wehrheim, Polyfolds: A First and Second Look, arXiv:1210.6670.
\bibitem{H2} J. Fish and  H. Hofer, Lectures on Polyfold Constructions in Symplectic Geometry {I}: The Polyfolds of Gromov-Witten  Theory, in preparation.
\bibitem{FHWZ} J. Fish, H. Hofer, K. Wysocki, and E. Zehnder, Applications of Polyfold Theory {II}:
The Polyfolds of  SFT, in preparation.
\bibitem{FH} A. Floer\ and\ H. Hofer, Coherent orientations for periodic orbit problems in symplectic geometry, Math. Z. {\bf 212} (1993), no.~1, 13--38.
    \bibitem{FO} K. Fukaya\ and\ K. Ono,
    Arnold conjecture and Gromov-Witten invariants.
    Topology,Vol. 38 No 5, 1999.pp. 933-1048.
\bibitem{FOOO} K. Fukaya, Y.-G. Oh, H. Ohta and K. Ono,
Lagrangian intersection Floer theory-anomaly and obstruction, Part I. AMS/IP Studies in Advanced Mathematics, 46.1. American Mathematical Society, Providence, RI; International Press, Somerville, MA, 2009. pp i-xii+396.
\bibitem{FOOO-II} K. Fukaya, Y.-G. Oh, H. Ohta and K. Ono,
Lagrangian intersection Floer theory-anomaly and obstruction, 
Part II. AMS/IP Studies in Advanced Mathematics, 46.2. American Mathematical Society, Providence, RI; International Press, Somerville, MA, 2009. pp. i-xii + 397-805.
\bibitem{GZ} P. Gabriel and M Zisman, Calculus of Fractions and
Homotopy Theory, Ergebnisse Vol. 35, Springer (1967).
\bibitem{G} M.~Gromov, Pseudoholomorphic Curves in
Symplectic Geometry, {\it Inv. Math.} Vol. 82 (1985), 307-347.
\bibitem{hae1} A. Haefliger,  Homotopy and integrability, in Manifolds (Amsterdam, 1970), 133--163, Springer Lecture Notes in Math., 197, 1971.
        \bibitem{hae2} A. Haefliger, Holonomie et classifiants, Asterisque 116 (1984), 70--97.
 \bibitem{hae3} A. Haefliger, Groupoids and foliations, Contemp. Math. 282 (2001), 83--100.
 \bibitem{H1} H. Hofer, A General Fredholm Theory and
Applications, Current Developments in Mathematics, edited by D.
Jerison, B. Mazur, T. Mrowka, W. Schmid, R. Stanley, and S. T. Yau,
International Press, 2006.
\bibitem{H2008} H. Hofer, Polyfolds And A General Fredholm Theory, arxiv 0809.3753.
\bibitem{HWZ-DM} H. Hofer, K. Wysocki\ and\ E. Zehnder,
Deligne-Mumford-Type spaces with a View Towards Symplectic Field
Theory, lecture note in preparation.
\bibitem{HWZ1} H. Hofer, K. Wysocki\ and\ E. Zehnder,
A General Fredholm Theory {I}: A Splicing-Based Differential
Geometry,  J. Eur. Math. Soc. (JEMS)  9  (2007),  no. 4, 841--876.
\bibitem{HWZ2} H. Hofer, K. Wysocki\ and\ E. Zehnder,
A General Fredholm Theory {II}: Implicit Function Theorems,
 Geom. Funct. Anal. 19 (2009), no. 1, 206-293. 
\bibitem{HWZ3} H. Hofer, K. Wysocki\ and\ E. Zehnder,
A General Fredholm Theory {III}: Fredholm Functors and Polyfolds,  Geom. Topol. 13 (2009), no. 4, 2279-2387.
\bibitem{HWZ5} H. Hofer, K. Wysocki\ and\ E.
Zehnder, Integration Theory for  zero sets of polyfold Fredholm
sections, Math. Ann. 346 (2010), no. 1, 139-198.
\bibitem{HWZ6} H. Hofer, K. Wysocki\ and\ E. Zehnder,
Applications of Polyfold Theory {I}: Gromov-Witten Theory, arxiv 1107.2097, Memoirs of the AMS, to appear.
\bibitem{HWZ7} H. Hofer, K. Wysocki\ and\ E. Zehnder, Polyfolds and Fredholm Theory {I}: Basic Theory in M-Polyfolds,  arXiv:1407.3185
\bibitem{HWZ8} H. Hofer, K. Wysocki\ and\ E. Zehnder, Polyfolds and Fredholm Theory {II}: Basic Theory in Polyfolds, in preparation.
\bibitem{HWZ8.7} H. Hofer, K. Wysocki\ and\ E. Zehnder, \textit{Sc-Smoothness, Retractions and New Models for Smooth Spaces}, 
Discrete and Continuous Dynamical Sytems, Vl 28 (No 2), (2010), 665--788.
\bibitem{Joyce} D. Joyce, Kuranishi bordism and Kuranishi homology,
preprint arxiv 0707.3572v1.
\bibitem{Ko} M. Konsevich, Enumeration of rational curves via torus action, in ``The Moduli Space of Curves'' (R. Dijkgraaf, C. Faber and G. van der Geer, eds.) Birkhauser (1995), 335-568.
\bibitem{LA} S. Lang, {\it
Introduction to differentiable manifolds}, Second edition, Springer,
New York, 2002.
\bibitem{LiT} J. Li and G. Tian, Virtual moduli cycles and Gromov-Witten
invariants of general symplectic manifolds, in: Topics in symplectic
4-manifolds (Irvine, CA, 1996), 47-83, Int. Press (1998).
\bibitem{LuT} G. Lu and G. Tian, Constructing virtual Euler cycles and classes,  Int. Math. Res. Surv. IMRS  2007,  2008 in electronic version, Art. ID rym001, 220 pp.
\bibitem{Man} Y. Manin, {\it Frobenius Manifolds, Quantum Cohomology, and
Moduli Spaces}, AMS Colloquium Publications, Volume 47.
\bibitem{Mc} D. McDuff, Groupoids, Branched Manifolds and
Multisection, J. Symplectic Geom. 4, 259-315 (2006).
\bibitem{McD} D. McDuff, Notes on Kuranishi Atlases, arxiv:1411.4306v1.
\bibitem{MW12} D. McDuff and K. Wehrheim, Smooth Kuranishi atlases with trivial isotropy,
arXiv:1208.1340; revision in progress.
\bibitem{MW14} D. McDuff and K. Wehrheim, Smooth Kuranishi atlases with isotropy, work in progress.
\bibitem{MWss} D. McDuff and K. Wehrheim, Stratied smooth Kuranishi atlases, work in progress.
\bibitem{MS} D. McDuff and D. Salamon, {\it Introduction
 to symplectic topology}, 2nd edition, Oxford University Press, 1998.
\bibitem{MS2} D. McDuff and D. Salamon, {\it J-holomorphic curves and
symplectic topology} , Colloquium Publications, vol. 52, Amer. Math.
Soc., Providence, RI, 2004, xii+669 pp.
\bibitem{Mj} I. Moerdijk, Orbifolds as Groupoids: An Introduction, Contemp. Math. 310, 205-222 (2002).
\bibitem{MM} I. Moerdijk and J. Mr\v cun, {\it Introduction to
Foliation and Lie Groupoids}, Cambridge studies in advanced
mathematics, Vol. 91, 2003.
\bibitem{RS} J. Robbin and D. Salamon, A construction of the Deligne-Mumford orbifold. J. Eur. Math. Soc. (JEMS) 8 (2006), no. 4, 611-699. 
\bibitem{R1} Y. Ruan, Topological sigma model and Donaldson type
invariants in Gromov theory, Duke Math. J. 83, 1996, 451-500.
\bibitem{R2} Y. Ruan, Symplectic Topology on Algebraic 3-folds, J.
Diff. Geom. 39, 1994, 215-227.
\bibitem{smale} S. Smale,
An infinite dimensional version of Sard's theorem.
Amer. J. Math. 87 (1965), 861--866.
\bibitem{SZ} S. Suhr and K. Zehmisch, Polyfolds, Cobordisms, and the strong Weinstein conjecture, arXiv:1411.5016.
\bibitem{Tian} G. Tian, Quantum cohomology and its associativity,
Current Developments in Mathematics, 1995, p. 361-397, International
Press
\bibitem{Tr} H. Triebel, {\it Interpolation theory, function spaces, differential operators}, North-Holland, Amsterdam, 1978.


\end{thebibliography}
\end{document}